\documentclass[reqno, 11pt]{amsart}
\pdfoutput=1

\usepackage{eucal}
\usepackage{outlines}
\usepackage{amssymb,amsfonts}
\usepackage{stmaryrd}

\usepackage{autobreak}
\usepackage{enumerate}
\usepackage{mathrsfs}
\usepackage{tikz-cd}
\usepackage{mathtools}
\usepackage{bm}
\usepackage{quiver} 
\usepackage{scalerel}
\usepackage{stackengine,wasysym}
\usepackage{csquotes}
\usepackage{amsmath,calligra,mathrsfs}
\usepackage[style=alphabetic,
minalphanames=1,
maxalphanames=7,
maxbibnames=99]{biblatex}
\usepackage{hyperref}
\usepackage[capitalise]{cleveref}

\newtheorem{thm}{Theorem}[section]

\newtheorem{cor}[thm]{Corollary}
\newtheorem{prop}[thm]{Proposition}
\newtheorem{lem}[thm]{Lemma}

\newtheorem{claim}[thm]{Claim}

\theoremstyle{definition}

\newtheorem{defn}[thm]{Definition}

\newtheorem{exmp}[thm]{Example}

\newtheorem{notn}[thm]{Notation}

\newtheorem{disc}[thm]{Discussion}

\newtheorem{conv}[thm]{Convention}

\newtheorem{interlude}[thm]{Interlude}
\newtheorem{importantrem}[thm]{Important Remark}
\newtheorem*{mainidea*}{Main Idea}
\newtheorem*{disc*}{Discussion}

\theoremstyle{remark}

\newtheorem{rem}[thm]{Remark}

\newtheorem{notation}[thm]{Notation}

\makeatletter
\let\c@equation\c@thm
\makeatother
\numberwithin{equation}{section}

\newcommand{\beqn}{\begin{equation}}
\newcommand{\eeqn}{\end{equation}}
\newcommand{\bclaim}{\begin{claim}}
\newcommand{\eclaim}{\end{claim}}
\newcommand{\blem}{\begin{lem}}
\newcommand{\elem}{\end{lem}}
\newcommand{\bproof}{\begin{proof}}
\newcommand{\eproof}{\end{proof}}
\newcommand{\bdef}{\begin{defn}}
\newcommand{\edefn}{\end{defn}}
\newcommand{\bprop}{\begin{prop}}
\newcommand{\eprop}{\end{prop}}
\newcommand{\bthm}{\begin{thm}}
\newcommand{\ethm}{\end{thm}}
\newcommand{\brem}{\begin{rem}}
\newcommand{\erem}{\end{rem}}
\newcommand{\bcor}{\begin{cor}}
\newcommand{\ecor}{\end{cor}}

\newcommand{\bbA}{\mathbb{A}}

\newcommand{\bbE}{\mathbb{E}}

\newcommand{\bbH}{\mathbb{H}}

\newcommand{\bbN}{\mathbb{N}}

\newcommand{\bbP}{\mathbb{P}}
\newcommand{\bbQ}{\mathbb{Q}}

\newcommand{\bbZ}{\mathbb{Z}}

\newcommand{\calA}{\mathcal{A}}
\newcommand{\calB}{\mathcal{B}}

\newcommand{\calD}{\mathcal{D}}
\newcommand{\calE}{\mathcal{E}}

\newcommand{\calG}{\mathcal{G}}
\newcommand{\calH}{\mathcal{H}}

\newcommand{\calJ}{\mathcal{J}}

\newcommand{\calL}{\mathcal{L}}
\newcommand{\calM}{\mathcal{M}}
\newcommand{\calN}{\mathcal{N}}
\newcommand{\calO}{\mathcal{O}}
\newcommand{\calP}{\mathcal{P}}

\newcommand{\calR}{\mathcal{R}}
\newcommand{\calS}{\mathcal{S}}
\newcommand{\calT}{\mathcal{T}}
\newcommand{\calU}{\mathcal{U}}

\newcommand{\calW}{\mathcal{W}}
\newcommand{\calX}{\mathcal{X}}
\newcommand{\calY}{\mathcal{Y}}

\newcommand{\scrC}{\mathscr{C}}
\newcommand{\scrD}{\mathscr{D}}
\newcommand{\scrE}{\mathscr{E}}

\newcommand{\scrK}{\mathscr{K}}
\newcommand{\scrL}{\mathscr{L}}
\newcommand{\scrM}{\mathscr{M}}
\newcommand{\scrN}{\mathscr{N}}

\newcommand{\sfU}{\mathsf{U}}
\newcommand{\sfV}{\mathsf{V}}

\renewcommand{\a}{\alpha}

\newcommand{\on}{\operatorname}
\newcommand{\wit}{\widetilde}

\newcommand{\un}{\underline}
\newcommand{\wih}{\widehat}

\newcommand{\G}{\mathbb{G}}
\newcommand{\C}{\mathbb{C}}

\newcommand{\Aone}{\mathbb{A}^1}

\newcommand{\hook}{\hookrightarrow}
\newcommand{\F}{\mathcal{F}}
\newcommand{\Four}{\mathscr{F}}

\newcommand{\llb}{\llbracket}
\newcommand{\rrb}{\rrbracket}
\newcommand{\llp}{(\!(}
\newcommand{\rrp}{)\!)}
\mathchardef\mhyphen="2D
\renewcommand{\c}{\mhyphen\mathrm{c}}

\newcommand{\heart}{\heartsuit}

\newcommand{\oblv}{\mathrm{oblv}}

\newcommand{\rh}{\mathrm{rh}}

\newcommand{\an}{\mathrm{an}}

\newcommand{\id}{\mathrm{id}}

\newcommand{\End}{\mathrm{End}}

\newcommand{\Alg}{\mathrm{Alg}}

\newcommand{\QCoh}{\on{\mathsf{QCoh}}}

\newcommand{\Vect}{\mathsf{Vect}}
\newcommand{\Sing}{\on{\mathsf{Sing}}}
\newcommand{\Coh}{\on{\mathsf{Coh}}}
\newcommand{\Perf}{\on{\mathsf{Perf}}}

\newcommand{\Ind}{\on{\mathsf{Ind}}}

\newcommand{\Shv}{\on{Shv}}

\newcommand{\pr}{\on{pr}}

\newcommand{\Fun}{\on{Fun}}

\newcommand{\op}{\on{op}}
\newcommand{\pt}{\on{pt}}

\newcommand{\res}{\on{res}}

\newcommand{\PreShv}{\on{PreShv}}

\newcommand{\coker}{\on{coker}}

\newcommand{\gr}{\on{gr}}
\newcommand{\RH}{\on{RH}}
\newcommand{\DR}{\on{DR}}

\newcommand{\Ab}{\on{Ab}}

\DeclareMathOperator{\colim}{colim}
\DeclareMathOperator{\Perv}{Perv}

\DeclareMathOperator{\Mod}{Mod}

\DeclareMathOperator{\Hom}{Hom}
\DeclareMathOperator{\Sp}{Sp}

\DeclareMathOperator{\sHom}{\mathscr{H}\text{\kern -3pt {\calligra\large om}}\,}

\renewcommand{\Pr}{\calP\mathrm{r}}

\newcommand{\Ch}{\on{Ch}}
\newcommand{\perf}{\mathrm{perf}}

\newcommand{\comp}{\mathrm{comp}}

\newcommand{\Top}{\calT\mathrm{op}}

\newcommand{\Ext}{\on{Ext}}
\newcommand{\Set}{\calS\mathrm{et}}

\renewcommand{\O}{\calO}
\newcommand{\invRH}{\RH^{\mu}}
\newcommand{\wE}{\wih{\scrE}}
\newcommand{\LMod}{\mathrm{LMod}}

\newcommand{\reghol}{\mathrm{reg.\, hol.}}

\newcommand{\Loc}{\calL\mathrm{oc}}
\newcommand{\Conn}{\mathrm{Conn}}

\newcommand{\reg}{\mathrm{reg}}
\renewcommand{\G}{\calG}

\newcommand{\ndalg}{\mathrm{alg}}
\newcommand{\wc}{\mathrm{w.c.}}
\newcommand{\Zar}{\mathrm{Zar}}
\newcommand{\bfAone}{\C}
\newcommand{\fil}{\mathsf{fil}}
\newcommand{\sEnd}{\calE nd}
\newcommand{\Fil}{\mathsf{Fil}}

\newcommand{\Der}{\calD er}
\newcommand{\witsfU}{\wit{\sfU}}
\newcommand{\ev}{\mathrm{ev}}
\newcommand{\mon}{\mathrm{mon.}}
\newcommand{\Ad}{\calA d}
\newcommand{\RFun}{\mathrm{RFun}}
\newcommand{\mult}{\mathrm{mult}}
\newcommand{\Seq}{\mathrm{Seq}}
\newcommand{\Mor}{\mathrm{Mor}}
\renewcommand{\comp}{\mathrm{comp}}
\newcommand{\LM}{\mathcal{LM}}
\newcommand{\Assoc}{\mathrm{Assoc}}

\renewcommand{\PreShv}{\calP\mathrm{rShv}}
\newcommand{\RegHol}{\calR\mathrm{eg}\calH\mathrm{ol}}
\newcommand{\Del}{\mathrm{Del}}
\newcommand{\Temp}{\mathcal{RH}}
\newcommand{\red}{\mathrm{red}}
\newcommand{\sExt}{\calE\mathrm{xt}}
\renewcommand{\sHom}{\calH\mathrm{om}}

\newcommand{\regconn}{\mathrm{reg.\, conn.}}
\newcommand{\lLoc}{\mathrm{Loc}}
\newcommand{\muRH}{\Temp^{\mu}}
\newcommand{\poly}{\,\mathrm{alg}}
\renewcommand{\res}{\mathrm{res}}

\DeclareMathOperator{\wstar}{\wih{\star}}

\DeclareFontFamily{U}{mathx}{\hyphenchar\font45}
\DeclareFontShape{U}{mathx}{m}{n}{
      <5> <6> <7> <8> <9> <10>
      <10.95> <12> <14.4> <17.28> <20.74> <24.88>
      mathx10
      }{}
\DeclareSymbolFont{mathx}{U}{mathx}{m}{n}
\DeclareFontSubstitution{U}{mathx}{m}{n}
\DeclareMathAccent{\widecheck}{0}{mathx}{"71}
\DeclareMathAccent{\wideparen}{0}{mathx}{"75}

\renewcommand{\colim}{\operatornamewithlimits{colim}}

\makeatletter
\newcommand*\bigcdot{\mathpalette\bigcdot@{.5}}
\newcommand*\bigcdot@[2]{\mathbin{\vcenter{\hbox{\scalebox{#2}{$\m@th#1\bullet$}}}}}
\makeatother

\usepackage{geometry}
\title{Derived $V$-filtrations and the Kontsevich--Sabbah--Saito theorem}

\author{Kendric Schefers}

\begin{document}

\begin{abstract}
Let $f: X \to \Aone$ be a regular function on a smooth complex algebraic
variety $X$. We formulate and prove an equivalence between the algebraic formal
twisted de Rham complex of $f$ and the vanishing cycles with respect to $f$
as objects in the category of sheaves valued in the derived $\infty$-category
of modules over $\wE_{\C,0}^{\poly}$, the ring of germs of algebraic
formal microdifferential operators. This is a direct generalization 
of Kontsevich's conjecture, proven in work by Sabbah and then Sabbah--Saito, 
of an algebraic formula computing vanishing cohomology.

The novelty in our approach is the introduction of a canonical $V$-filtration
on the derived $\infty$-category of regular holonomic $\scrD_{\C,0}$-modules,
and the use of various techniques from the theory of higher categories and higher algebra
in the context of the subject of microdifferential calculus.
\end{abstract}

\maketitle

\setcounter{tocdepth}{2}
\tableofcontents



\section{Introduction}
The Kontsevich--Sabbah--Saito theorem is a formula that
computes the vanishing cohomology of a regular function
$f: X \to \Aone$ on a smooth complex variety $X$ in terms
of the algebraic formal twisted de Rham complex.
The formula was originally conjectured by Maxim Kontsevich
in connection with a question appearing in \cite{KoSo}
about how to define the notion of 
vanishing cycles of a regular function
on a smooth formal scheme. Kontsevich's conjecture
was initially proven by Claude Sabbah in \cite{Sabbah},
whose arguments and results were refined and improved upon
in joint work by Sabbah and Morihiko Saito in \cite{SabbahSaito}.

The exact statement of the Kontsevich--Sabbah--Saito theorem may
be found, of course, in \cite{Sabbah} as Theorem 1.1, or later
in this present work as \cref{thm: Sabbah's theorem}.
The formula proven therein is an isomorphism of the $k$\,th
vanishing cohomology of $f$ with the $k$\,th hypercohomology
of the algebraic formal twisted de Rham complex as ordinary
$\C\llp u \rrp$-modules with connection. In other words, it 
determines an isomorphism of ordinary modules over $\wE_{\C,0}^{\poly}$,
the ring of germs of algebraic (i.e. polynomial) formal microdifferential 
operators on $\C$ at $(0;1) \in T^*\C$.\footnotemark
\footnotetext{We provide a description of $\wE_{\C,0}$ and
its subring $\wE_{\C,0}^{\poly}$ later in \cref{sec: wE_{C,0}-modules}.
We indicate how a $\C\llp u \rrp$-module with connection determines
a $\wE_{\C,0}^{\poly}$-module in
\cref{rem: going from u-connections to micromodules}.}
In this work, we formulate and prove a sheaf-theoretic generalization
of the Kontsevich--Sabbah--Saito theorem. Namely, our generalization establishes
an equivalence of Zariski sheaves 
valued in the bounded derived $\infty$-category of $\wE_{\C,0}^{\poly}$-modules.
The following is the precise statement of our main theorem. 

\bthm[{\cref{thm: main theorem}}]
\label{thm: intro main theorem}
Let $X$ be a smooth complex algebraic variety,
and $f: X \to \Aone$ a regular function on $X$.
Let $\scrM$ be a locally free $\O_X$-module of finite rank
equipped with a flat connection $\nabla$ having regular singularity
at infinity, and let $\scrL$ denote its local system of analytic flat sections. 
Let $\tau: X^{\an} \to X^{\Zar}$ denote the canonical
continuous function from the analytic topology on $X$ to the underlying
Zariski topology on $X$.
Then there is an equivalence of objects in
$\Shv^b(X^{\Zar}; \wE_{\C,0}^{\poly})$,
\beqn
\label{eqn: intro main theorem}
\DR(\wE_X^{-f/u} \otimes_{\O_X} \scrM) 
	\simeq \bigoplus_{c \in \C} \un{\wih{\scrE}}^{-c/u} \otimes_{\C\llp u \rrp} 
		\tau_*\muRH_{\Shv}\left(\varphi_{f-c}(\scrL), T_{f-c} \right),
\eeqn
which recovers the original statement
of the Kontsevich--Sabbah--Saito theorem
upon taking hypercohomology.
\ethm

The object appearing on the left-hand side of
\labelcref{eqn: intro main theorem} is the algebraic
formal twisted de Rham complex on $X$, while $\muRH_{\Shv}$
is a functor from the $\infty$-category of constructible $\calD^b(\Loc(\C \setminus 0))$-valued 
sheaves to the $\infty$-category of constructible\footnotemark 
\footnotetext{Constructibility for us will always include the condition that
stalks belong to some nice subcategory. See \cref{sec: Appendix A}
for details. In this case, the stalks of a constructible $\calD^b(\wE_{\C,0})$-valued
sheaf are required have regular holonomic cohomology $\wE_{\C,0}$-modules.}
$\calD^b(\wE_{\C,0})$-valued sheaves.
Precise definitions of all these terms are given later in the text.

A novelty in our approach is the introduction of a derived
$V$-filtration (also called the Kashiwara-Malgrange filtration), 
on objects of $\calD^b_{\reghol}(\scrD_{\C,0})$,
the bounded derived $\infty$-category of regular holonomic
modules over $\scrD_{\C,0}$, the ring of germs of differential operators at $0 \in \C$.
This novelty, combined with a suitable interpretation of objects appearing
in the Kontsevich--Sabbah--Saito theorem, allow us to bootstrap Sabbah's
original arguments in \cite{Sabbah} to the setting of sheaves
valued in higher categories.

\subsection{A confusing point}
The original Kontsevich--Sabbah--Saito theorem is
formulated as an isomorphism of various $\C\llp u \rrp$-modules
with connection, but the proof found in \cite{Sabbah} 
uses the formalism of $\wE_{\C,0}$-modules. 
At various points throughout this work, we work with
$\wE_{\C,0}$-modules, $\wE_{\C,0}^{\poly}$-modules,
and $\C\llp u \rrp$-modules with connection. We clarify 
here the relationship between these three types of objects, and
explain the sometimes subtle passage from one object to another.

A $\C\llp u \rrp$-module with connection is a module $\C\llp u \rrp$-module
$E$ equipped an operator $\nabla: \C\llp u \rrp\partial_u \times E \to E$
which is $\C\llp u \rrp$-linear in the first variable and which satisfies the Leibniz
rule with respect to the $\C\llp u \rrp$-action in the second.
As explained in \cref{rem: going from u-connections to micromodules},
the connection on $E$ allows us to define an action of $\wE_{\C,0}^{\poly}$,
that is, Laurent series in the formal variable $u$ whose coefficients are polynomials
in the affine coordinate $t$. Conversely, given a $\wE_{\C,0}^{\poly}$-module $M$,
its underlying $\C\llp \xi \rrp$-module $E$ is a $\C\llp u \rrp$-module with a canonical connection
after relabeling $u := \xi$ and setting $\nabla_{\partial_u}(e) := - t \cdot e\xi^{-2}$
for $e \in E$. By means of the resultant equivalence of categories, 
$\Conn_{\C\llp u \rrp} \simeq \calA_{\wE_{\C,0}^{\poly}}$, we use
$\wE_{\C,0}^{\poly}$-modules and $\C\llp u \rrp$-modules with connection
interchangeably throughout the text.

As its notation suggests, $\wE_{\C,0}^{\poly}$ is a subring of
$\wE_{\C,0}$, the ring of germs at $(0;1) \in T^*\C$ of formal
microdifferential operators on $\C$. The abelian categories
of modules over these two rings are not equivalent, but because
we would like to use existing results in the literature, it will be
convenient to work with $\wE_{\C,0}$-modules at times. This
is facilitated by the fact that the categories of \emph{regular holonomic}
modules over $\wE_{\C,0}$ and $\wE_{\C,0}^{\poly}$, respectively,
\emph{are} equivalent (see \cref{lem: poly}). Because the central
objects of our proof are constructible sheaves with regular holonomic
stalks, this equivalence provides enough leverage to put the theory of 
$\wE_{\C,0}$-modules to use in proving the main theorem.

\subsection{Structure of the paper}
In \cref{sec: Notation and Conventions}, we list a few
conventions that hold throughout the body of this work,
and we include a glossary of commonly used categories 
with possibly nonstandard notation, for the benefit
of the reader.

In \cref{sec: Preliminaries on D_X-modules}, we recall some
preliminaries definitions and notions in the theory of analytic
$\scrD$-modules, including the $V$-filtration on holonomic modules
and the statement of the Riemann--Hilbert correspondence.
We end the section with a particular formulation of Deligne's correspondence
that we use throughout the body of the paper.

In \cref{sec: D_{C,0}-modules} we establish the basic
notion and facts regarding modules over the ring of germs
at $0$ of differential operators on $\C$.

In \cref{sec: The V-filtration on D_{C,0}-modules}, we introduce
a $V$-filtration functor on $\calD^b_{\reghol}(\scrD_{\C,0})$
taking values in the category of $V^{\bullet}\scrD_{\C,0}$-module
objects in the category of filtered spectra.

In \cref{sec: wE_{C,0}-modules}, we recall the ring of
germs of formal microdifferential operators at $(0;1)
\in T^*\C$ and establish an equivalence between
the formal microlocalization and vanishing cycles
in $\calD^b_{\reghol}(\scrD_{\C,0})$.

In \cref{sec: Sheaves of D_{C,0}-modules}, we show
that the constructions of the previous sections
induce corresponding constructions on the $\infty$-categories
of constructible sheaves valued of $\scrD_{\C,0}$-modules
and $\wE_{\C,0}$-modules. In particular, we show the
existence and uniqueness of $V$-filtration on 
objects of $\Shv_c(X;\scrD_{\C,0})$.

In \cref{sec: The Main Theorem}, we recall the precise
statement of the Kontsevich--Sabbah--Saito theorem, and formulate
our main theorem. We also state the proof of the main theorem
using results from the following section.

In \cref{sec: Proof of the Main Theorem}, we state and prove 
the components used in the proof of the Main Theorem of
the previous section. Our results and presentation in this section
follow Sabbah's in \cite{Sabbah} very closely. The only originality we
contribute here is the application of the foundational results of the
previous sections that allows us to adapt Sabbah's arguments to
the setting of sheaves in higher categories.

In \cref{sec: Appendix A}, we establish our notion
and conventions for constructible sheaves taking values
in a wide variety of $\infty$-categories, as well as prove
some results needed to prove our Main Theorem and
show the existence of a $V$-filtration on sheaves. This
section is a mix of known and (as far as we know)
new results, and we have tried to indicate 
throughout which is which.

In \cref{sec: Appendix B}, we collect some results on
derived local systems on the punctured complex
plane for use in stating a derived version of Deligne's
correspondence.

In \cref{sec: Appendix C}, we establish our notation
and conventions for the category of filtered objects
in an $\infty$-category admitting sequential limits,
which are taken directly from the paper of
Gwilliam and Pavlov on the subject. We also
prove a few results on filtered objects in the
context of categories of sheaves.

\subsection{Acknowledgments} 
I would like to thank Sam
Raskin for his help with a number of technical
details in this work, as well as
David Ben-Zvi for suggesting a major 
simplification of my original approach to this problem.

Special thanks are due to Pavel Safronov for reading and 
pointing out several errors in an initial draft of this paper.

I would also like to thank Tom Gannon, Rok Gregoric,
and Sam Gunningham for helpful discussions relating
to this work.

I am very grateful to Mads Bach Villadsen for patiently
explaining aspects of $V$-filtrations to me
over several days during the Summer School
on Motivic Integration theory hosted at the
HHU D\"{u}sseldorf, where he and I met. 

Lastly, I warmly thanks Claude
Sabbah for answering questions about his own
work, and Pierre Schapira for
answering several questions about
microdifferential operators.

\section{Notation and Conventions}
\label{sec: Notation and Conventions}
\subsubsection{}
Derived functors are denoted using the same notation as their underived
counterparts, unless otherwise written.

\subsubsection{}
Our grading conventions (homological vs. cohomological)
in this paper are at least locally constant, but possibly globally
nonconstant.

\subsubsection{}
We freely use the language of $\infty$-categories.
Our main reference for higher categories and
higher algebra will be \cite{HTT} and \cite{HA}, respectively.
Our main reference for sheaves in the higher categorical
context, in addition to the two already mentioned, will be
\cite{SAG}.

\subsubsection{}
$\pt$ will always denote the terminal map
from a topological space to a point. 
When the domain of $\pt$ is not clear from
context, we denote its domain using a subscript,
e.g. $\pt_X$.

\subsubsection{}
We use $\Aone$ to denote the algebraic
affine line, and $\bfAone$ to denote the analytic
affine line.

\subsubsection{}
$f_*$ and $f^*$ denote the sheaf-theoretic
pushforward and pullback along the map $f$
(as opposed to potentially the $\scrD$-module pushforward
and pullback).

\subsection{Glossary of categories}
${}$\\

\begin{outline}

\0 \noindent General notation for categories.
	\1 $\Ch(\calA)$ denotes the abelian category of chain complexes of objects in the abelian
		category $\calA$.
	\1 $\calD(\calA)$ and its bounded variations denotes the derived $\infty$-categories
	of the abelian category $\calA$.
	\1 If $\calD$ is a stable $\infty$-category with a t-structure:
		\2 $\calD_{\geq 0}$ denotes the connective part of the t-structure.
		\2 $\calD_{\leq 0}$ denotes the coconnective part of the t-structure.
		\2 $\calD^{\heart} := \calD_{\geq 0} \cap \calD_{\geq 0}$ denotes the heart
		of the t-structure.
	\1 $\Seq(\scrC)$ the category of sequences in the stable $\infty$-category $\scrC)$. See \cref{def: Seq(C)}.
		\2 $\Fil(\scrC)$ is the filtered category of the $\infty$-category $\scrC$ admitting
		sequential limits. It is a localization of $\Seq(\scrC)$ by graded equivalences. See \cref{def: Fil(C)}.
	\1 $\Seq(\calA)$ is the category of sequences in the abelian category $\calA$. See \cref{def: Seq(A)}.	
		\2 $\Fil(\calA)$ is the filtered category of the abelian category $\calA$. Its objects consists of sequences in $\calA$
		whose structure maps are monomorphisms. See \cref{def: Fil(A)}.

\0 Special categories.
	\1 $\Sp$ denotes the stable $\infty$-category of spectra.
		\2 $\Sp^*$, where $* = +,-,b$, denote the category of
		left-bounded, right-bounded, and bounded spectra, respectively.
	\1 $\Mod_R$, for $R$ a classical ring, denotes either the stable $\infty$-category of $R$-module spectra if
	$R$ is a commutative ring, or left $R$-module spectra otherwise.
	\1 $\calA_R$, for $R$ a classical ring, denotes the abelian category of $R$-modules if $R$ is commutative, 
	left $R$-modules otherwise. Note that $\calA_R \simeq \Mod_R^{\heart}$ with the canonical t-structure coming
	from spectra on the latter.
	\1 $\calD(R)$ and its bounded variations denote the derived $\infty$-categories of $\calA_R$. They are equivalent
	$\Mod_R$ and its bounded variations. 
	\1 $\Vect$ is special notation for $\Mod_{\C}$ (or $\calD(\C)$).
		\2 $\Perf$ denotes the $\infty$-subcategory of compact objects in $\Vect$,
		also know as perfect complexes.

\0 Categories of sheaves.
	\1 $\Shv(\calX; \scrC)$ denotes the $\infty$-category of $\scrC$-valued sheaves
	on an $\infty$-topos $\calX$.
		\2 $\Shv(X; \scrC)$ denotes the $\infty$-category of sheaves on the $\infty$-topos
		associated to the topological space $X$.
	\1 $\Shv(\calX; R)$ denotes the $\infty$-category of $\Mod_R$-valued sheaves on the $\infty$-topos $\calX$.
		\2 $\Shv^*(\calX; R)$, where $* = -,+,b$, denotes the $\infty$-category of $\Mod_R^*$-valued sheaves on $\calX$.		
	\1 $\Shv_{\wc}(X; \scrC)$ denotes the $\infty$-category of weakly constructible sheaves. See \cref{def: weakly constructible sheaves}. 
		\2 $\Shv_c(X; \scrC)$ denotes the full $\infty$-subcategory of constructible sheaves with respect
		to some coefficient pair $(\scrC, \scrN)$. See \cref{def: constructible sheaves}.
		\2 $\Shv_c^{f.s}(X;\scrC)$ denotes the full $\infty$-category spanned by sheaves that
		are constructible with respect to some finite stratification of $X$.
		\2 $\Shv^A(X;\scrC)$ denotes the full $\infty$-subcategory of sheaves that
		are $A$-constructible, for the stratification $X \to A$. See definition \cref{def: A-constructible sheaves}.
	\1 $\calA(X; R)$ denotes the abelian category of sheaves of $R$-modules on $X$, for $R$ a
		ring.
		\2 $\lLoc(X)$ denotes the full abelian subcategory of $\calA(X;R)$ spanned by locally constant sheaves
		with projective stalks.
	\1 $\calD(X;R)$ and its bounded variations denote the derived $\infty$-categories of $\calA(X;R)$. Note that $\calD(X;R)$ (and
	its bounded variations) is equivalent to $\Shv(X;R)$ (and its bounded variations).
	\1 $\Perv(X)$ denotes the abelian category of perverse sheaves on the complex analytic space $X$. It is the heart of
	 the middle perverse t-structure on $\Shv_c(X; \C)$.
		$\Loc(X)$ denotes the full subcategory of $\Perv(X)$ spanned by local systems
		of finite rank, shifted by the dimension of their supports, so as to be perverse.
	\1 $\calD^b_{\lLoc}(X; \C)$ denotes the full subcategory of $\calD^b(X;\C)$ spanned
	by objects whose cohomology sheaves belong to $\lLoc(X)$.
	\1 $\calD^b_{\Loc}(X; \C)$ denotes the full subcategory of $\Shv_c(X; \C)$ spanned
	by objects whose perverse cohomology sheaves belong to $\Loc(X)$.

\end{outline}

\section{Preliminaries on $\scrD_X$-modules}
\label{sec: Preliminaries on D_X-modules}

\subsection{Definitions}
Let $X$ be a complex manifold, and let $\O_X$ denote
the sheaf of holomorphic functions on $X$, which we regard
as the structure sheaf on $X$. Let $\scrD_X$ denote
the sheaf of differential operators on $X$. By definition,
$\scrD_X$ is the subsheaf of $\C$-algebras of $\sEnd(\O_X)$, the sheaf
of $\C$-linear endomorphisms of the $\O_X$, generated by
$\O_X$, which acts by multiplication, and the sheaf of $\C$-linear
derivations $\Der(X) \subset \sEnd(\O_X)$.

\bdef
A sheaf of $\O_X$-modules $M$ is a (left) $\scrD_X$-module if for every open
set $U \subseteq X$, $M(U)$ has a left $\scrD_X(U)$-module structure, 
compatible with restrictions.
\edefn

There is a definition of right $\scrD_X$-module in which
the sections of $M$ over $U$ are required to have a right $\scrD_X(U)$-module
structure. We will work principally with left $\scrD_X$-modules in
this note, and in the sequel we will refer to them simply as ``$\scrD_X$-modules."

We let $\calA_{\scrD_X}$ denote the category of $\scrD_X$-modules.
The category $\calA_{\scrD_X}$ is abelian, and we denote by
$\calD^b(\scrD_X)$ its bounded derived $\infty$-category.

\bdef
We let $\RegHol_{\scrD_X}$ denote the full subcategory of
$\calA_{\scrD_X}$ spanned by regular holonomic $\scrD_X$-modules.
\edefn

The precise definitions of regularity and holonomicity
are not important to us, so we do not recall them. On the other hand,
it is important to note that the category
$\RegHol_{\scrD_X}$ is abelian, and, by a theorem of Beilinson
(\cite[Theorem 1.3]{BeilinsonPerverse}), the natural map,
\[\calD^b(\RegHol_{\scrD_X}) \xrightarrow{} \calD^b_{\reghol}(\scrD_X)\]
is an equivalence, where the left-hand side denotes the bounded
derived $\infty$-category of $\RegHol_{\scrD_X}$ and the right-hand
side denotes the full $\infty$-subcategory of $\calD^b(\scrD_X)$ spanned
by objects with regular holonomic cohomology.

\subsection{The $V$-filtration}
Given a complex manifold $X$ of dimension $n$, 
and a smooth hypersurface $H \subset X$, 
there is a canonical decreasing filtration on $\scrD_X$ 
called the \emph{$V$-filtration} (also called
the \emph{Kashiwara-Malgrange filtration}).
We recall its definition below.

\bdef
Let $I_H \subset \O_X$ denote the ideal of definition of
the smooth hypersurface $H \subset X$.
The $V$-filtration along $H$ on $\scrD_X$, denoted
$V^{\bullet}\scrD_X$ is a decreasing $\bbZ$-filtration given by,
\[V^{\ell}\scrD_X := \left\{P \in \scrD_X | P \cdot {I_H}^j \subset {I_H}^{j+\ell} \, \text{for all} \, j \in \bbZ \right\}.\] 
\edefn 

In local coordinates $(x,t) := (x_1, \dots x_{n-1}, t)$ 
on $X$, where locally $H = \{t=0\}$,
the $V$-filtration takes the following
form,
\begin{align*}
V^0\scrD_X 	    		&= \left\{\sum_{\a,k} a_{\a,k}(x,t){\partial_x}^{\a} (t\partial_t)^k \big| \a \in \bbN^n, k \in \bbN, a_{\a,k} \in \O_X\right\}, \\
V^k\scrD_X 			&= \begin{cases}
							t^k V^0\scrD_X 							&\text{if} \, k \geq 0, \\ 
							\sum_{j=0}^{-k} {\partial_t}^j V^0\scrD_X 		&\text{if} \, k \leq 0.
					      \end{cases}
\end{align*}

\brem
Observe that $t V^{k}\scrD_X \subset V^{k+1}\scrD_X$, 
and ${\partial_t}V^{k}\scrD_X \subset V^{k-1}\scrD_X$.
On the other hand, all operators in the $x$-direction (i.e. of the form 
$\sum_{j \in \bbN^{n-1}}^m a_j(x) {\partial_x}^j$, for $a_j(x) \in \O_X$)
preserve the pieces of the filtration.
\erem

There is also a related notion of $V$-filtration
along $H$ for $\scrD_X$-modules.

\bdef
\label{def: V-filtration on D-modules}
Given a $\scrD_X$-modules $M$, a 
\emph{Kashiwara-Malgrange} or \emph{$V$-filtration
along $H$} on $M$ is any exhaustive, decreasing $\bbZ$-filtration
on $M$ that satisfies the following properties:
\begin{enumerate}[(i)]
	\item Each $V^k M$ is a coherent module over $V^0\scrD_X$.
	\item $V^k M = t^k V^0 M$ for $k \geq 0$.
	\item $V^k M = \sum_{j=0}^{-k-1} {\partial_t}^j V^{-1} M$ for $k > 0$.
	\item Locally on $X$, there exists a polynomial $b(s) \in \C[s]$ with roots in $\bbQ \cap [0,1)$ such
		that $b(t\partial_t -k)$ vanishes identically on $\gr^k_V M$ for each $k$. 
\end{enumerate}
\edefn

A well-known result of Kashiwara (\cite[Theorem 1]{Kashiwara83}) states
that if a $V$-filtration exists, it is unique. It turns out that
holonomic $\scrD_X$-modules always admit Kashiwara-Malgrange filtrations
along any smooth hypersurface. 
In fact, any morphism $f: M \to N$ between $\scrD_X$-modules
that possess Kashiwara-Malgrange filtrations is strictly compatible with the filtrations
in the sense that $f(V^iM) = f(M) \cap V^iN$.\footnotemark
\footnotetext{In his \href{https://people.math.harvard.edu/~mpopa/notes/DMBG-posted.pdf}{online notes}, 
Mihnea Popa leaves this as an exercise for the reader, as will we.}
We note that the $V$-filtration is compatible with
the standard $V$-filtration on $\scrD_X$

Using the notation established in \cref{sec: Appendix C},
we might summarize the above discussion by
saying that the $V$-filtration 
is a functor on regular holonomic $\scrD_X$-modules
taking values in modules over the filtered sheaf of algebras
$V^{\bullet}\scrD_X$,
\[\RegHol_{\scrD_X} \xrightarrow{V^{\bullet}} \Mod_{V^{\bullet}\scrD_X}(\Fil(\Coh(X))),\]
where $\Coh(X)$ denotes the abelian category of
coherent sheaves on $X$. As discussed in \cref{sec: Appendix C},
$\Fil(\Coh(X))$ is not abelian, though it is additive and
admits kernels and cokernels. In general, however, the
the natural map from the coimage of a morphism to the
image is not an isomorphism. Nonetheless, we may still
talk about short exact sequences in $\Mod_{V^{\bullet}\scrD_X}(\Fil(\Coh(X)))$ and
exact functors to and from $\Mod_{V^{\bullet}\scrD_X}(\Fil(\Coh(X)))$.

\blem
The $V$-filtration functor, $V^{\bullet}$
is an exact functor.
\elem

\bproof
It is clear that $V^{\bullet}$ is an additive functor.
That it preserves short exact sequences follows
from the strict compatibility of maps of $\scrD_X$-modules
with the $V$-filtration.
\eproof

\subsection{de Rham functors}
Let $X$ be a complex manifold of dimension
$n$. For $0 \leq p \leq n$, we let $\Omega_X^p$
denote the sheaf of degree $p$ holomorphic differential
forms on $X$. There is a canonical right $\scrD_X$-module 
structure on $\Omega_X^n$. We will use $\Omega_X$ to
denote $\Omega_X^n$ as a $\scrD_X$-module.

\bdef
The \emph{de Rham functor} on $\calD^b(\scrD_X)$,
denoted $\DR_X$, is defined to be
\[\Omega_X \otimes_{\scrD_X} -: \calD^b(\scrD_X) \to \Shv(X; \C).\]
\edefn

The right $\scrD_X$-module $\Omega_X$ is also the
transfer bimodule $\scrD_{\ast \xleftarrow{} X}$. There is
a canonical locally free resolution of $\Omega_X$ as a right $\scrD_X$-module
given by the complex
\[0 \to \Omega_X^0 \otimes_{\O_X} \scrD_X \to \cdots \to \Omega_X^n \otimes_{\O_X} \scrD_X \to \Omega_X \to 0.\]
More generally, given a decomposition $X = Y \times Z$,
we may consider the transfer bimodule, $\scrD_{Y \xleftarrow{} X}$.

\bdef
The \emph{relative de Rham} functor for $\pr: X \to Y$
on $\calD^b(\scrD_X)$, denoted
$\DR_{X/Y}$, is defined to be
\[\scrD_{Y \xleftarrow{} X} \otimes_{\scrD_X} -: \scrD^b(\scrD_X) \to \calD^b(\pr^*\scrD_Y),\]
where $\pr^*$ denote the pullback of sheaves
(rather than the pullback of $\scrD_X$-modules).
\edefn

The bimodule $\scrD_{Y \xleftarrow{} X}$
has a resolution as a right $\scrD_X$-module,
similarly given by tensoring with the \emph{relative} differential forms on $X$.
Set $d = \dim Z$, and let $\Omega_{X/Y}^k$ denote the sheaf
of relative differential forms on $X$, for $0 \leq k \leq d$.
Then the complex,
\[0 \to \Omega_{X/Y}^0 \otimes_{\O_X} \scrD_X \to \cdots \to \Omega_{X/Y}^d \otimes_{\O_X} \scrD_X \to \scrD_{Y \xleftarrow{} X} \to 0\]
is a locally free resolution of $\scrD_{Y \xleftarrow{} X}$.
Using this locally free resolution, we have, for any $M \in \calA_{\scrD_X}$,
\beqn
\label{eqn: resolution of relative de Rham complex}
\DR_{X/Y}(M) \simeq \left[0 \to \Omega_{X/Y}^0 \otimes_{\O_X} M \to \cdots \to \Omega_{X/Y}^d \otimes_{\O_X} M \to 0\right]
\eeqn
as objects of the derived category $\calD^b(\pr^*\scrD_Y)$.

\subsection{The Riemann--Hilbert correspondence}
Let $X$ be a complex manifold of dimension $n$.
Kashiwara's constructibility theorem (\cite[Theorem 4.6.3]{Hotta})
states that the image of a holonomic $\scrD_X$-module
under the de Rham functor $\DR_X$ lies in the 
subcategory $\Shv_c(X;\C) \subset \Shv(X;\C)$.
In fact, $\DR_X$ induces an equivalence between
$\calD^b_{\reghol}(\scrD_X)$ and the
stable $\infty$-category of constructible sheaves.

\bthm[The Riemann--Hilbert correspondence]
The de Rham functor induces a t-exact equivalence
of stable $\infty$-categories,
\[\DR_X: \calD^b_{\reghol}(\scrD_X) \xrightarrow{\simeq} \Shv_c(X;\C)\]
with respect to the standard t-structure on the left-hand
category and the middle perverse t-structure on the right-hand
category.
\ethm

\subsubsection{Deligne's correspondence}
The Riemann--Hilbert correspondence generalizes
the classical correspondence of Deligne between
integrable connections and local systems.

\bdef
We say $M \in \calA_{\scrD_X}$ is an integrable connection
if its underlying $\O_X$-module is locally free of finite rank.
We denote the full subcategory of $\calA_{\scrD_X}$ spanned
by integrable connections by $\Conn(X)$.
\edefn

\bdef
A local system on $X$ is a perverse 
sheaf on $X$ of the form $\scrL[n]$, where
$\scrL$ is a locally constant sheaf of 
$\C$-vector spaces of finite rank. We denote
by $\Loc(X)$ the full subcategory of $\Perv(X)$ spanned by
local systems.
\edefn

Let $D \subset X$ be a divisor (complex hypersurface) on $X$.
We briefly recall the notion of \emph{regular meromorphic
connection on $X$ along $D$}. 

A meromorphic connection on $X$ along $D$ is
a coherent $\O_X(D)$-module with a $\C$-linear connection.
Heuristically, a meromorphic connection $M$ on the unit open
disk $B \subset \C$ along $0$
is regular if its associated connection matrix has poles of
at most order $1$ at $0$. A meromorphic connection $M$ on $X$ 
along $D$ is regular if for any map $u: B \to X$ such
that $u^{-1}(D) = \{0\}$, the germ of $i^*M$ is 
regular along $0 \subset B$. We refer the
reader to \S5.1.2 and \S5.2 in
\cite{Hotta} for precise definitions.

Let $\Conn(X;D)$ denote the category of
meromorphic connections on $X$ along $D$,
and let $\Conn^{\reg}(X;D)$ denote the full
subcategory of $\Conn(X;D)$ spanned by regular meromorphic
connections along $D$. Restriction of a meromorphic connection $M$ on $X$
along $D$ to the complement of the divisor defines a functor,
$\Del_{(X,D)}: \Conn(X;D) \to \Conn(X \setminus D)$,
i.e. the restriction $M|_{X \setminus D}$ is locally free.
A landmark result of Deligne (\cite{Deligne}) states that
the restriction of $\Del_{(X;D)}$ to the subcategory of regular connections
induces an equivalence,
\[\Del_{(X;D)}: \Conn^{\reg}(X;D) \xrightarrow{\simeq} \Conn(X \setminus D).\]

\brem
Regular meromorphic connections on $X$ along $D$ define
regular holonomic $\scrD_X$-modules. In fact,
$\Conn^{\reg}(X; D)$ is an abelian
subcategory of $\RegHol_{\scrD_X}$.  
\erem

On the other hand, a simple computation\footnotemark shows that
the de Rham functor induces an equivalence of abelian
categories,
\[\DR_X: \Conn(X) \xrightarrow{\simeq} \Loc(X).\]
\footnotetext{See the discussion
preceding \cite[Theorem 4.2.4]{Hotta} for details.}
This observation allows us to state what is 
known as \emph{Deligne's correspondence}.

\bthm[Deligne's correspondence]
Let $X$ be a complex manifold, and let
$D \subset X$ be a divisor. Then there is an equivalence,
\[\DR_X \circ \, {\Del_{(X;D)}}: \Conn^{\reg}(X;D) \xrightarrow{\simeq} \Loc(X \setminus D).\] 
\ethm

We will be interested in the special case when 
$X = \bfAone$ and $D = \{0\}$, in which case Deligne's
correspondence produces an equivalence,
$\Conn^{\reg}(\bfAone; 0) \xrightarrow{\simeq} \Loc(\bfAone \setminus 0)$.
Taking bounded derived $\infty$-categories,
this equivalence induces a t-exact equivalence,
\[\calD^b(\Conn^{\reg}(\bfAone; 0)) \xrightarrow{\simeq} \calD^b(\Loc(\bfAone \setminus 0)).\]
By \cref{cor: beeg corollary}, the right-hand category is
t-exact equivalent to $\calD^b_{\Loc}(\bfAone \setminus 0;\C)$,
where the latter is the full subcategory of $\Shv_c(\bfAone \setminus 0;\C)$
spanned by objects whose perverse cohomology sheaves
belong to $\Loc(\bfAone \setminus 0)$.
This fact combined with the Riemann--Hilbert correspondence
yields the following lemma.

\blem
\label{lem: derived category of conn reg}
Let $\calD^b_{\regconn}(\scrD_{\bfAone})$ denote the
full subcategory of $\calD^b_{\reghol}(\scrD_{\bfAone})$
spanned by objects whose cohomology
lies in $\Conn^{\reg}(\bfAone; 0)$.
Then there is a t-exact equivalence 
\[\calD^b(\Conn^{\reg}(\bfAone; 0)) \xrightarrow{\simeq} \calD^b_{\regconn}(\scrD_{\bfAone}).\]
\elem

Combining \cref{lem: derived category of conn reg} above with 
\cref{lem: pairs} and \cref{cor: beeg corollary} in \cref{sec: Appendix B}, 
we obtain the following t-exact equivalence, 
which is a kind of derived version of Deligne's correspondence.

\bcor[Derived Deligne's correspondence]
\label{cor: derived Deligne's correspondence}
There is a t-exact equivalence of stable $\infty$-categories,
\beqn
\label{eqn: local Riemann--Hilbert}
\calD^b_{\regconn}(\scrD_{\bfAone}) \xrightarrow{\simeq} \Fun(S^1, \Perf),
\eeqn
where the t-structure on the right-hand category
is the one coming from the perverse t-structure on
$\Shv_{\Loc}(\bfAone \setminus 0;\C)$ under the
equivalence obtained in \cref{lem: pairs}.
\ecor

\brem
Note that the above discussion holds for \emph{any} subset
of $\bfAone$ whose homotopy type is $K(\bbZ, 1)$.
In particular, \labelcref{eqn: local Riemann--Hilbert}
holds for the complex analytic disk 
$D_{\varepsilon}(0) \subset \bfAone$ of radius
$\varepsilon > 0$ in place of $\bfAone$, and
$\mathring{D}_{\varepsilon}(0)$ in place of $\bfAone \setminus 0$.
\erem

\section{$\scrD_{\C,0}$-modules}
\label{sec: D_{C,0}-modules}

\bdef
Let $\scrD_{\C,0}$ denote the stalk at $0$ of
the sheaf of differential operators $\scrD_{\bfAone}$.
As a module, it is $\C\{t\}\langle \partial_t \rangle$,
where $\C\{t\}$ denotes the ring of germs of holomorphic
functions in the variable $t$; its product
structure is given by the Leibniz rule.
\edefn

Since $\scrD_{\C,0}$ is an ordinary associative ring,
it can be regarded as a discrete $\bbE_1$-ring. 
Moreover, the $\infty$-category of left $\scrD_{\C,0}$-module
spectra is equivalent
to the unbounded derived $\infty$-category, 
$\calD(\scrD_{\C,0})$, of the abelian
category of left $\mathscr{D}_{\C,0}$-modules, $\calA_{\scrD_{\C,0}}$
(see \cite[Remark 7.1.1.16]{HA}).

\brem
\label{rem: monoidal structure on D_{C,0}-modules}
Both $\calA_{\scrD_{\C,0}}$ and its derived category
$\calD^b(\scrD_{\C,0})$ inherit symmetric monoidal
structure products (the latter the derived functor of the former) 
coming from the canonical tensor product\footnotemark of $\scrD_{\bfAone}$-modules,
since every $\scrD_{\C,0}$-module (i.e. object in $\calA_{\scrD_{\C,0}}$)
has a representative in $\calA_{\scrD_{\bfAone}}$ of which it is the
germ at $0 \in \bfAone$.
\footnotetext{See e.g. \cite[Proposition 1.2.9]{Hotta} and \cite[pg. 38-40]{Hotta} for the definition
of tensor products of $\scrD$-modules.} 
\erem

\subsection{Regular holonomic $\scrD_{\C,0}$-modules}
Given an ordinary $\scrD_{\bfAone}$-module $M$,
its germ at $0 \in \C$ is naturally a $\scrD_{\C,0}$-module.
This feature provides a convenient way of
defining the notion of regular holonomic 
$\scrD_{\C,0}$-module.

\bdef
Given an ordinary $\scrD_{\C,0}$-module, $M \in \calA_{\scrD_{\C,0}}$,
we say that $M$ is \emph{regular holonomic} if
there a regular holonomic $\scrD_{\bfAone}$-module 
$\wit{M}$ such that the germ of $\wit{M}$ 
at $0 \in \C$ is $M$.
\edefn

We denote the abelian category of regular holonomic
$\scrD_{\C,0}$-modules by $\RegHol_{\scrD_{\C,0}}$.
This category turns out to be equivalent to the category
of holonomic $\scrD$-modules on a disc around
the origin that are regular at $0 \in \bfAone$ and restrict
to an integrable connection on the punctured disc. 
Let $D_{\varepsilon} \subset \bfAone$ denote
the disc of radius $\varepsilon$ around $0$ as above.

\bprop[{e.g.\ \cite[\S5.3]{Dimca}}]
\label{prop: descriptions of D-modules}
The following categories are naturally equivalent,
\begin{enumerate}[(i)]
\item $\RegHol_{\scrD_{\C,0}}$
\item The category of holonomic $\scrD_{D_{\varepsilon}}$-modules
$M$ such that $M|_{\mathring{D}_{\varepsilon}}$ is an integrable
connection that is regular along $0 \in D_{\varepsilon}$.
\item The category of holonomic $\scrD_{\bbP^1}$-modules
$M$ (either algebraic or analytic) 
such that $M|_{\bbP^1 \setminus \{0, \infty\}}$ is an
integrable connection regular at $0$ and $\infty$.
\item The category of holonomic $A_t$-modules $M$,
regular at $0$ and at $\infty$ such that $M|_{\Aone \setminus 0}$
is an integrable connection. Here $A_t = \C[t]\langle \partial_t \rangle$
is the Weyl algebra of polynomial linear differential
operators.
\end{enumerate}
The equivalences between (iii) and (i),
and (ii) and (i), are induced by taking germs; the
equivalence between (iii) and (ii) is induced by
restriction. The equivalence between (i) and (iv)
is induced either by restriction along the ring
inclusion $A_t \hook \scrD_{\C,0}$, or base change
$\scrD_{\C,0} \otimes_{A_t} -$, which adjoint equivalences.
\eprop

\subsection{Riemann--Hilbert for $\scrD_{\C,0}$-modules}
Using Deligne's correspondence we obtain 
a functor,
\beqn
\label{eqn: Deligne's correspondence for D_{C,0}}
\Loc(\bfAone \setminus 0) \xrightarrow{\simeq} \Conn^{\reg}(\C;0) \subset \RegHol_{\scrD_{\C,0}}.
\eeqn

\brem
\label{rem: calculation of RH}
Let $(V,T) \in \Loc(\bfAone \setminus 0)$ be a finite-dimensional 
$\C$-linear representation
of $\pi_1(\bfAone \setminus 0)$. Under the equivalence \labelcref{eqn: Deligne's
correspondence for D_{C,0}}, $(V,T)$ is sent to the
$\scrD_{\C,0}$-module $\C\{t\} \otimes_{\C} V$,
with the action of $\partial_t$ given by 
$\nabla_{\partial_t}(f \otimes v) = df(t) \otimes v + \frac{f}{t} \otimes Mv$, where $M: V \to V$
is such that $\exp(-2\pi iM) = T$. In order words,
under Deligne's correspondence, $(V,T)$ corresponds
to the $\C\{t\}$-module $\C\{t\} \otimes_{\C} V$
equipped with connection $\nabla = d + \frac{M}{t}dt$.
It is easy to see that the monodromy of this
connection is $T$ (see, e.g. \cite[\S4.b]{SabbahBook}).
\erem

\bdef
We denote by $\calD^b_{\reghol}(\scrD_{\C,0}) \subset
\calD^b(\scrD_{\C,0})$ the full subcategory of $\calD^b(\scrD_{\C,0})$
spanned by objects whose coholomogies are regular holonomic
$\scrD_{\C,0}$-modules.
\edefn

The following lemma is a direct consequence of
\cref{prop: descriptions of D-modules} and \cref{lem: derived
category of conn reg}. It can be thought of a
local version on $\bfAone$ of Beilinson's famous
\cite[Theorem 1.3]{BeilinsonPerverse} showing that
the derived category of regular holonomic modules
is equivalent to the category of complexes with
regular holonomic cohomologies.

\blem
\label{lem: modified Beilinson}
There is a t-exact equivalence,
$\calD(\RegHol_{\scrD_{\C,0}}) 
	\xrightarrow{\simeq} \calD_{\reghol}(\scrD_{\C,0})$.
\elem

Combined with \cref{lem: modified Beilinson}, 
the derived Deligne's correspondence
\labelcref{eqn: local Riemann--Hilbert} 
induces the following functor on stable $\infty$-categories.

\blem[Deligne's correspondence for $\scrD_{\C,0}$-modules]
\label{lem: Deligne's correspondence for D_{C,0}-modules}
There is a t-exact functor of stable
$\infty$-categories,
\beqn
\label{eqn: Riemann--Hilbert for D_{C,0}}
\Temp: \Fun(S^1, \Perf) \to \calD^b_{\reghol}(\scrD_{\C,0}),
\eeqn
where the t-structure on
the left-hand category is the one used in the statement
of \cref{cor: derived Deligne's correspondence}.
\elem

\brem
The functor $\Temp$ is given by composition
of an inverse functor to the derived Deligne's
correspondence with the inclusion
$\calD^b_{\regconn}(\scrD_{C,0}) \subset \calD^b_{\reghol}(\scrD_{\C,0})$.
The notation chosen to denote this functor 
alludes to the notation ``$\text{RH}$" used by Kashiwara to denote
the explicit inverse to the de Rham functor he constructs
in his proof of the Riemann--Hilbert correspondence.
(see \cite{KashiwaraTempered}).
\erem

\subsection{Fourier transform}
The category of $\scrD_{\C,0}$-modules
has a Fourier transform coming from its
description as holonomic modules over the
Weyl algebra in \cref{prop: descriptions of D-modules}.

\bdef[{e.g.\ \cite{GarciaLopez}}]
The Fourier anti-involution is the map
of $\C$-algebras $A_t \to A_{\eta}$ given
by the assignment on generators,
\[t \mapsto -\partial_{\eta}, \hspace{2cm} \partial_t \mapsto \eta.\]
The \emph{Fourier transform} of an object
$M \in \Mod_{A_t}$ is defined as the object,
\[\Four(M) := A_{\eta} \otimes_{A_t} M\]
in $\Mod_{A_{\eta}}$,
where $A_{\eta}$ is regarded as a right
$A_t$-module via the Fourier anti-involution.
\edefn

\brem
We will regard the Fourier transform as an
involution of the category of $A_t$-module
spectra, $\Four: \Mod_{A_t} \to \Mod_{A_t}$.
\erem

By \cite[Proposition 3.2.7]{Hotta}, a coherent, ordinary $A_t$-module
is holonomic if and only if its Fourier transform is holonomic, as well.
On the other hand, the Fourier transform of a regular
$A_t$-module is not necessarily regular. Fortunately,
the Fourier transforms of the highly structured sort of $A_t$-modules
that correspond to regular holonomic $\scrD_{\C,0}$-modules
\emph{are} regular, by the following lemma.

\blem
\label{lem: Fourier is regular}
Suppose that $M$ is holonomic
$A_t$-module that is regular at $0$
and $\infty$ and such that $M|_{\Aone \setminus 0}$
is an integrable connection. Then $\Four(M)$ is also
a holonomic $A_t$-module that is
regular at $0$ and $\infty$ and such that
$\Four(M)|_{\Aone \setminus 0}$ is an integrable connection.
\elem

\bproof
As mentioned, $\Four(M)$ is holonomic
if and only if $M$ is. It remains to show
that $\Four(M)$ is regular at $0$ and
$\infty$.

By a theorem of Brylinski (\cite{Brylinski}),
$\Four(M)$ is regular if $M$ is a
\emph{monodromic}\footnotemark
\footnotetext{See \cite[Definition 5.3.8]{Dimca}.
Alternatively, a $A_t$-module is monodromic
if its image under Riemann--Hilbert is locally
constant along $\C^*$-orbits.}
$A_t$-module.\footnotemark
\footnotetext{In fact the converse is true. See \cite[Theorem 1.2]{Ito}.}
On the other hand, by \cite[Proposition 5.3.9]{Dimca},\footnotemark
\footnotetext{Also \cite{DimcaSaito}.}
an $A_t$-module $M$ is monodromic 
if and only if it is regular holonomic 
and its restriction to $\Aone \setminus 0$ is an integrable connection.
Thus, under our hypotheses,
$\Four(M)$ is regular holonomic.

We note that $\Four(M)$ is also monodromic by \cite[Theorem 1.2]{Ito},
so it follows from another application of \cite[Proposition 5.3.9]{Dimca}
that $\Four(M)$ is an integrable connection on
$\Aone \setminus 0$. It follows that $\Four(M)$
has regular singularities contained
in the complement of $\Aone \setminus 0 \subset \bbP^1$,
at $0$ and $\infty$.
\eproof

It follows from \cref{lem: Fourier is regular}
that the Fourier transform restricts to
an involution,\footnotemark
\footnotetext{See also the definition 
in the Introduction of \cite{StationaryPhase}
of what the author calls the local formal Laplace
transform.}  
\[\Four: \calD_{\reghol}(\scrD_{\C,0}) 
\to \calD_{\reghol}(\scrD_{\C,0}).\]

\subsection{Vanishing cycles}
We recall the definition of vanishing 
cycles for regular holonomic\footnotemark $\scrD_{\bfAone}$-modules.
\footnotetext{The definition we give is valid more generally
for all $\scrD_{\bfAone}$-modules that are specializable
along $0 \in \bfAone$, but we will not need this.}
Let $M \in \RegHol_{\scrD_{\bfAone}}$,
and consider the Kashiwara-Malgrange
filtration with respect to the hypersurface $\{0\} \subset \C$.
The vanishing cycles of $M$, denoted $\varphi_{\scrD}(M)$ 
is defined to be the $-1$st associated graded piece of this filtration.
The result is a local system on the punctured plane: 
a vector space equipped with
a monodromy automorphism. It is equivalent
to $\varphi_t \DR_{\bfAone}(\wit{M})$, where $\wit{M}$
is a $\scrD_{\bfAone}$-module representative for
$M$ and $\varphi_t$ is the usual vanishing cycles functor\footnotemark
\footnotetext{See \cite[\S4.2]{Dimca} for the definition. We use 
the convention that $\varphi_t$ preserves perversity.}
on perverse sheaves. The vanishing cycles functor on
$\RegHol_{\scrD_{\C,0}}$ induces the vanishing cycles
functor on the derived category, 
$\varphi_{\scrD}: \calD^b_{\reghol}(\scrD_{\C,0}) 
\to \Fun(S^1, \Perf)$, which, by abuse of 
notation, we also denote by $\varphi_{\scrD}$.

\section{The $V$-filtration on $\scrD_{\C,0}$-modules}
\label{sec: The V-filtration on D_{C,0}-modules}

The standard $V$-filtration on $\scrD_{\bfAone}$ along the hypersurface 
$\{0\} \subset \C$ induces a filtration on $\scrD_{\C,0}$, also called the
standard $V$-filtration, which we we denote by $V^{\bullet}\scrD_{\C,0}$.
Concretely, the filtration is given by
\begin{align*}
V^0\scrD_{\C,0} 	    		&= \C\{t\}\langle t\partial_t\rangle, \\
V^k\scrD_{\C,0} 			&= \begin{cases}
							t^k V^0\scrD_{\C,0} 								&\text{if} \, k \geq 0, \\ 
							\sum_{j=0}^{-k} {\partial_t}^j V^0\scrD_{\C,0} 			&\text{if} \, k \leq 0.
					      \end{cases}
\end{align*}

There is a corresponding notion of $V$-filtration
on $\scrD_{\C,0}$-modules as well, which differs slightly from what
one might expect by comparison with \cref{def: V-filtration on D-modules}.

\bdef
\label{def: usual KM filtration}
Given a $\scrD_{\C,0}$-module $M$, a 
\emph{Kashiwara-Malgrange} or \emph{$V$-filtration} on $M$ 
is any exhaustive, decreasing $\bbZ$-filtration
on $M$ that satisfies the following properties:
\begin{enumerate}[(i)]
	\item Each $V^k M$ is a coherent module over $\C\{t\}$.
	\item $V^k M = t^k V^0 M$ for $k \geq 0$.
	\item $V^k M = \sum_{j=0}^{-k-1} {\partial_t}^j V^{-1} M$ for $k > 0$.
	\item There exists a polynomial $b(s) \in \C[s]$ with roots in $\bbQ \cap [0,1)$ such
		that $b(t\partial_t -k)$ vanishes identically on $\gr^k_V M$ for each $k$. 
\end{enumerate}
\edefn

Just as for the $V$-filtration on $\scrD_{\bfAone}$-modules,
if a $V$-filtration on $M \in \calA_{\scrD_{\C,0}}$
exists, then it is unique, so we speak of \emph{the} Kashiwara-Malgrange
filtration. Crucially, the $V$-filtration on
$M$ is compatible with the standard $V$-filtration on $\scrD_{\C,0}$, making
$V^{\bullet}M$ a module over $V^{\bullet}\scrD_{\C,0}$.

If $M \in \RegHol_{\scrD_{\C,0}}$, then
the $V$-filtration on 
$M$ exists, so we obtain an additive functor,
\[V^{\bullet}: \RegHol_{\scrD_{\C,0}} \to \Mod_{V^{\bullet}\scrD_{\C,0}}(\Fil(\Ab)).\]
The additive category $\Mod_{V^{\bullet}\scrD_{\C,0}}(\Fil(\Ab))$
is in fact quasi-abelian, so it makes sense the
state the following lemma.

\blem
\label{lem: V-filtration is exact}
The above functor $V^{\bullet}$ is fully faithful, 
exact, and preserves filtered colimits. 
\elem

\bproof
We recall that any morphism $f: M \to N$ of $\scrD_{\C,0}$-modules
induces a strict map of filtered $V^{\bullet}\scrD_{\C,0}$-modules, 
$V^{\bullet}M \to V^{\bullet}N$, if the Kashiwara-Malgrange 
filtration exists on both. It follows that $V^{\bullet}$ is fully faithful.
Because $V^{\bullet}$ takes morphisms of $\scrD_{\C,0}$-modules
to strict morphisms, it suffices to show that it preserves injections and
surjections in order to show it is exact. But this is clear from the fact
that, given $f: M \to N$, $f(V^iM) = f(M) \subset V^iN$. 
Finally, we show that $V^{\bullet}$ preserves filtered---in fact,
all---colimits. On the one hand, it follows from the definition of the 
$V$-filtration that $V^{\bullet}$ preserves coproducts;
indeed, one can check that setting $V^k(M \oplus N) := V^kM \oplus V^kN$
defines a filtration which satisfies all the properties of the Kashiwara-Malgrange
filtration. On the other hand, a functor preserves all colimits if and only if it preserves
coequalizers and coproducts. Since $V^{\bullet}$ is additive, it preserves coequalizers
if and only if it preserves cokernels, which was shown above.
\eproof

Taking $R = \scrD_{\C,0}$ in \cref{prop: filtered R-modules
into filtered R-module spectra}, we obtain a fully faithful
functor $i_{\heart}: \Mod_{V^{\bullet}\scrD_{\C,0}}(\Fil(\Ab)) \hook \Mod_{V^{\bullet}\scrD_{\C,0}}(\Fil(\Sp))$
which takes exact sequences to fiber sequences and preserves filtered colimits.
Since $\Mod_{V^{\bullet}\scrD_{\C,0}}(\Fil(\Sp))$
admits filtered (in fact, all) colimits,
the composition $i_{\heart} \circ V^{\bullet}$
extends to a functor,
\[\Ind(\RegHol_{\scrD_{\C,0}}) \to \Mod_{V^{\bullet}\scrD_{\C,0}}(\Fil(\Sp)),\]
that preserves filtered colimits, which by
abuse of notation, we denote by $V^{\bullet}$.

If $C$ is an ordinary $1$-category, then $\Ind(C)$ is
also a $1$-category, and coincides with the standard ind-category
of $C$. Moreover, it is well-known that if $\calA$ is an abelian
category, $\Ind(\calA)$ is a \emph{Grothendieck} abelian category.
In particular, $\Ind(\RegHol_{\scrD_{\C,0}})$ is a Grothendieck
abelian category.

\blem
The functor $V^{\bullet}: \Ind(\RegHol_{\scrD_{\C,0}}) \to 
\Mod_{V^{\bullet}\scrD_{\C,0}}(\Fil(\Sp))$
takes exact sequences to fiber sequences.
\elem

\bproof
This is a direct consequence of the fact that $V^{\bullet}$
is induced by an exact functor out of $\RegHol_{\scrD_{\C,0}}$.
\eproof

We recall a proposition of Lurie that characterizes that will
allows us to extend $V^{\bullet}$ to a functor on the unseparated
derived $\infty$-category of $\RegHol_{\scrD_{\C,0}}$.

\bprop
Let $\calA$ be a Grothendieck abelian category, and let
$\scrC$ be a stable presentable $\infty$-category. Then
then there is an equivalence between 
colimit-preserving functors $\widecheck{\calD}(\calA)
\to \scrC$, on the one hand, and functors $\calA \to
\scrC$ that preserve filtered colimits and send exact sequences
to fiber sequences, on the other.
\eprop

\bproof
See \cite{LurieMOanswer}.
\eproof

By the above lemma and proposition, we obtain a colimit-preserving
functor,
\[\widecheck{\sfV}^{\bullet}: \widecheck{\calD}(\Ind(\RegHol_{\scrD_{\C,0}}))
\to \Mod_{V^{\bullet}\scrD_{\C,0}}(\Fil(\Sp)).\] 
Restricting $\widecheck{\sfV}^{\bullet}$ to the full subcategory
of bounded complexes, $\widecheck{\calD}^b(\Ind(\RegHol_{\scrD_{\C,0}}))$,
we obtain another colimit-preserving functor, which we also denote
by $\widecheck{\sfV}^{\bullet}$.

\blem
\label{lem: bounded unseparated = bounded separated}
Let $\calA$ be a Grothendieck abelian category.
Then there is a t-exact equivalence 
$\widecheck{\calD}^+(\calA) \xrightarrow{\simeq} \calD^+(\calA)$.
\elem

\bproof
The following argument is due to Sam Raskin.
More generally, let $\scrC$ be a presentable stable $\infty$-category 
with a t-structure compatible with filtered colimits, and let $\calA$
denote its heart. In order to check that the induced 
functor $\calD^+(\calA) \to \scrC^+$ is an equivalence
it suffices to check that for every pair of objects $F, I \in A$ 
with $I$ injective, $H^i\Hom_{\scrC}(F,I) = 0$ for all $i>0$,
where $\calA$ is considered as a subcategory of $\scrC$.
Since this is obviously true taking $\scrC$ to be either
$\widecheck{\calD}(\calA)$ or $\calD(\calA)$ itself, we are done.
\eproof

Using \cref{lem: bounded unseparated = bounded separated},
we view $\widecheck{\sfV}^{\bullet}$ as a functor on $\calD^b(\Ind(\RegHol_{\scrD_{\C,0}}))$. 
Finally, composition
with the functor $\calD^b(\RegHol_{\scrD_{\C,0}}) \to
\calD^b(\Ind(\RegHol_{\scrD_{\C,0}}))$, induced by the
exact embedding $\RegHol_{\scrD_{\C,0}} \hook 
\Ind(\RegHol_{\scrD_{\C,0}})$, gives the functor
\beqn
\label{eqn: derived KM filtration}
\sfV^{\bullet}: \calD^b_{\reghol}(\scrD_{\C,0}) \simeq \calD^b(\RegHol_{\scrD_{\C,0}}) \to 
\Mod_{V^{\bullet}\scrD_{\C,0}}(\Fil(\Sp)),
\eeqn
which we call the \emph{derived Kashiwara-Malgrange
filtration}.
Given an object $M \in \calD^b_{\reghol}(\scrD_{\C,0})$,
we will refer to $\sfV^{\bullet}M$ as the Kashiwara-Malgrange
filtration on $M$.

\bdef
A \emph{filtered derived $\scrD_{\C,0}$-module} is
an object of $\Mod_{V^{\bullet}\scrD_{\C,0}}(\Fil(\Sp))$.
\edefn

\brem
Note that the structure of a left module over $V^{\bullet}\scrD_{\C,0}$
on an object $M^{\bullet} \in \Fil(\Sp)$
includes the data of a factorization of the multiplication map
as follows,
\[\begin{tikzcd}
	& {M^{i+j}} \\
	{V^i\scrD_{\C,0} \times M^j} && {M^{\bullet}}
	\arrow["\mult", from=2-1, to=2-3]
	\arrow[from=1-2, to=2-3]
	\arrow[dashed, from=2-1, to=1-2]
\end{tikzcd}\]
\erem

We let $\sfV^k: \calD^b_{\reghol}(\scrD_{\C,0}) \to \Mod_{V^0\scrD_{\C,0}}$ 
denote the functor obtained obtained similarly
from the universal property of the unseparated derived
$\infty$-category using the composition of 
$(i_{\heart} \circ V^{\bullet}): 
\Ind(\RegHol_{\scrD_{\C,0}})
\to \Mod_{V^{\bullet}\scrD_{\C,0}}(\Fil(\Sp))$ with 
with the evaluation functor $\ev_k: \Mod_{V^{\bullet}\scrD_{\C,0}}(\Fil(\Sp))
\to \Mod_{V^0\scrD_{\C,0}}$ (see \cref{rem: ev_k on
modules}), noting that the latter functor also commutes with
filtered colimits and takes fiber sequences to fiber sequences.
By definition, $\sfV^k$ is a t-exact functor with respect to the
canonical t-structures on both the source and the target.

\blem
\label{lem: kth piece of filtration}
The functors $\sfV^k$ and $\ev_k \circ \sfV^{\bullet}$
are naturally equivalent.
\elem

\bproof
Note that the functors,
\[\widecheck{\sfV}^k, \ev_k \circ \widecheck{\sfV}^{\bullet}: 
	\widecheck{\calD}(\Ind(\RegHol_{\scrD_{\C,0}})) 
		\to \Mod_{V^0\scrD_{\C,0}},\]
both preserve colimits, and both restrict to the
functor, $\ev_k \circ V^{\bullet}: \Ind(\RegHol_{\scrD_{\C,0}}) \to
\Mod_{V^0\scrD_{\C,0}}$, which preserves
filtered colimits and takes exact sequences to fiber
sequences. By the universal property of $\widecheck{\calD}$,
it follows that $\widecheck{\sfV}^k \simeq \ev_k \circ \widecheck{\sfV}^{\bullet}$.
Restriction the full subcategory $\calD^b_{\reghol}(\scrD_{\C,0})$ yields the
result.
\eproof

In the same way as we obtain $\sfV^k$, we may obtain
a t-exact functor $\gr_{\sfV}^k: \calD^b_{\reghol}(\scrD_{\C,0}) 
\to \Vect_{\C}$.

\blem
The functors $\gr_{\sfV}^k$ and $\pr_k \circ \gr_{\sfV}^{\bullet}$
are naturally equivalent.
\elem

\bproof
Proven in the same manner as \cref{lem: kth piece of filtration}.
\eproof

We show in the following proposition that 
the derived $V$-filtration produces a filtered
derived $\scrD_{\C,0}$-module which
satisfies a list of properties analogous to those
satisfied by the
$V$-filtration on the abelian category of regular
holonomic modules.

\bprop
\label{prop: sfV properties}
Suppose that $M \in \calD^b_{\reghol}(\scrD_{\C,0})$. 
Then $\sfV^{\bullet}M$ satisfies the
following properties:
\begin{enumerate}[(i)]
	\item \label{I} $\sfV^k M$ is a perfect $\C\{t\}$-module.
	\item \label{II} The map, $t \cdot: \sfV^k M \to \sfV^{k+1} M$ is an equivalence for each $k \geq 0$.
	\item \label{III} The induced map $\sum_{j=0}^{-k-1} {\partial_t}^j \sfV^{-1} M \to \sfV^k M$
	is an equivalence for each $k < 0$.
	\item \label{IV} For each $k \in \bbZ$, $\gr^k_{\sfV}M$ is a perfect
	$\C$-module, and there exists a nonzero polynomial $b(s) \in \C[s]$ 
	whose roots have real parts in $[0,1)$ and a factorization,
\[\begin{tikzcd}
	& {\sfV^{k+1}M} \\
	{\sfV^k M} && {M}
	\arrow["{b(t\partial_t-k)}", from=2-1, to=2-3]
	\arrow["{i_{k+1}}", from=1-2, to=2-3]
	\arrow[dashed, from=2-1, to=1-2].
\end{tikzcd}\]
\end{enumerate}
\eprop

\bproof
We show that $\sfV^{\bullet}$ satisfies the enumerated
properties one-by-one.

\begin{enumerate}[(i)]

\item Since $\sfV^k$ is t-exact by the discussion preceding the
proposition, we have an natural equivalence $\pi_i(\sfV^kM)
\simeq V^k(\pi_i M)$, for $M \in \calD^b_{\reghol}(\scrD_{\C,0})$,
where $V^k$ denote the $k$th piece of the ordinary Kashiwara-Malgrange
filtration on the discrete $\scrD_{\C,0}$-module $\pi_i M$. It follows from
the properties of the ordinary $V$-filtration that
$\pi_i(\sfV^k M)$ is a finite type $\C\{t\}$-module. Moreover, t-exactness
also implies that $\sfV^k M \in \Mod_{V^0\scrD_{\C,0}}^b$. Therefore,
noting that $\C\{t\}$ is a regular ring, we find that
$\sfV^k M \in \Mod_{\C\{t\}}^{\perf}$ under the forgetful 
functor $\Mod_{V^0\scrD_{\C,0}}^b \to \Mod_{\C\{t\}}^b$.

\item We show that, each $k \geq 0$, the morphism 
$t \cdot: \sfV^k M \to \sfV^{k+1} M$ induced isomorphisms
of $V^0\scrD_{\C,0}$-modules on all homotopy groups.
Let $i \in \bbZ$, and consider the induced map,
$\pi_i(\sfV^k M) \to \pi_i(\sfV_{k+1}M)$. Since $\sfV^n$ is
t-exact for any $n \in \bbZ$, this is the map
$t \cdot: V^k(\pi_i M) \to V^{k+1}(\pi_i M)$, which is
an isomorphism by the properties of the ordinary Kashiwara-Malgrange
filtration.

\item This is proven in the exact same manner as the
item above, except we note additionally that $\pi_i$
commutes with finite coproducts.

\item We show that $\gr_{\sfV}^k M$ is a perfect
complex of $\C$-vector spaces in the same manner as
item (i), using that $\gr_{\sfV}^k$ is exact. Since $\gr_{\sfV}^k M$
is perfect, it has finitely many nonzero homotopy groups. Let $I \subset \bbZ$
be the collection of indices $i$ such that $\pi_i(\gr_{\sfV}^k M) \not\simeq 0$.
By the properties of the ordinary Kashiwara-Malgrange
filtration, for each $i \in I$ there exists a nonzero polynomial
$b_i(s)$ with roots in $\bbQ \cap [0,1)$ such 
that $b_i(t \partial_t - k)$ vanishes identically
on $\gr_V^k(\pi_i M)$. Let $b(s) := \prod_{i \in I} b_i(s)$. 
Since $\gr_{\sfV}^k$ is t-exact,
$\gr_V^k(\pi_i M) \simeq \pi_i(\gr_{\sfV}^k M)$,
so we obtain that the morphism,
\[b(s) \cdot: \gr_{\sfV}^k M \to \gr_{\sfV}^k M,\]
induces the zero map on all homotopy groups.
Thus, the action of $b(s)$ on $\gr_{\sfV}^k M$ is
nullhomotopic. The choice of nullhomotopy furnishes a
factorization as desired.
\end{enumerate}
\eproof

Moreover, the properties enumerated in \cref{prop: sfV properties}
uniquely characterize the filtered derived $\scrD_{\C,0}$-module 
$\sfV^{\bullet}M$.

\bprop
\label{prop: uniqueness of derived KM filtration}
Suppose that $M \in \calD^b_{\reghol}(\scrD_{\C,0})$, and
suppose that $\sfU^{\bullet} M$ and $\witsfU^{\bullet} M$ are
two objects in $\Mod_{V^{\bullet}\scrD_{\C,0}}(\Fil(\Sp))$
whose underlying derived $\scrD_{\C,0}$-module is $M$.
Suppose additionally that each of $\sfU^{\bullet}M$ and
$\witsfU^{\bullet}M$ satisfies the enumerated properties 
in \cref{prop: sfV properties}. Then
there is an equivalence $\sfU^{\bullet}M \xrightarrow{\simeq} 
\witsfU^{\bullet}M$.
\eprop

\bproof
We begin by producing, for each $k$, a factorization 
\beqn
\label{factorization triangle}
\begin{tikzcd}
	& {\witsfU^jM} \\
	{\sfU^kM} && {M}
	\arrow["i_k", from=2-1, to=2-3]
	\arrow["\imath_j", from=1-2, to=2-3]
	\arrow[dashed, from=2-1, to=1-2],
\end{tikzcd}
\eeqn
for some $j_k \leq k$. By assumption, $\sfU^kM$
is a perfect $\C\{t\}$-module. By \cite[Proposition 7.2.4.11]{HA}, any object
$M \in \Mod_{\C\{t\}}^{\perf}$ can be obtained by as a colimit over a finite 
diagram $P: I \to \Mod_{\C\{t\}}$ such that
each $P_i$ is a finitely generated free $\C\{t\}$-module.
Since we have a map $\sfU^kM \to M$, we also obtain
maps $P_i \to M$ that exhibit $M$ as a cocone over $P$. 
Since $P_i$ is finitely generated free,
say of rank $n_i$, such a map is determined by an $n_i$-tuple
$(x_1, \dots, x_{n_i})$ of maps $x_{\ell}: \ast \to \pi_0(M)$. Since
$M \simeq \varinjlim_{j} \witsfU^jM$, the image
of $x_{\ell}$ lies in the image of $\pi_0(\witsfU^{j_{\ell}}M)
\to \pi_0(M)$ for some $j_{\ell}$; hence $x_{\ell}$ factors through
$\pi_0(\witsfU^{j_{\ell}}M) \to \pi_0(M)$. Set $j_k := \max_{\ell} j_{\ell}$.
Then by the universal property of free modules, $P_i \to M$
factors through $\witsfU^{j_k}M$. It is easy to check that this
exhibits $\witsfU^{j_k}M$ as a cocone over $P$, so by the
universal property, we obtain a map $\sfU^kM \to \witsfU^{j_k}M$
that factorizes $\sfU^kM \to M$.

We proceed to show that there exists an $N \geq 0$ such that
we can take $j_k$ to be $k - N$ for \emph{any} $k$. In fact, we may
take $N = \max(-j_0, -(1+j_1))$. Indeed, if $k \geq 0$, $\sfU^kM \to M$
factors through $\witsfU^{k-j_0}M \to M$, and, if $k \leq 0$,
$\sfU^kM \to M$ factors through $\witsfU^{k-j_1}M \to M$,
using the equivalences furnished by properties (\labelcref{II}) and
(\labelcref{III}).

Continuing, we use the existence of the $b$-functions and factorizations
guaranteed by property (\labelcref{iv}).
Indeed, let $b(s)$ and $\wit{b}(s)$ be the functions, and
\[\begin{tikzcd}
	& {\sfU^{k+1}M} &&&& {\witsfU^{k+1}M} \\
	{\sfU^k\F} && M && {\witsfU^kM} && M
	\arrow["{\wit{b}(t\partial_t-k)}", from=2-5, to=2-7]
	\arrow["{b(t\partial_t - k)}", from=2-1, to=2-3]
	\arrow[dashed, from=2-5, to=1-6]
	\arrow["{\imath_{k+1}}", from=1-6, to=2-7]
	\arrow[dashed, from=2-1, to=1-2]
	\arrow["{i_{k+1}}", from=1-2, to=2-3]
\end{tikzcd}\]
the factorizations, for $\sfU^{\bullet}$
and $\witsfU^{\bullet}$, respectively, guaranteed by property (\labelcref{IV}). Combining 
the left-hand diagram above with \labelcref{factorization triangle} 
we obtain the commutative diagram,
\[\begin{tikzcd}
	& {\sfU^{k+1}M} && {\witsfU^{k-N+1}M} \\
	{\sfU^kM} && {M}
	\arrow["{b(t\partial_t-k) \cdot}", from=2-1, to=2-3]
	\arrow["{(1)}", dashed, from=2-1, to=1-2]
	\arrow["i_{k+1}", from=1-2, to=2-3]
	\arrow["\imath_{k-N+1}", from=1-4, to=2-3]
	\arrow["{(2)}", dashed, from=1-2, to=1-4].
\end{tikzcd}\]
For similar reasons, we also obtain the commutative
diagram,
\[\begin{tikzcd}
	& {\witsfU^{k-N}M} && {\witsfU^{k-N+1}M} \\
	{\sfU^kM} && {M}
	\arrow["{\wit{b}(t\partial_t-k+N) \cdot}", from=2-1, to=2-3]
	\arrow["{(1')}", dashed, from=2-1, to=1-2]
	\arrow["{\wit{b}(t\partial_t-k+N) \cdot}", from=1-2, to=2-3]
	\arrow["\imath_{k-N+1}", from=1-4, to=2-3]
	\arrow["{(2')}", dashed, from=1-2, to=1-4].
\end{tikzcd}\]
Since $b(s-k)$ and $\wit{b}(s-k+N)$ have no
common root, there exist $p(s), q(s) \in \C[s]$ such that $p(s)b(s-k) +
q(s)\wit{b}(s-k+N) = 1$ by B\'ezout's lemma. 
Hence, using the diagrams above,
we obtain a chain of homotopies, 
\begin{align*}
i_k 	&= p(t\partial_t)b(t\partial_t-k) + q(t\partial_t)\wit{b}(t\partial_t-k+N) \cdot \\
	&\simeq {p(t\partial_t) \cdot} \circ \imath_{k-N+1} \circ (2) \circ (1) + {q(t\partial_t) \cdot} \circ \imath_{k-N+1} \circ (2') \circ (1') \\
	&\simeq \imath_{k-N+1} \circ \left({p(t\partial_t) \cdot} \circ (2) \circ (1) + {q(t\partial_t) \cdot} \circ (2') \circ (1')\right),
\end{align*}
whose composition gives a factorization
$i \simeq \imath \circ F$, where we have
let $F := {p(t\partial_t) \cdot} \circ (2) \circ (1) + 
{q(t\partial_t) \cdot} \circ (2') \circ (1')$.
Repeating this process, we obtain a map
$\sfU^kM \to \witsfU^kM$ and a factorization,
\[\begin{tikzcd}
	& {\witsfU^{k}M} \\
	{\sfU^kM} && M
	\arrow["i_k", from=2-1, to=2-3]
	\arrow[dashed, from=2-1, to=1-2]
	\arrow["\imath_k", from=1-2, to=2-3].
\end{tikzcd}\]

Tracing through the constructions, 
one sees that the maps obtained in this way
are natural in the index $k$, so it obtains
a map of filtered objects,
$\sfU^{\bullet}M \to \witsfU^{\bullet}M$.
In order to show that this map is an equivalence,
we shows that it induces an equivalence on associated
graded objects. Let $n \in \bbZ$. By the exactness
of $\gr^n$, the induced map on homotopy 
groups, for each $i \in \bbZ$,
$\pi_i \gr_{\sfU}^nM \to \pi_i \gr_{\witsfU}^nM$
equivalently gives a map,
\beqn
\label{eqn: 1}
\gr_U^n(\pi_i M) \to \gr_{\wit{U}}^n(\pi_i M),
\eeqn
where $U^{\bullet}$ and $\wit{U}^{\bullet}$ are the filtrations on 
$\pi_i M$ induced by $\sfU^{\bullet}$ and $\witsfU^{\bullet}$. 
It is easy to see that both $U^{\bullet}$ and $\wit{U}^{\bullet}$
satisfy the enumerated properties
of the usual $V$-filtration in
\cref{def: usual KM filtration}. It is also clear that
the application of $\pi_i(\gr^n(-))$ to the objects
and maps in each step of the above proof 
yields an argument identical to Kashiwara's original proof of the
uniqueness of the $V$-filtration
(c.f. \cite[Theorem 1.1]{Kashiwara83}).
It follows that \labelcref{eqn: 1} is an equivalence,
which completes the proof.
\eproof

\section{$\wE_{\C,0}$-modules}
\label{sec: wE_{C,0}-modules}

Heuristically, the sheaf of formal 
microdifferential operators $\wE_X$
on a complex manifold $X$ is
a sheaf on $T^*X$ containing differential
operators on $X$, such that if the
principal symbol of a differential operator
$P$ is invertible on an open set
$U \subset T^*X$, then $P$
itself is invertible in $\wE_X$.
By construction there is a canonical
map of sheaves,
\[\pi^*\scrD_X \to \wE_X,\]
where $\pi: T^*X \to X$ denotes the
canonical projection.
When $X = \bfAone$, taking the stalk at $(0;1) \in T^*\bfAone$
induces a canonical inclusion of rings,
\beqn
\label{eqn: canonical inclusion of ring}
\scrD_{\C,0} \hook \wE_{\C,0},
\eeqn
where $\wE_{\C,0}$ denotes the ring 
of germs at $(0;1)$ of formal
microdifferential operators on $\C$.
The inclusion $\scrD_{\C,0} \hook \wE_{\C,0}$
makes $\wE_{\C,0}$ into both a left and 
a right $\scrD_{\C,0}$-module.

\bdef
\label{def: formal microlocalization}
We let $\wih{(-)}: \calD_{\reghol}(\scrD_{\C,0}) 
\to \calD_{\reghol}(\wE_{\C,0})$
denote the functor given by $N \mapsto \wE_{\C,0} 
\otimes_{\scrD_{\C,0}} N$. We will call $\wih{N}$ the
\emph{formal microlocalization} of $N$.
\edefn

\brem
The inclusion $\scrD_{\C,0} \hook \wE_{\C,0}$ is a flat map of rings.
\erem

Concretely, $\wE_{\C,0}$ is given as a $\C$-vector space by
the subspace of $\C\{t\}\llp \xi \rrp$ spanned by elements of
form $\sum_{k \geq k_0} a_k(t)\xi^k$, where
there is some fixed $r>0$ such that the power
series $a_k(t)$ all have radius of convergence
greater than or equal to $r$. That is, elements
of $\wE_{\C,0}$ are Laurent series in the formal variable $\xi$ whose
coefficients all converge on some fixed
disc of radius $r$, independent of $k$. The product structure
on $\wE_{\C,0}$ is not the standard commutative one,
but rather is noncommutative. It is written explicitly, for example,
in \cite[\S1.2]{Schapira}.
Using this concrete description of $\wE_{\C,0}$, 
the canonical inclusion $\scrD_{\C,0} \hook \wE_{\C,0}$ is given
by sending $\partial_t \mapsto \xi^{-1}$.

\brem
We recommend either \cite{Schapira}
or \cite{microcalc} as introductions to the
theory of microdifferential operators and microdifferential systems.
\erem

\subsubsection{Microdifferential operators with polynomial coefficients}
Elements of $\C\{t\}$ that are polynomials trivially converge on a disc
of any positive radius, so elements of $\C\{t\}\llp \xi \rrp$ whose coefficients 
are all polynomials in $t$ form a subring of $\wE_{\C,0}$.

\bdef
\label{def: polynomial microdifferential operators}
Let $\wE_{\C,0}^{\poly}$ denote the subring of
$\wE_{\C,0}$ spanned by elements of $\C\{t\}\llp \xi \rrp$
whose coefficients are polynomials in $t$.
\edefn

The inclusion of rings $\wE_{\C,0}^{\poly} \hook \wE_{\C,0}$
induces a canonical restriction functor,
\[\res: \Mod_{\wE_{\C,0}} \to \Mod_{\wE_{\C,0}^{\poly}},\]
as well as a functor on the hearts, $\res: \calA_{\wE_{\C,0}} \to \calA_{\wE_{\C,0}}$.

\subsubsection{Regular holonomic modules}
Just as for $\scrD_X$-modules, there is a notion of regular 
holonomic $\wE_X$-module. We will not recall the definition,
which can be found in \cite{microcalc} as Definition 8.27.

Regular holonomic $\wE_{\C,0}$-modules
are defined as those objects in $\calA_{\wE_{\C,0}}$
which are germs of regular holonomic $\wE_{\C}$-modules
defined in a neighborhood of $(0;1) \in T^*\C$.
We denote by $\RegHol_{\wE_{\C,0}}$ the full 
abelian subcategory of $\calA_{\wE_{\C,0}}$
spanned by such objects.

We define $\calD^b_{\reghol}(\wE_{\C,0})$ to be the
full subcategory of $\calD^b(\wE_{\C,0})$ on objects
whose cohomology modules are regular holonomic. 

%
%

The abelian category of regular holonomic $\wE_{\C,0}$-modules
and the formal microlocalization functor admit the following convenient
description. It is well known (see e.g. \cite[Theorem 5.3.12]{Dimca}) 
that $\RegHol_{\scrD_{\C,0}}$ 
is equivalent to the category whose objects are diagrams
\[\begin{tikzcd}
	E & F
	\arrow["c", shift left, from=1-1, to=1-2]
	\arrow["v", shift left, from=1-2, to=1-1],
\end{tikzcd}\]
where $E$ and $F$ are finite dimensional $\C$-vector spaces
and $c \circ v + \id_F$ and $v \circ c + \id_E$ are isomorphisms.
Similarly, $\RegHol_{\wE_{\C,0}}$ admits a description
as the category of finite dimensional $\C$-vector spaces $V$
equipped with an automorphism $T: V \to V$.
In terms of these descriptions, the functor of formal microlocalization $\wih{(-)}:
\RegHol_{\scrD_{\C,0}} \to \RegHol_{\wE_{\C,0}}$ is given
objects by the assignment,
\[\begin{tikzcd}
	{(E,F, c,v)} & {(F, c \circ v)}
	\arrow[maps to, from=1-1, to=1-2].
\end{tikzcd}\]

Alternatively, $\RegHol_{\scrD_{\C,0}}$ is equivalent to 
the category $\Perv(\C;0)$ of germs of perverse sheaves on $\C$ constructible
with respect to stratification consisting of $\{0\} \subset \C$
and its complement, and $\RegHol_{\wE_{\C,0}}$ is equivalent
to the category of germs of local systems on $\C\setminus 0$.
In these terms, the formal microlocalization is given
by taking the vanishing cycles of an object in $\Perv(\C;0)$
equipped with its monodromy automorphism. 

Either of the above descriptions of $\wih{(-)}: \RegHol_{\scrD_{\C,0}}
\to \RegHol_{\wE_{\C,0}}$ makes clear that it is essentially surjective,
which allows us to utilize the following working definition of a regular
holonomic $\wE_{\C,0}$-module. 

\bdef
\label{def: regular holonomic E-module}
An $\wE_{\C,0}$-module $M \in \calA_{\wE_{\C,0}}$
is \emph{regular holonomic} if there exists an object
$N \in \RegHol_{\scrD_{\C,0}}$ such that
$M = \wih{N}$.
\edefn

\brem
It follows from \cref{lem: modified Beilinson}
that $\calD^b(\RegHol_{\wE_{\C,0}}) 
\xrightarrow{\simeq} \calD^b_{\reghol}(\wE_{\C,0})$.
\erem

\brem
As pointed out to the author by Pavel Safronov, the
functor taking an object of $\Perv(\C;0)$ to its
\emph{nearby cycles} equipped with the monodromy
automorphism corresponds to the functor in Deligne's
correspondence if we restrict ourselves to considering
the image of $\Conn^{\reg}(\C;0)$ inside $\Perv(\C;0)$.
\erem

We will cheat and define the notion of regular holonomic
$\wE_{\C,0}^{\poly}$-module as follows.

\bdef
$\RegHol_{\wE_{\C,0}^{\poly}}$ is defined to be the essential
image of $\RegHol_{\wE_{\C,0}}$ under the functor $\res$.
The $\infty$-category $\calD^b_{\reghol}(\wE_{\C,0}^{\poly})$
is defined to be the full subcategory of $\calD^b(\wE_{\C,0}^{\poly})$ spanned
by objects whose cohomology modules are regular holonomic.
\edefn

\bdef
\label{def: muRH}
We denote by $\muRH: \Fun(S^1, \Perf) \to \calD^b_{\reghol}(\wE_{\C,0}^{\poly})$ the functor 
obtained as the composition of the
following chain of functors,
\[\Fun(S^1, \Perf) \xrightarrow{\Temp} 
	\calD^b_{\reghol}(\scrD_{\C,0}) \xrightarrow{\Four}
		\calD^b_{\reghol}(\scrD_{\C,0}) \xrightarrow{\wE_{\C,0} \otimes -} 
			\calD^b_{\reghol}(\wE_{\C,0}). \xrightarrow{\res}
				\calD^b_{\reghol}(\wE_{\C,0}^{\poly}).
\]
\edefn

\bprop
The restriction of $\muRH$ to $\Loc(\bfAone \setminus 0)$,
viewed as a functor $\Loc(\bfAone \setminus 0) \to \RegHol_{\wE_{\C,0}^{\poly}}$,
is naturally isomorphic to the functor $\invRH$
defined in \cref{def: inverse RH}.
\eprop

\bproof
Let $(V,T) \in \Loc(\bfAone \setminus 0)$ be a pair
of a finite dimensional vector space
and an automorphism representing a
local system on $\bfAone \setminus 0$. By \cref{rem:
calculation of RH}, the image of $(V,T)$
under Deligne's correspondence is the
$\scrD_{\C,0}$-module given by
$\C\{t\} \otimes_{\C} V$ with the
action of $\partial_t$ given by
$\partial_t \cdot (f \otimes v) := \partial_t f \otimes v
+ \frac{f}{t} \otimes \frac{\log(T)v}{-2\pi i}$.
Let us denote this $\scrD_{\C,0}$-module
by $M_{V,T}$.
The Fourier transform of $M_{V,T}$ is
computed, essentially by definition, as follows. 
First, consider $M_{V,T}$
as a $\C[t]\langle \partial_t\rangle$-module
by restricting along the ring inclusion
$\C[t]\langle \partial_t\rangle
\hook \C\{t\}\langle \partial_t\rangle$,
which we denote by $M^{\ndalg}_{V,T}$.
Compute the Fourier transform of $M^{\ndalg}_{V,T}$
as an $A_t$-module to obtain an
$A_{\eta}$-module, $\Four(M^{\ndalg}_{V,T})$.
Now $\Four(M_{V,T}) = \C\{\eta\}\langle \partial_{\eta} \rangle 
\otimes_{\C[\eta]\langle \partial_{\eta} \rangle} \Four(M^{\ndalg}_{V,T})$.
The object $\Four \circ \Temp(M) \in \calA_{\C\{\eta\}((\partial_{\eta}^{-1}))}$ is therefore 
equivalent to $\C\{\eta\}((\partial_{\eta}^{-1})) 
\otimes_{\C[\eta]\langle \partial_{\eta} \rangle} \Four(M^{\ndalg}_{V,T})$.
Unraveling the definitions, it is clear that
this the restriction of the module structure of the latter object 
to the subring $\C[\eta]((\partial_{\eta}^{-1})) \subset \C\{\eta\}((\partial_{\eta}^{-1}))$ 
is precisely the $\C\llp u \rrp$-module
$V\llp u \rrp$ with connection $\nabla = d + \frac{\log(T)}{-2\pi i u}du$,
once we set $u := \partial_{\eta}^{-1}$.
\eproof

Recall from \cref{prop: descriptions of D-modules} that the
abelian category of regular holonomic $\scrD_{\C,0}$-modules
is equivalent to a certain subcategory of the category of regular 
holonomic modules over the Weyl algebra. In other words,
when a $\scrD_{\C,0}$-module is regular holonomic it
is actually determined by its module structure over the much
smaller subring of polynomial differential operators inside $\scrD_{\C,0}$. 
This observation motivates the following lemma.

\blem
\label{lem: poly}
The functor $\res$ induces an equivalence
of categories,
\[\RegHol_{\wE_{\C,0}} \xrightarrow[\res]{\simeq} \RegHol_{\wE_{\C,0}^{\poly}}.\]
\elem

\bproof[Proof sketch.]
The lemma can be deduced by a GAGA argument
similar to the one showing the equivalence between
$\RegHol_{\scrD_{\C,0}}$ and the abelian
category of holonomic $A_t$-modules $M$
regular at $0$ and at infinity such that $M|_{\Aone \setminus 0}$
of \cref{prop: descriptions of D-modules}. Details
of the argument showing the latter are found
in \cite[\S I.4]{MalgrangeTextbook}. 
\eproof

\brem
In light of \cref{lem: poly}, 
we regard $\muRH$ as a functor from $\calD^b_{\reghol}(\scrD_{\C,0})$
$\calD^b_{\reghol}(\wE_{\C,0})$. It follows that $\muRH$
may be defined alternatively as $\muRH := (\wE_{\C,0} \otimes -) \circ \Four \circ \Temp$.
\erem

\subsection{Vanishing cycles and microlocalization}

\bdef
\label{def: mu}
Let $\mu: \calD^b_{\reghol}(\scrD_{\C,0}) \to 
\calD^b_{\reghol}(\wE_{\C,0})$
denote the composition
$\muRH \circ \varphi_{\scrD}$.
\edefn

The following proposition can be viewed as
describing the relationship between vanishing
cycles and the formal microlocalization.
It is the derived version of a result appearing
in \cite{Sabbah}.

\bprop
\label{prop: mu = wih}
There is a natural isomorphism,
$\wih{(-)} \simeq \mu(-)$.
\eprop

\bproof
By \cite[Corollary 3.6]{Sabbah}, we have a natural
isomorphism, 
\beqn
\label{eqn: mu=wih}
\wih{(-)}|_{\RegHol_{\scrD_{\C,0}}}
\simeq \mu|_{\RegHol_{\scrD_{\C,0}}}.
\eeqn
It follows that, since 
$\wih{(-)}|_{\RegHol_{\scrD_{\C,0}}}$ is right exact,
$\mu|_{\RegHol_{\scrD_{\C,0}}}$ is too.
Thus, $\mu$ is the left derived functor of its
restriction $\mu|_{\RegHol_{\scrD_{\C,0}}}$
by \cite[Theorem 1.3.3.2]{HA}.
and the natural isomorphism \labelcref{eqn: mu=wih}
induces a natural equivalence of derived functors,
$\wih{(-)} = \mu$. 
\eproof

\subsection{Monoidal structure on $\wE_{\C,0}$-modules}
Although $\wE_{\C,0}$ is a noncommutative ring,
there is a canonical symmetric monoidal structure on
$\calD^b(\wE_{\C,0})$ in the same way that there
is a canonical symmetric monoidal structure on
$\scrD_{\C,0}$-modules.

\begin{exmp}
The tensor product of two $\C\llp u \rrp$-modules with
connection\footnotemark 
\footnotetext{See \cref{rem: going from u-connections to micromodules}.}
in $\calA_{\wE_{\C,0}^{\poly}}$ is easy to write down explicitly.
Let $\calM = (V, \nabla)$ and $\calM' = (V', \nabla')$ denote
two such objects, which we view as $\C\llp u \rrp$-modules
$V, V'$ equipped with $u$-connections $\nabla, \nabla'$.
Then 
\[\calM \otimes \calM' = (V \otimes_{\C\llp u \rrp} V', \nabla + \nabla').\]
\end{exmp}

\begin{notation}
We denote the symmetric monoidal product on
$\calD^b(\wE_{\C,0})$ by $\otimes$.

We denote the restricted symmetric monoidal product on
$\calD^b(\wE_{\C,0}^{\poly})$ by $\otimes_{\C\llp u \rrp}$
in order to emphasize that
$\C\llp u \rrp$ with the trivial $\wE_{\C,0}^{\poly}$-module
structure is the monoidal unit for this structure and 
to more closely match the notation in \cite{Sabbah}.
\end{notation}

\brem
The subcategory $\calD^b_{\reghol}(\wE_{\C,0})$ is stable
under the monoidal product $\otimes$.
\erem

\section{Sheaves of $\scrD_{\C,0}$-modules}
\label{sec: Sheaves of D_{C,0}-modules}

Constructible sheaves of $\scrD_{\C,0}$-modules
play a key role in Sabbah's proof of Kontsevich's
conjecture. In this section, we collect some results
about these objects in the $\infty$-categorical setting
that will allow us to formulate our main theorem.

\subsection{Constructible sheaves of $\scrD_{\C,0}$-modules}
The subcategory $\calD^b_{\reghol}(\scrD_{\C,0}) \subset
\calD(\scrD_{\C,0})$ determines a coefficient pair in the sense
of \cref{def: coefficient pair}. As such, given a complex analytic
space $X$, there is the associated category 
$\Shv_c(X; \calD(\scrD_{\C,0}))$ of constructible 
sheaves on $X$ with respect to $(\calD(\scrD_{\C,0}), 
\calD^b_{\reghol}(\scrD_{\C,0}))$ as in \cref{def: constructible sheaves}.
Likewise, $(\calD(\wE_{\C,0}), \calD^b_{\reghol}(\wE_{\C,0}))$
is a coefficient pair, and we might consider the category 
$\Shv_c(X; \calD(\wE_{\C,0}))$.

\begin{notn}
We will use $\Shv_c(X; \scrD_{\C,0})$ and $\Shv_c(X; \wE_{\C,0})$,
respectively, to denote $\Shv_c(X; \calD(\scrD_{\C,0}))$ 
and $\Shv_c(X; \calD(\wE_{\C,0}))$, respectively.
\end{notn}

\brem
The triangulated category $\pi_0(\Shv_c(X; \scrD_{\C,0}))$
is equivalent to $D^b_{\C\c, \rh}(X;\scrD_{\C,0})$ the bounded derived
category of sheaves of $\scrD_{\C,0}$-modules
whose cohomology sheaves are constructible
sheaves of regular holonomic $\scrD_{\C,0}$-modules,
defined in \cite[\S3.a]{Sabbah}.
\erem

\begin{exmp}
\label{example1}
Let $\F \in \Shv_c(X; {\scrD_{\C,0}})$.
Although the stalks of $\F$ are bounded with regular holonomic
cohomology, for a general open $U$, it is not the case that
$\F(U) \in \calD^b_{\reghol}(\scrD_{\C,0})$.
Indeed, consider the complex analytic space, 
$\C \setminus \bbZ$. We have that $\bbH^1(\C\setminus\bbZ; \un{\bbZ})
\simeq H^1(\C\setminus\bbZ; \bbZ) \simeq \bbZ^{\oplus_{i \in \bbZ}}$. 
Since $\scrD_{\C,0}$ is torsion-free,
it is flat over $\bbZ$, so this implies that $\bbH^1(\C\setminus\bbZ; 
\un{\scrD_{\C,0}}) \simeq \scrD_{\C,0}^{\oplus_{i \in \bbZ}}$ 
Thus the derived global sections of the constant sheaf 
$\un{\scrD_{\C,0}}_{\C \setminus \bbZ}$
is not coherent,\footnotemark so in particular not
regular holonomic.
\footnotetext{Recall that a $\scrD_{\C,0}$-module $M$ is
regular holonomic if and only if there exists a $\scrD_{\bfAone}$-module
$\wit{M}$ whose germ at $0$ is $M$. 
If a $\scrD_{\bfAone}$-module $\wit{M}$
is coherent, there exists a neighborhood
$U$ around each point $z \in \C$  such 
that the restriction $\wit{M}|_U$ admits 
a surjection ${\scrD_{\bfAone}^n}|_U
\to \wit{M}|_U$ for some natural number $n$. This clearly implies that the
germ of $\wit{M}$ at any point must be finitely generated.}
\end{exmp}

In light of \cref{example1}, we see that the category
$\Shv_c(X; {\scrD_{\C,0}})$ is not equivalent to
$\Shv_c(X; {\scrD_{\C,0}^{b,\reghol}})$; i.e. the
sections assigned to any open $U \subset X$ do not
necessarily have regular holonomic cohomology.
As such, a limit-preserving functor $F: \calD^b_{\reghol}(\scrD_{\C,0})
\to \calD^b_{\reghol}(\wE_{\C,0})$
does not naively define a functor,
\[\Shv_c(X; {\scrD_{\C,0}}) \to \Shv_c(X; {\wE_{\C,0}})\]
via composition of functors.

Observe, however, that, as a complex analytic space, $X$
admits a basis of open subsets which are homotopy
equivalent to finite CW complexes (finite homotopy type). For each such
open $U$ and each constructible sheaf
$\F \in \Shv_c(X; {\scrD_{\C,0}})$, $\F(U)$ in
fact \emph{is} an object of $\calD^b_{\reghol}(\scrD_{\C,0})$.
This observation suggests that, if $X$ admits a basis
of opens (including itself) that have finite homotopy
type, we might instead consider working with sheaves defined
on this basis.

In \cref{sec: Appendix A}, we have collected a
number of results that allow us to obtain functors
on $\Shv_c(X; \scrD_{\C,0})$ from functors on
$\calD_{\reghol}(\scrD_{\C,0})$ when $X$ has
finite homotopy type. In particular, \cref{prop: basis shv = real shv}
shows that, by restriction, working with basis sheaves on $X$, as suggested
above, is equivalent to working with constructible sheaves
on the entirety of $\calU(X)$.
In the examples that follow, 
let $X$ and $\calB(X)$ be as in
\cref{conv: finite homotopy type}.

\begin{exmp}
\label{exmp: mu_shv = wih_shv}
By \cref{cor: induces shv functor},
the functors $\mu, \wih{(-)}: \calD^b_{\reghol}(\scrD_{\C,0}) 
\to \calD^b_{\reghol}(\wE_{\C,0})$
induce functors,
\[\mu_{\Shv}, \wih{(-)}_{\Shv}: \Shv_c^{f.s.}(X; {\scrD_{\C,0}}) 
	\to \Shv_c^{f.s.}(X; {\wE_{\C,0}}),\]
which are naturally isomorphic by
\cref{prop: mu = wih} and \cref{cor: induces natural iso}.
\end{exmp}

\begin{exmp}
The functor $\muRH: \Fun(S^1, \Vect) \to \calD^b_{\reghol}(\wE_{\C,0})$
similarly defines a functor,
\[\muRH_{\Shv}: \Shv_c^{f.s.}(X; \Fun(S^1, \Vect)) \to \Shv_c^{f.s.}(X; \wE_{\C,0}).\]
Note that, using the equivalence \labelcref{eqn: S1 commutes with Shv},
we may also regard $\muRH_{\Shv}$ as a functor $\Fun(S^1, \Shv_c^{f.s.}(X; \C)) \to
\Shv_c^{f.s.}(X; \wE_{\C,0})$.
\end{exmp}

\begin{exmp}
\label{exmp: V-filtration on sheaves}
Using \cref{cor: non-c induces shv functor}, the derived functor of the $V$-filtration
on regular holonomic $\scrD_{\C,0}$-modules, $\sfV^{\bullet}$, induces a functor,
\[\sfV^{\bullet}: \Shv_c^{f.s.}(X; \scrD_{\C,0}) \to \Shv(X; \Mod_{V^{\bullet}\scrD_{\C,0}}(\Fil(\Sp))).\]
Sheaves valued in $\Mod_{V^{\bullet}\scrD_{\C,0}}(\Fil(\Sp))$ are equivalently
left modules in $\Fil(\Sp)$-valued sheaves over the constant
sheaf, $\un{V^{\bullet}\scrD_{\C,0}}$.
At the same time, $\Shv(X; \Fil(\Sp))
\simeq \Fil(\Shv(X; \Sp))$ by \cref{cor: Fil commutes with Shv}.
These facts taken all together,
we may regard $\sfV^{\bullet}$ as 
a functor $\Shv_c^{f.s.}(X; \scrD_{\C,0}) \to
\Mod_{V^{\bullet}\un{\scrD_{\C,0}}}(\Fil(\Shv(X; \Sp))$,
where $V^{\bullet}\un{\scrD_{\C,0}}$ denotes the
filtration on the constant sheaf induced by the
standard $V$-filtration
on $\scrD_{\C,0}$.
\end{exmp}

\subsection{The $V$-filtration on sheaves}
In \cref{exmp: V-filtration on sheaves}, we built a filtration on 
constructible sheaves $\Shv_c^{f.s.}(X; \scrD_{\C,0})$ using the
$V$-filtration on objects in $\calD^b_{\reghol}(\scrD_{\C,0})$.
Just as the classical $V$-filtration is unique,
we might hope that the $\sfV^{\bullet}$-filtration on constructible
sheaves is unique, as well. In this section, we show that this
is indeed the case by imitating Kashiwara's original proof
of the uniqueness of the $V^{\bullet}$-filtration on regular holonomic
$\scrD$-modules.

Kashiwara's proof relies on the structure maps in the $V^{\bullet}$-filtration
being inclusions, so in order to imitate his proof, we establish the following
lemma showing that the structure morphisms of the $\sfV^{\bullet}$-filtration
are monomorphisms.

\bprop
\label{prop: uniqueness of sfV-filtration}
Let $\F \in \Shv_c^{f.s.}(X; \scrD_{\C,0})$, and let $X \to A$
denote a finite Whitney stratification with respect to which $\F$ is
constructible. Suppose that $\sfU^{\bullet}\F$ and $\witsfU^{\bullet}\F$
are two objects of $\Mod_{V^{\bullet}\un{\scrD_{\C,0}}}\Fil(\Shv^A(X;\Sp))$
whose underlying $\un{\scrD_{\C,0}}$-module each
is $\F$, that each satisfy the following properties:
\begin{enumerate}[(i)]
	\item \label{i} For each $a \in A$, $\sfU^k \F|_{X_a}$ is a locally constant sheaf
	of perfect $\C\{t\}$-modules.
	\item \label{ii} The map, $t \cdot: \sfU^k \F \to \sfU^{k+1} \F$ is an equivalence for each $k \geq 0$.
	\item \label{iii} The induced map $\sum_{j=0}^{-k-1} {\partial_t}^j \sfU^{-1} \F \to \sfU^k \F$
	is an equivalence for each $k < 0$.
	\item \label{iv} For each $a \in A$ and each $k \in \bbZ$, $\gr^k_{\sfU}\F|_{X_a}$ is a locally constant sheaf of perfect
	$\C$-modules, and, locally on $X$, there exists a nonzero polynomial $b(s) \in \C[s]$ 
	whose roots have real parts in $[0,1)$ and a factorization,
\[\begin{tikzcd}
	& {\sfU^{k+1}\F} \\
	{\sfU^k\F} && {\F}
	\arrow["{b(t\partial_t-k)}", from=2-1, to=2-3]
	\arrow["{i_{k+1}}", from=1-2, to=2-3]
	\arrow[dashed, from=2-1, to=1-2].
\end{tikzcd}\]
\end{enumerate}
Then there is a natural map $\sfU^{\bullet}\F \to \witsfU^{\bullet}\F$
which is an equivalence of objects in $\Mod_{V^{\bullet}\un{\scrD_{\C,0}}}\Fil(\Shv^A(X;\Sp))$.
\eprop

\bproof
Suppose that $\sfU^{\bullet}\F$ and $\witsfU^{\bullet}\F$ are
two such filtrations of $\F$ satisfying the enumerated properties.
Using recollement for sheaves and induction on the strata, 
we demonstrate that it suffices to produce an equivalence
${\sfU^{\bullet}\F}|_{X_a} \xrightarrow{\simeq} 
{\witsfU^{\bullet}\F}|_{X_a}$, for each $a \in A$.

Without loss of generality, we assume that each $X_a$ is connected.
Let $X_{a_0}$ be one of the closed strata in $X$, and denote by $U_{a_0}
= X \setminus X_{a_0}$ its complement. We denote the inclusion into $X$
of the former by $i$ and of the latter by $j$.
Since $\Mod_{V^{\bullet}\scrD_{\C,0}}(\Fil(\Sp))$ is a closed
symmetric monoidal presentable stable $\infty$-category,
$\Shv(X; \Mod_{V^{\bullet}\scrD_{\C,0}}(\Fil(\Sp)))$ is a recollement of
$\Shv(X_{a_0}; \Mod_{V^{\bullet}\scrD_{\C,0}}(\Fil(\Sp)))$
and $\Shv(U_{a_0}; \Mod_{V^{\bullet}\scrD_{\C,0}}(\Fil(\Sp)))$
in the sense of \cite[Definition A.8.1]{HA} (see \cite[Remark 4.12]{Volpe}).
This means, in particular, that we obtain a fiber sequence,
\[i_*i^*\calG \to \calG \to j_*j^*\calG,\]
for any $\calG \in \Shv(X; \Mod_{V^{\bullet}\scrD_{\C,0}}(\Fil(\Sp)))$,
and that, moreover, a morphism $u$
is an equivalence if and only if $i^*u$ and $j^*u$ are.
Thus, to obtain an equivalence $\sfU^{\bullet}\F \xrightarrow{\simeq} \witsfU^{\bullet}\F$,
it suffices to exhibit equivalences $i^*\sfU^{\bullet}\F \simeq i^*\witsfU^{\bullet}\F$
and $j^*\sfU^{\bullet}\F \simeq j^*\witsfU^{\bullet}\F$.
Repeating this argument for a stratum 
$X_{a_1} \subset X$, contained in and closed in $U_{a_0}$ and its
complement $U_{a_1} := U_{a_0} \setminus X_{a_1}$, we find that
in order to obtain an equivalence, $j^*\sfU^{\bullet}\F \simeq 
j^*\witsfU^{\bullet}\F$, it suffices to exhibit equivalences $(j^*\sfU^{\bullet}\F)|_{X_{a_1}} 
\simeq (j^*\witsfU^{\bullet}\F)|_{X_{a_1}}$ and $(j^*\sfU^{\bullet}\F)|_{U_{a_1}} \simeq 
(j^*\witsfU^{\bullet}\F)|_{U_{a_1}}$. Repeating this process, we see by induction
that to produce an equivalence $\sfU^{\bullet}\F \xrightarrow{\simeq} \witsfU^{\bullet}\F$,
it suffices to exhibit equivalences ${\sfU^{\bullet}\F}|_{X_a} \xrightarrow{\simeq}
{\witsfU^{\bullet}\F}|_{X_a}$ for each $a \in A$.

Since $\F$ is assumed constructible with respect
to $X \to A$, $\F|_{X_a}$ is the constant sheaf given by $\pt^*(\F_x)$,
for any fixed $x \in X_a$.
It therefore suffices to exhibit a factorization,
\[\begin{tikzcd}
	& {\witsfU^{k}\F_x} \\
	{\sfU^k\F_x} && \F_x
	\arrow["i_k", from=2-1, to=2-3]
	\arrow[dashed, from=2-1, to=1-2]
	\arrow["\imath_k", from=1-2, to=2-3],
\end{tikzcd}\]
for each $k$, where the notation $\sfU^k\F_x$
is unambiguous because $(\ev_k \circ \sfU^{\bullet}\F)_x 
\simeq \ev_k \circ (\sfU^{\bullet}\F)_x$, as $\ev_k$ commutes
with filtered colimits. To exhibit such a factorization,
we note that because
$\sfU^{\bullet}\F$ and $\witsfU^{\bullet}\F$
satisfy properties enumerated in the statement
of the proposition above, the stalks 
$\sfU^{\bullet}\F_x, \witsfU^{\bullet}\F)_x \in \Mod_{V^{\bullet}\scrD_{\C,0}}(\Fil(\Sp))$
satisfy the properties for the derived Kashiwara-Malgrange
filtration enumerated in \cref{prop: sfV properties}. The result now follows.
\eproof

We now verify that $\sfV^{\bullet}\F$ satisfies
all the enumerated properties in \cref{prop: uniqueness of
sfV-filtration}.

\bprop
\label{prop: derived KM filtration satisfies properties}
The $V$-filtration on sheaves, $\sfV^{\bullet}$,
satisfies all of the properties enumerated in \cref{prop: uniqueness of sfV-filtration}.
\eprop

\bproof
Let $\F \in \Shv_c^{f.s.}(X; \scrD_{\C,0})$, and let $X \to A$
be a finite Whitney stratification of $X$ with respect to which
$\F$ is constructible. It is clear from its construction
that $\sfV^{\bullet}\F \in 
\Mod_{V^{\bullet}\un{\scrD_{\C,0}}}\Fil(\Shv^A(X; \Sp))$.

We begin by showing that $\sfV^{\bullet}\F$ satisfies
properties (\labelcref{ii}) and (\labelcref{iii}).
Observe that, by \cref{cor: non-c induces shv functor}, 
the following diagram is commutative:
\[\begin{tikzcd}
	{\Shv_c^{f.s.}(X;\scrD_{\C,0})} & {} & {\Shv(X;\scrD_{\C,0})} \\
	{\Shv_c^{f.s.}(\calB(X);\scrD_{\C,0})} & {} & {\Shv(\calB(X);\scrD_{\C,0})}
	\arrow["{\sfV^n}", from=1-1, to=1-3]
	\arrow["{(\ev_n)_{\calB(X)} \circ {\sfV^{\bullet}}_{\calB(X)}}", from=2-1, to=2-3]
	\arrow["\theta_2"', "{\simeq}", from=1-3, to=2-3]
	\arrow["\theta_1"', "{\simeq}", from=1-1, to=2-1],
\end{tikzcd}\]
where we have used the notation therein.
Note that $(\ev_n)_{\calB(X)} \circ {\sfV^{\bullet}}_{\calB(X)} 
\simeq (\ev_n \circ \sfV^{\bullet})_{\calB(X)}$ is given
by pointwise composition with the the $n$th piece of the
$\sfV^{\bullet}$-filtration on $\calD^b_{\reghol}(\scrD_{\C,0})$.
It follows from the properties of the $\sfV^{\bullet}$-filtration 
that $t: ({\sfV^{\bullet}}_{\calB(X)}\theta_1(\F))^k \to ({\sfV^{\bullet}}_{\calB(X)}\theta_1(\F))^{k+1}$ 
is an equivalence for each $k \geq 0$, and that $({\sfV^{\bullet}}_{\calB(X)}\theta_1(\F))^{-k-1} 
\simeq \sum_{j=0}^k \partial_t^j ({\sfV^{\bullet}}_{\calB(X)}\theta_1(\F))^{-1}$ for each $k \geq 0$. 
It now follows that $\sfV^{\bullet}\F$ satisfies properties 
(\labelcref{ii}) and (\labelcref{iii}) by the commutative square.

By \cref{cor: F_Shv compatible with pullback}, $\sfV^{\bullet}$ 
is compatible with restriction along the inclusion of 
a complex analytic subset $Z \subset X$ in the sense that,  
$\sfV^k(\F|_Z) \simeq (\sfV^k\F)|_Z$. Thus, in order to show
that it satisfies properties (\labelcref{i}) and (\labelcref{iv}), 
we may assume that $\F$ is a local system. Without loss of generality,
we may also assume $X$ is connected. In this case, two more
applications of \cref{cor: F_Shv compatible with pullback} yields
an equivalence $\sfV^{\bullet}\F \simeq \pt^*(\sfV^{\bullet}\F_x) \simeq \pt^*(\sfV^{\bullet}(\F_x))$,
for any $x \in X$. But the latter sheaf is the constant sheaf on $X$
with value $\sfV^{\bullet}(\F_x)$, so properties (\labelcref{i}) and (\labelcref{iv}) 
follow from the analogous properties for $\sfV^k(\F_x)$ and
$\gr_{\sfV}^k(\F_x)$ coming from $\sfV^{\bullet}$-filtration on
objects of $\calD^b_{\reghol}(\scrD_{\C,0})$.
\eproof

\section{The Main Theorem}
\label{sec: The Main Theorem}

\subsection{The Kontsevich--Sabbah--Saito theorem}

We recall the set-up of \cite{Sabbah} in order
to state the Kontsevich--Sabbah--Saito theorem. For clarification
or details, we refer the reader to the introduction
of \textit{loc.\ cit.}

\subsubsection{}
Let $X$ be either a smooth complex algebraic variety
or a complex manifold. In each case, we let $\O_X$
denote the structure sheaf of $X$.

\begin{notation}
We use the following notation: given any $\C$-vector space $E$, 
we denote by $E\llb u\rrb$ the $\C\llb u\rrb$-module of formal power 
series in the formal variable $u$ 
with coefficients in $E$ and by $E\llp u \rrp$
the $\C\llp u \rrp$-vector space of formal Laurent
series with coefficient in $E$. For a sheaf $\F$
on $X$, $\F\llp u \rrp$ denotes the sheaf associated 
to the presheaf $U \mapsto \F(U)\llp u \rrp$.
Note that the sheaf $\O_X\llp u \rrp$ is flat over $\O_X$.
\end{notation}

\bdef
\label{def: vs with connection}
Given a $\C\llp u \rrp$-vector space $E$, a $u$-connection 
on $E$ is a map $\nabla: \C\llp u \rrp\partial_u \times E \to E$ which is
$\C\llp u \rrp$-linear in the first factor and satisfies the Leibniz rule
with respect to scalar multiplication by $\C\llp u \rrp$ in the following
sense: if $\ell(u) \in \C\llp u \rrp$ and $v \in E$, then $\nabla_{\partial_u}(v \cdot \ell(u))
= \nabla_{\partial_u}(v) \cdot \ell(u) + v \cdot \partial_u(\ell)$.
We denote the category of $\C\llp u \rrp$-modules with connection
by $\Conn_{\C\llp u \rrp}$.
\edefn

\bdef
More generally, a $u$-connection on a strict complex $(E^{\bullet},d)$
of $\C\llp u \rrp$-vector spaces is a collection of
maps $\nabla: \C\llp u \rrp \partial_u \times E^i \to E^i$,
$\C\llp u \rrp$-linear in the first factor and
satsifying the Leibniz rule with respect to
$\C\llp u \rrp$-multiplication in the second,
such that $[\nabla_{\partial_u}, d] = 
p(u) \cdot d$ for $p(u)$ a unit in $\C\llp u \rrp$
(i.e. $\nabla_{\partial_u}$ preserves cohomology
groups).
\edefn

\brem
\label{rem: going from u-connections to micromodules}
A $\C\llp u \rrp$-module $E$ with connection has a canonical
$\wE_{\C,0}^{\poly}$-module structure given by setting
$\xi \cdot e := eu$ and
$t \cdot e := \nabla_{\partial_u}(e)u^2$,
for $e \in E$. One may check that the action so-defined
satsifies the ring relations for $\wE_{\C,0}^{\poly}$.
Conversely, given a $\wE_{\C,0}^{\poly}$-module $M$,
its underlying $\C\llp \xi \rrp$-module $E$ is a $\C\llp u \rrp$-module
with a canonical connection after relabeling $u := \xi$ and
setting $\nabla_{\partial_u}(e) :=t \cdot e\xi^{-2}$.
There two constructions are obviously inverse to each
other, and given an equivalence of categories,
\[\Conn_{\C\llp u \rrp} \simeq \calA_{\wE_{\C,0}^{\poly}}.\]
\erem

\bdef
\label{def: sheaf with connection}
Given a topological space $Y$
and a sheaf $\F$ of $\C\llp u \rrp$-vector spaces on $Y$,
a $u$-connection on $\F$ is a $\C$-linear map of sheaves
$\nabla_{\partial_u}: \F \to \F$ 
which satisfies the Leibniz rule with respect to
multiplication by $\C\llp u \rrp$ on sections.
More generally, given a strict complex of sheaves of $\C\llp u \rrp$-vector spaces
$(\F^{\bullet},d)$, a $u$-connection on $(\F^{\bullet},d)$ is defined similarly. 
\edefn

\brem
Given a sheaf $\F$ (or complex of sheaves $\F^{\bullet}$)
of $\C$-vector spaces on $Y$, we will often talk about 
``a $u$-connection on $\F$," by which we mean a $u$-connection
on the sheaf $\F\llp u \rrp$ (resp. complex of sheaves $\F\llp u \rrp$).
\erem

\brem
Given a bounded strict complex $(E^{\bullet},d)$
of $\C\llp u \rrp$-vector spaces with connection, we will often
identify $(E^{\bullet},d)$ with its image
under the localization map $\Ch^b(\calA_{\wE_{\C,0}})
\to \calD^b(\wE_{\C,0})$. Similarly, given a sheaf $\F$ of $\C$-vector
spaces, a $u$-connection on $\F\llp u \rrp$ determines an object
of $\calA(Y; \wE_{\C,0})$. We will
often identify a bounded complex $(\F^{\bullet}\llp u \rrp,d)$ of such sheaves 
with $u$-connection with its image under the 
localization map $\Ch^b(\calA(Y;\wE_{\C,0}))
\to \Shv^b(Y; \wE_{\C,0})$.
\erem

\brem
When it is clear from context that we are
working with a $u$-connection on an object 
with a $\C\llp u \rrp$- (or $\un{\C\llp u \rrp}$-)
action, we simply talk about ``the connection" on
that object.
\erem

Let $f \in \Gamma(X; \O_X)$ be a globally
defined function on $X$ (i.e.\ either an algebraic
function $X \to \Aone$ or holomorphic function
$X \to \C$. In this case, we denote 
by $\wih{\scrE}_X^{-f/u}$ the sheaf $\O_X\llp u \rrp$ 
equipped with the $\C\llp u \rrp$-linear connection given by
$d - df/u$. That is, given a section  $s \in \O_X(U)$,
$U \subset X$, and a formal Laurent series $\ell(u)
\in \C\llp u \rrp$,
\[\nabla(s \otimes \ell(u)) := ds \otimes \ell(u) - s \cdot df \otimes \frac{\ell(u)}{u}.\]

\subsubsection{The formal twisted de Rham complex}
Let $\scrM$ be a locally free $\O_X$-module of finite rank
equipped with a flat connection $\nabla$ having regular singularity
at infinity. Then $\scrM \otimes_{\O_X} \wih{\scrE}_{X}^{-f/u}$ is a locally
free $\O_X\llp u \rrp$-module with a canonical connection given by
$u \cdot \nabla -df \otimes {\id_{\scrM}}$. Sabbah's theorem
concerns the \emph{formal twisted de Rham complex} of $\scrM$.

\bdef
\label{def: formal twisted de Rham complex}
The formal twisted de Rham complex of $\scrM$ by $f$
is defined to be the complex of $\O_X\llp u \rrp$-modules
given by
\[\DR(\wih{\scrE}_X^{-f/u} \otimes_{\O_X} \scrM) := 
	(\Omega_X^{\bullet}\llp u \rrp \otimes_{\O_X} \scrM, u \cdot \nabla - df \otimes \id_{\scrM}).\]
This complex comes with a connection,
defined by $\nabla_{\partial_u} := \partial_u + f/u^2$, that commutes
with the differential. 
\edefn

\brem
The global sections,  $\Gamma(X; \DR(\wih{\scrE}_{X}^{-f/u} \otimes_{\O_X} \scrM))$,
give a complex of $\C\llp u \rrp$-vector spaces with connection $\nabla_{\partial_u}$. Its cohomology
groups are $\C\llp u \rrp$-vector spaces with connection.
\erem

\brem
We note that, because $\scrM$ is a coherent $\O_X$-module,
there is an isomorphism $\scrM\llp u \rrp = \scrM \otimes_{\O_X} \O_X\llp u \rrp$,
as well as quasi-isomorphisms, $\Omega_X^{\bullet}\llp u \rrp \otimes_{\O_X} \scrM
\simeq \Omega_X^{\bullet} \otimes_{\O_X} \scrM\llp u \rrp \simeq 
(\Omega_X^{\bullet} \otimes_{\O_X} \scrM)\llp u \rrp \simeq
\Omega_X^{\bullet} \otimes_{\O_X} \scrM \otimes_{\O_X} \O_X\llp u \rrp$,
of complexes in $\Ch^b(\QCoh^{\heart}(X))$.
\erem

\subsubsection{Vanishing cycles and microlocal formal Riemann--Hilbert}
If $X$ is a complex algebraic variety, it has an
associated complex analytic variety 
called the \emph{analytification} of $X$, which is a
complex manifold in the case that $X$ is smooth.
All of the algebraic objects that we have discussed above 
therefore have analytic counterparts. Denote the analytification of
$X$ by $X^{\an}$. Given an algebraic vector bundle
with connection, $(\scrM, \nabla)$, its analytification
is analytic vector bundle with flat connection
$(\scrM^{\an}, \nabla^{\an})$ on $X^{\an}$.
Consider the analytic function $f^{\an}: X^{\an} \to \C$
associated to the regular function $f$ above.

Let $\scrL = \ker \nabla^{\an}$ denote the local
system on $X^{\an}$ of flat sections of $\scrM^{\an}$.
This is a locally constant sheaf (in degree $0$) on
$X^{\an}$ with coefficients in $\C$.

\bdef[{e.g.\ \cite[Definition 4.2.4]{Dimca}}]
For each $c \in \C$, let $\varphi_{f-c} \scrL$ denote the
vanishing cycles of $\scrL$ with respect to
the function $f^{\an}-c$. It is an object of
$\Shv_c(X^{\an}_c;\C)$, where $X^{\an}_c = (f^{\an})^{-1}(c)$, 
and it is constructible with respect
to an \emph{algebraic} stratification of $X^{\an}$.
Moreover, we may assume this stratification is finite
(c.f.\ \cite[Corollary 4.1.8]{Dimca}).
\edefn

\brem
In this work, the shift in our definition for the vanishing cycles
functor is normalized so that the functor preserves perversity. 
\erem

The complex of vanishing cycles, $\varphi_{f-c}(\scrL)$,
has a canonical monodromy automorphism that we
denote by $T_c$. The hypercohomology groups,
\[H^k(f^{-1}(c); \varphi_{f-c}(\scrL))\]
are finite dimensional $\C$-vector spaces
equipped with an automorphism $T_{f-c}$.

\bdef
\label{def: inverse RH}
Let $E$ be a finite dimensional $\C$-vector
space equipped with an automorphism $T$. Writing
$T = \exp(-2\pi i M)$ for some $M: E \to E$, define
\[\invRH(E, T) := (E\llp u \rrp, d + Mdu/u),\]
where the right-hand side is a $\C\llp u \rrp$-vector
space with connection.\footnotemark
\edefn
\footnotetext{This functor is denoted $\wih{\RH}^{-1}$ in \cite{Sabbah}.}

\bdef
We define the following $\C\llp u \rrp$-vector space
with connection,
\[\wih{\scrE}^{-c/u} := (\C\llp u \rrp, d+ cdu/u^2).\]
\edefn

\begin{importantrem}
\label{important remark}
While $\invRH$ is a $u$-connection that determines
a regular holonomic $\wE_{\C,0}^{\poly}$-module,
hence a $\wE_{\C,0}$-module by \cref{lem: poly}, 
$\wih{\scrE}^{-c/u}$ is an \emph{irregular}
connection, so determines just a $\wE_{\C,0}$-module.
The Kontsevich--Sabbah--Saito
theorem therefore obtains an isomorphism
of $\wE_{\C,0}^{\poly}$-modules (or $\C\llp u \rrp$-modules
with connection) rather than an isomorphism
of $\wE_{\C,0}$-modules as one might hope.
\end{importantrem}

\subsubsection{The Kontsevich--Sabbah--Saito theorem}
The following theorem was conjectured by Maxim Kontsevich
as a possible answer to a question raised in \cite[\S7.4]{KoSo}
about how to define the vanishing cycles of a function on a
smooth formal scheme. It was proven by Claude Sabbah 
in the following form in \cite{Sabbah}.

\bthm[{\cite[Theorem 1.1]{Sabbah}}]
\label{thm: Sabbah's theorem}
Let $X$ be a smooth complex algebraic variety,
and $f: X \to \Aone$ a regular function on $X$.
Let $\scrM$ be a locally free $\O_X$-module of finite rank
equipped with a flat connection $\nabla$ having regular singularity
at infinity, and let $\scrL := \ker \nabla^{\an}$ denote its locally constant
sheaf of flat sections.
Then for each $k$, we have an isomorphism of $\C\llp u \rrp$-vector spaces
with connection,
\begin{multline}
\label{eqn: Sabbah's theorem}
\left(H^{k}(X; \DR(\wih{\scrE}_X^{-f/u} \otimes_{\O_X} \scrM)), \nabla_{\partial_u} \right)  
	\simeq \\ \bigoplus_{c \in \C} \wih{\scrE}^{-c/u} \otimes_{\C\llp u \rrp} 
		\invRH\left(H^k(f^{-1}(c); \varphi_{f-c}(\scrL)), T_c \right),
\end{multline}
where we have used the notation established above.
\ethm

\brem
Note that the direct sum on the right-hand side
of \labelcref{eqn: Sabbah's theorem} is finite
since $\varphi_{f-c}(\scrL) \simeq 0$ when
$c$ is a non-critical value of $f^{\an}$. Since $f^{\an}$ 
is holomorphic, it has only finitely many 
critical values.
\erem

\brem
Since $X^{\an}$ is a complex manifold, the underlying $\O_{X^{\an}}$-module
of $\scrM$ is given by $\scrL \otimes_{\C} \O_{X^{\an}}$, so we may also
write $\DR(\wE_X^{-f/u} \otimes_{\O_{X^{\an}}} \scrM^{\an}) \simeq 
\scrL \otimes_{\C} (\Omega_{X^{\an}}^{\bullet}\llp u \rrp, u \cdot \nabla - df \otimes \id_{\scrM^{\an}})$.
\erem

\subsection{Refinement of the Kontsevich--Sabbah--Saito theorem}

Let $X$ be a complex algebraic variety.
There is a continuous map of topological spaces
$\tau: X^{\an} \to X^{\Zar}$, where $X^{\Zar}$ denotes
the Zariski topology on $X$. Furthermore, 
$X^{\an}$ is of finite homotopy type,\footnotemark
\footnotetext{See \cite{CisinskiMO} for justification
of this claim.}  
and we can complete $\{X^{\an}\}$ 
to a basis $\calB(X)$ as in \cref{conv: finite homotopy type}
by virtue of the local topology of $X^{\an}$.

\bthm[Main Theorem]
\label{thm: main theorem}
Let $X$ be a smooth complex algebraic variety,
and $f: X \to \Aone$ a regular function on $X$.
Let $\scrM$ be a locally free $\O_X$-module of finite rank
equipped with a flat connection $\nabla$ having regular singularity
at infinity, and let $\scrL := \ker \nabla^{\an}$ denote its 
locally constant sheaf of flat sections.
Then there is an equivalence of objects in
$\Shv^b(X^{\Zar}; \wE_{\C,0}^{\poly})$,
\beqn
\label{eqn: main theorem}
\DR(\wE_X^{-f/u} \otimes_{\O_X} \scrM) 
	\simeq \bigoplus_{c \in \C} \un{\wih{\scrE}}^{-c/u} \otimes_{\C\llp u \rrp} 
		\tau_*\muRH_{\Shv}\left(\varphi_{f-c}(\scrL), T_{f-c} \right),
\eeqn
which recovers the isomorphisms 
\labelcref{eqn: Sabbah's theorem} 
upon taking hypercohomology.
\ethm

\bproof
By taking a covering of $X$ by quasi-projective
Zariski open sets, we reduce to the case where $X$
is quasi-projective.
By \cref{cor: major corollary} below, there is
an equivalence, 
\[{\DR(\wE_X^{-f/u} \otimes_{\O_{X^{\an}}} 
\scrM^{\an})}|_{X^{\an}_c} \simeq \un{\wE}^{-c/u} \otimes_{\C\llp u \rrp} \muRH_{\Shv}(\varphi_{f-c}(\scrL), T_{f-c}),\]
of objects in $\Shv_c(X^{\an}_c; \wE_{\C,0})$.
Since, $\varphi_{f-c}(\scrL) \simeq 0$ for all but finitely many
$c \in \C$, we find that
\[\DR(\wE_X^{-f/u} \otimes_{\O_{X^{\an}}} \scrM^{\an}) \simeq 
\bigoplus_{c \in \C} \un{\wih{\scrE}}^{-c/u} \otimes_{\C\llp u \rrp} 
		\muRH_{\Shv}\left(\varphi_{f-c}(\scrL), T_{f-c} \right).\]
Apply $\tau_*$ to the above
equivalence. Observing that $\tau_*$ commutes
with finite direct sums and using the projection formula,
we obtain
\begin{align*}
\tau_*\left(\bigoplus_{c \in \C} \un{\wih{\scrE}}^{-c/u} \otimes_{\C\llp u \rrp} \muRH_{\Shv}\left(\varphi_{f-c}(\scrL), T_{f-c} \right)\right) 	&\simeq \bigoplus_{c \in \C} \tau_*\left(\un{\wih{\scrE}}^{-c/u} \otimes_{\C\llp u \rrp} \muRH_{\Shv}\left(\varphi_{f-c}(\scrL), T_{f-c} \right)\right) \\
																									&\simeq \bigoplus_{c \in \C} \tau_*\left(\tau^*\un{\wE}_{X^{\Zar}_c}^{-c/u} \otimes_{\C\llp u \rrp} \muRH_{\Shv}\left(\varphi_{f-c}(\scrL), T_{f-c} \right)\right) \\
																									&\simeq \bigoplus_{c \in \C} \un{\wih{\scrE}}^{-c/u} \otimes_{\C\llp u \rrp} \tau_*\muRH_{\Shv}\left(\varphi_{f-c}(\scrL), T_{f-c} \right).
\end{align*}
On the other hand, since we have reduced to the
case where $X$ is quasi-projective, \cref{prop: major prop 1} below
establishes a natural algebraic-analytic comparison equivalence,
$\DR(\wE_X^{-f/u} \otimes_{\O_X} \scrM) \xrightarrow{\simeq} 
{\tau_X}_*\DR(\wE_X^{-f/u} \otimes_{\O_{X^{\an}}} \scrM^{\an})$,
from which the result follows.
\eproof

\section{Proof of the Main Theorem}
\label{sec: Proof of the Main Theorem}

The proof of \cref{thm: main theorem} we have given above 
closely imitates Sabbah's proof of \cite[Theorem 1.1]{Sabbah}. 
Namely, it proceeds in two broad parts:

\begin{enumerate}[(i)]
\item 
\label{Sabbah item 1}
A comparison of the 
\textit{algebraic} formal twisted
de Rham complex $\DR(\wE^{-f/u} \otimes_{\O_X} \scrM)$
to the \textit{analytic} formal twisted
de Rham complex, $\DR(\wE^{-f/u} \otimes_{\O_{X^{\an}}} \scrM^{\an})$
(\cref{prop: major prop 1}).

\item 
\label{Sabbah item 2}
A comparison of the
analytic formal twisted de Rham
complex to the right-hand side
of \labelcref{eqn: main theorem}
(\cref{prop: major prop 2} and \cref{cor:
major corollary}). 
\end{enumerate}

In the remaining sections of this paper, before the
appendices, we establish
each of these component parts.

\subsection{Algebraic-to-analytic comparison}
Let $F: Y \to \bbP^1$ be regular function
from a smooth projective variety $Y$ to $\bbP^1$
extending $f$, meaning that there exists a commutative
diagram,
\[\begin{tikzcd}
	X & Y \\
	\Aone & {\bbP^1}
	\arrow["f"', from=1-1, to=2-1]
	\arrow["F", from=1-2, to=2-2]
	\arrow["j", hook, from=1-1, to=1-2]
	\arrow[hook, from=2-1, to=2-2],
\end{tikzcd}\]
and such that $D := Y \setminus X$ is a (not necessarily
normal crossings) divisor
in $Y$. Then $(\scrM, \nabla)$ extends as a
coherent $\O_Y(\ast D)$-module\footnotemark 
with connection having regular singularity along $D$, that we
continue to denote by $(\scrM, \nabla)$.
In particular, $\scrM$ is $\scrD_Y$-holonomic,
hence $\scrD_Y$-coherent.
There are several sheaves we define on $Y$:
\begin{itemize}
\item $\wih{\O}_{(Y,D)} := \O_Y(\ast D)\llp u \rrp$.
\item $\wih{\scrD}_{(Y,D)} := \wih{\O}_{(Y,D)} \otimes_{\O_Y} \scrD_Y$.
\item $\wih{\scrE}^{-F/u}_{(Y,D)}$ denotes the 
$\wih{\scrD}_{(Y,D)}$-module, $\O_Y(\ast D)\llp u \rrp$,
equipped with the connection $(-dF \wedge -) + ud$.
\end{itemize}

\footnotetext{$\O_Y(\ast D)$ denotes the sheaf of rational
functions on $Y$ that are regular on $Y \setminus D$.}

With the above set-up and notation, we recall
the \textit{algebraic formal twisted de Rham complex}
associated to $(\scrM, D)$:
\beqn
\label{eqn: algebraic dR complex}
	\left(\Omega_Y^{\bullet}(\ast D)\llp u \rrp 
			\otimes_{\O_Y} \scrM, u \cdot \nabla + 
				(-dF \wedge -) \otimes \id_{\scrM}\right)
\eeqn
Similarly, recall the
\textit{analytic formal twisted de Rham complex}
associated to $(\scrM, D)$:
\beqn
\label{eqn: analytic dR complex}
	\left(\Omega_{Y^{\an}}^{\bullet}(\ast D)\llp u \rrp 
		\otimes_{\O_{Y^{\an}}} \scrM^{\an}, u \cdot \nabla + 
			(-dF \wedge -) \otimes \id_{\scrM^{\an}}\right).
\eeqn

We denote \labelcref{eqn: algebraic dR complex} and
\labelcref{eqn: analytic dR complex} by 
$\DR(\wE_{Y,D}^{-F/u} \otimes_{\O_Y} \scrM)$ and 
$\DR(\wE_{Y,D}^{-F/u} \otimes_{\O_{Y^{\an}}} \scrM^{\an})$,
respectively, and regard them as 
objects of $\Shv(Y;\wE_{\C,0})$
and $\Shv(Y^{\an}; \wE_{\C,0})$, respectively.

\subsubsection{Reminder on GAGA}
We have already introduced the associated analytic
space $X^{\an}$ to a complex algebraic variety $X$.
The definition of this associated analytic space may be
found in \S1 of \cite{GAGA}. As alluded to above, there
is a map on underlying topological spaces,
\[\tau: X^{\an} \to X^{\Zar},\]
which is continuous. The pullback of the structure
sheaf $\O_X$ on $X$ is a subsheaf of the sheaf of
holomorphic functions,
\[\tau^*\O_X \subset \O_{X^{\an}}.\]

\bdef[{\cite[\S3.9 Definition 2]{GAGA}}]
The \emph{analytification} of an $\O_X$-module
$\F$ is defined to be
\[\F^{\an} := \tau^*\F \otimes_{\tau^*\O_X} \O_{X^{\an}}.\]
\edefn

There is a canonical morphism 
\[\tau^*\F := \tau^*\F 
	\otimes_{\tau^*\O_X} \tau^*\O_X \to \tau^*\F 
		\otimes_{\tau^*\O_X} \O_{X^{\an}} =: \F^{\an}\]
induced by the inclusion $\tau^*\O_X \subset \O_{X^{\an}}$.
Composition of this morphism with the unit
of the $(\tau^*, \tau_*)$ adjunction yields
a natural map,
\beqn
\label{eqn: alg-an map}
\F \to \tau_*\F^{\an}.
\eeqn


\bprop[{\cite[Proposition 1]{SabbahSaito}}]
Let $\scrM$ be a coherent $\scrD_Y$-module.
Then the natural morphism \labelcref{eqn: alg-an map} induces
an equivalence on global sections
\beqn
\label{eqn: alg-an global sections}
\Gamma(Y; \DR(\wE_{Y,D}^{-F/u} \otimes_{\O_Y} \scrM)) \xrightarrow{\simeq} 
		\Gamma(Y^{\an}; \DR(\wE_{Y,D}^{-F/u} \otimes_{\O_{Y^{\an}}} \scrM^{\an})),
\eeqn
as objects in $\calD(\wE_{\C,0})$.
\eprop

\bproof
This follows immediately from the proof of
Proposition 1 in \cite{SabbahSaito}.
\eproof

We record the following proposition of Sabbah, which
relates the twisted de Rham complex of $X^{\an}$
and the twisted de Rham complex of its compactification
$Y^{\an}$.

\bprop[{\cite[Proposition 4.1]{Sabbah}}]
\label{prop: Proposition 4.1}
Let $X$ be a smooth complex variety, and $f: X \to \Aone$
a regular function. Let $F: Y \to \bbP^1$ be a map from
a smooth projective variety to $\bbP^1$, extending $f$.
We denote by $j^{\an}: X^{\an} \hook Y^{\an}$ the open inclusion, and by
$D = Y^{\an} \setminus X^{\an}$ the divisor given by the complement of $X^{\an}$.
Let $\scrM$ be a coherent $\O_{X^{\an}}(\ast D)$-module equipped with a
flat connection $\nabla$ with regular singularities along $D$, where $\O_{X^{\an}}(\ast D)$
denotes the sheaf of meromorphic functions with poles at most along $D$. Then
the natural morphism of sheaves,
\[\DR(\wE_{Y,D}^{-F/u} \otimes_{\O_X(\ast D)} \scrM) \xrightarrow{\simeq} Rj_*\DR(\wE_X^{-f/u} \otimes_{\O_X} \scrM|_X)\]
is an equivalence.
\eprop

Although our standing assumption is that all functors
are implicitly derived, we have written ``$Rj_*$" in the
proposition above to emphasize that Sabbah's result is
an equivalence of objects in the derived category.

\bprop
\label{prop: major prop 1}
Let $X$ be a smooth complex algebraic variety,
and let $f: X \to \Aone$ be a regular function. 
Let $\scrM$ be a coherent $\scrD_X$-module,
and let $\tau_X: X^{\an} \to X$ denote the
canonical map on topological spaces. Then the
natural morphism \labelcref{eqn: alg-an map} induces
an equivalence
\beqn
\label{eqn: alg-an comparison map}
\DR(\wE_X^{-f/u} \otimes_{\O_X} \scrM) \xrightarrow{\simeq} 
		{\tau_X}_*\DR(\wE_X^{-f/u} \otimes_{\O_{X^{\an}}} \scrM^{\an}),
\eeqn
between objects in $\Shv(X;\wE_{\C,0})$.
\eprop

\bproof
It suffices to show that the natural morphism
\labelcref{eqn: alg-an comparison map} induces
an equivalence on sections over each $U \subset X$,
\beqn
\label{eqn: alg-an sections map}
\Gamma(U;\DR(\wE_X^{-f/u} \otimes_{\O_X} \scrM)) \xrightarrow{\simeq} 
		\Gamma(U;\tau_*\DR(\wE_X^{-f/u} \otimes_{\O_{X^{\an}}} \scrM^{\an})).
\eeqn
The left-hand side of \labelcref{eqn: alg-an sections map}
is clearly equivalent to global sections of the sheaf,
$\DR(\wE_U^{-f|_U/u} \otimes_{\O_U} \scrM|_U)$. On the
other hand, by smooth base change\footnotemark 
\footnotetext{See \cite[Lemma 3.3]{Volpe} for 
a reference for smooth base change
in the setting of sheaves valued in 
presentable $\infty$-categories.}
using the following
pullback square of topological spaces,
\[\begin{tikzcd}
	{U^{\an}} & {X^{\an}} \\
	U & X
	\arrow["{\tau_X}", from=1-2, to=2-2]
	\arrow["{\tau_U}"', from=1-1, to=2-1]
	\arrow["{\jmath^{\an}}", hook, from=1-1, to=1-2]
	\arrow["\jmath", hook, from=2-1, to=2-2]
	\arrow["\lrcorner"{anchor=center, pos=0.125}, draw=none, from=1-1, to=2-2],
\end{tikzcd}\]
the right-hand side of \labelcref{eqn: alg-an sections map}
is equivalent to global sections of the sheaf 
$\DR(\wE_U^{-f|_U/u} \otimes_{\O_{U^{\an}}} (\scrM|_U)^{\an})$.
Observe that, as an open subset of $X$,
$U$ is also a smooth complex variety. Thus,
without loss of generality, it suffices to prove that
\labelcref{eqn: alg-an sections map} is an equivalence
for $U = X$.

In order to show the equivalence for global sections
on $X$, we invoke \cite[Proposition 4.1]{Sabbah}, which we have reproduced
above as \cref{prop: Proposition 4.1}. Consider the
following pullback square
\[\begin{tikzcd}
	{X^{\an}} & {Y^{\an}} \\
	X & Y
	\arrow["{j^{\an}}", hook, from=1-1, to=1-2]
	\arrow["j", from=2-1, to=2-2]
	\arrow["{\tau_X}"', from=1-1, to=2-1]
	\arrow["{\tau_Y}", from=1-2, to=2-2]
	\arrow["\lrcorner"{anchor=center, pos=0.125}, draw=none, from=1-1, to=2-2]
\end{tikzcd}\]
of topological spaces. It follows that
\begin{align}
\Gamma(Y^{\an}; \DR(\wE_{Y,D}^{-F/u} \otimes_{\O_{Y^{\an}}(\ast D)} j^{\an}_+\scrM^{\an}))&\simeq \pt_*{\tau_Y}_*\DR(\wE_{Y,D}^{-F/u} \otimes_{\O_{Y^{\an}}(\ast D)} j^{\an}_+\scrM^{\an}) \label{alg-an 1}\\
																		&\simeq \pt_*{\tau_Y}_*j^{\an}_*\DR(\wE_X^{-f/u} \otimes_{\O_{X^{\an}}} \scrM^{\an}) \label{alg-an 2}\\
																		&\simeq \pt_*j_*{\tau_X}_*\DR(\wE_X^{-f/u} \otimes_{\O_{X^{\an}}} \scrM^{\an}) \label{alg-an 3}\\
																		&\simeq \Gamma(X^{\an}; {\tau_X}_*\DR(\wE_X^{-f/u} \otimes_{\O_{X^{\an}}} \scrM^{\an})) \label{alg-an 6},
\end{align}
where \labelcref{alg-an 2} follows from \cref{prop: Proposition 4.1}.
On the other hand, by \cite[Remark 2.8]{Sabbah},
we have the equivalence,
\[\DR(\wE_{Y,D}^{-F/u} \otimes_{\O_Y} j_+\scrM) \xrightarrow{\simeq} 
	\left(\Omega_Y^{\bullet} \otimes_{\O_Y} j_*\scrM\llp u \rrp, 
		u \cdot \nabla + (-dF \wedge -)\right).\]
Because we are in the algebraic setting, $j_*\scrM\llp u \rrp \simeq (j_*\scrM)\llp u \rrp$.
It follows that 
\begin{align*}
\Omega_Y^{\bullet} \otimes_{\O_Y} (j_*\scrM)\llp u \rrp &\simeq (\Omega_Y^{\bullet} \otimes_{\O_Y} j_*\scrM)\llp u \rrp \\
										&\simeq (j_*(j^*\Omega_Y^{\bullet} \otimes_{\O_X} \scrM)\llp u \rrp \\
										&\simeq j_*(\Omega_X^{\bullet} \otimes_{\O_X} \scrM)\llp u \rrp,
\end{align*}
where we have used the projection
formula to obtain the penultimate equivalence.
Thus, we obtain
\begin{align*}
\Gamma(Y;\DR(\wE_{Y,D}^{-F/u} \otimes_{\O_Y(\ast D)} j_+\scrM)) &\simeq \Gamma\left(Y; \left(j_*(\Omega_X^{\bullet} \otimes_{\O_X} \scrM)\llp u \rrp, u \cdot \nabla + (-dF \wedge -) \right)\right) \\
													&\simeq \Gamma\left(X; \left(\Omega_X^{\bullet} \otimes_{\O_X} \scrM)\llp u \rrp, u \cdot \nabla + (-dF \wedge -) \right)\right) \\
													&= \Gamma(X; \DR(\wE_X^{-f/u} \otimes_{\O_X} \scrM)).
\end{align*}
The result now follows from \cite[Proposition 2.6]{Sabbah},
as illustrated in the following commutative diagram,
\[\begin{tikzcd}
	{\Gamma(Y;\DR(\wE_{Y,D}^{-F/u} \otimes_{\O_Y} j_+\scrM))} & {} & {\Gamma(Y^{\an}; \DR(\wE_{Y,D}^{-F/u} \otimes_{\O_{Y^{\an}}} j^{\an}_+\scrM^{\an}))} \\
	{\Gamma(X; \DR(\wE_X^{-f/u} \otimes_{\O_X} \scrM))} & {} & {\Gamma(X^{\an}; {\tau_X}_*\DR(\wE_X^{-f/u} \otimes_{\O_{X^{\an}}} \scrM^{\an}))}
	\arrow["\text{\labelcref{eqn: alg-an sections map}}", from=2-1, to=2-3]
	\arrow["\text{\cite[Proposition 2.6]{Sabbah}}", "\simeq"', curve={height=-18pt}, from=1-1, to=1-3]
	\arrow["\simeq"', from=1-3, to=2-3]
	\arrow["{\simeq }"', from=1-1, to=2-1].
\end{tikzcd}\]
\eproof

\subsection{Analytic-to-vanishing cycles comparison}
Now let $X$ be a complex analytic manifold
of dimension $n$, and $f: X \to \C$ a holomorphic
function on $X$. Let $\scrM$ be a locally free $\O_X$-module
of finite rank $d$ equipped with a flat holomorphic
connection $\nabla$, whose associated local system
of flat sections, $\ker \nabla$ is denoted by $\scrL$.

Let $\pr_2: X \times \C \to \C$ denote the canonical projection.
Using $\scrM$ and its connection $\nabla$, Sabbah defines the following
object of $\calD^b(\pr_2^*\scrD_{\C})$.

\bdef
Denote by $i_f: X \to X \times \C$
the graph of $f$, and let $\DR_{X \times \C/\C}$
denote the relative de Rham functor. Then $\scrK_f$
is defined to be $\DR_{X \times \C/\C}(\scrM_f)$,
where $\scrM_f := {i_f}_+\scrM$.\footnotemark
\edefn
\footnotetext{Here, ${i_f}_+$ denotes the pushforward
of $\scrD$-modules. Because $i_f$ is a closed immersion,
${i_f}_+$ is exact on ordinary $\scrD_X$-modules, and 
therefore $\scrM_f \in \RegHol_{\scrD_{X \times \C}}$.}

Let $X_0 := f^{-1}(0)$. We will be interested
in the restriction of $\scrK_f$ to the analytic subspace
$X_0 \subset X$, and denote this restriction by $\scrK_{f,0}
:= {\scrK_f}|_{X_0}$. It follows from the definition
that $\scrK_{f,0}$ is an object of $\Shv^b(X_0; \scrD_{\C,0})$.

Using the locally free resolution
of the relative de Rham complex 
recalled earlier in \labelcref{eqn: resolution of relative de Rham complex},
Sabbah works with an explicit model for $\scrK_{f,0}$
in $\Ch^b(X_0; \scrD_{\C,0})$, the category
of bounded complexes of $\un{\scrD_{\C,0}}$-modules
on $X_0$, in order to prove the results of \cite{Sabbah}.
We denote the complex in $\Ch^b(X_0;\scrD_{\C,0})$ representing
$\scrK_{f,0}$ by $\scrK_{f,0}^{\bullet}$. 
The complex $\scrK_{f,0}^{\bullet}$ admits a filtration
which is compatible with the filtration, $V^{\bullet}\scrD_{\C,0}$,
making it into a ``$V$-filtered complex" in the words of \cite{Sabbah}.

\bdef 
We set $U^k \scrK_{f,0}^{\bullet} 
:= {\DR_{X \times \C/\C}(V^k\scrM_f)}|_{X_0}$,
where $V^k\scrM_f$ denotes the Kashiwara-Malgrange
filtration along the hypersurface $X \times 0 \subset X \times \C$.
\edefn

\brem
As noted in \cite{Sabbah}, it makes sense to apply
$\DR_{X \times \C/\C}$ to $V^k\scrM_f$ because 
it is a coherent $V^0\scrD_{X \times \C}$-module,
and hence a $\scrD_{X \times \C/\C}$-module.
\erem

The main properties of the filtration
$U^{\bullet}\scrK_{f,0}^{\bullet}$ are established 
in \cite[Theorem 3.9]{Sabbah}:
\begin{itemize}
\item[$-$]
The complex, $\gr^{-1}_U\scrK_{f,0}^{\bullet}$ equipped with
the operator $\exp(-2\pi it\partial_t)$ is quasi-isomorphic to
the complex $\varphi_f(\scrL)$ equipped with its
monodromy operator, $T_f$.

\item[$-$] 
$\scrK_{f,0}^{\bullet}$ satisfies the assumptions of
\cite[Corollary 3.4]{Sabbah}. Namely,
\vspace{3pt}
\begin{enumerate}[(i)]
\item the image of $\scrK_{f,0}^{\bullet}$ under the canonical
localization functor, $\Ch^b(X_0; \scrD_{\C,0}) 
\to \Shv^b(X; \scrD_{\C,0})$, where $K^b(X_0; \scrD_{\C,0})$
denotes the ordinary category of bounded complexes 
of $\un{\scrD_{\C,0}}_{X_0}$-modules,
lies in $\Shv_c(X_0; \scrD_{\C,0})$,
\item $t: U^k(\scrK_{f,0}^i) \to U^{k+1}(\scrK_{f,0}^i)$ is an isomorphism of $\C\{t\}$-modules for each $k > 0$ and each $i$,
\item $U^{-k-1}(\scrK_{f,0}^i) = \sum_{\ell=0}^k \partial_t^{\ell} U^{-1}(\scrK_{f,0}^i)$ for each $k > 0$ and each $i$,
\item there exists a nonzero polynomial $b(s) \in \C[s]$ with roots having their real part
in $[0, 1)$ such that $b(t\partial_t - k)$ vanishes on each $\gr_U^k(\scrK_{f,0}^i)$ for each $i$, and
\item for each $i$, each germ $H^i(U^{-1}\scrK_{f,0}^{\bullet})_y$ is finite type over $\C\{t\}$.
\end{enumerate}
\end{itemize}

It is clear from the properties listed above that
the $V$-filtered complex $U^{\bullet}\scrK_{f,0}^{\bullet}$ endows $\scrK_{f,0}$
with the structure of a module in $\Fil(\Shv(X_0; \Sp)$ over
$V^{\bullet}\un{\scrD_{\C,0}}$ satisfying the properties enumerated
in \cref{prop: uniqueness of sfV-filtration}.
Moreover, it is clear from the definition of $\scrK_{f,0}^{\bullet}$
that it is constructible with respect to a finite stratification of $X_0$.
Altogether, we have the following corollary.

\bcor
\label{cor: scrK is filtered module}
The object $\scrK_{f,0}$ is an object
of $\Shv_c^{f.s.}(X_0; \scrD_{\C,0})$ and
has the canonical enhancement to an
object of $\Mod_{V^{\bullet}\un{\scrD_{\C,0}}}
(\Fil(\Shv_c^{f.s.}(X_0;\Sp))$
satisfying the properties enumerated 
in \cref{prop: uniqueness of sfV-filtration},
coming from the $V$-filtered complex, 
$U^{\bullet}\scrK_{f,0}^{\bullet}$.
\ecor

\bprop
\label{prop: two kinds of vanishing cycles}
There is a natural equivalence,
\[\varphi_{\Shv}(\scrK_{f,0}) \simeq (\varphi_f(\scrL), T_f)\]
of objects in $\Fun(S^1, \Shv_c^{f.s}(X_0; \C))$.
\eprop

\bproof
Recall that $\varphi_f(\calL) \simeq \DR_X(\gr^{-1}_V \scrM_f)$
by a theorem independently obtained 
by Kashiwara and Malgrange. On the other hand,
$U^i\scrK_{f,0}^{\bullet}$ is by definition the 
image of $V^i({\scrM_f}|_{X_0 \times 0})$
under the relative de Rham functor $\DR_{X \times \C/\C}(-)$.
Since $\DR_{X \times \C/\C}(-)$ commutes with colimits
(given by tensor with relative differentials),
we obtain,
\begin{align*}
\gr^{-1}_U(\scrK_{f,0}^{\bullet}) 	&\simeq \DR_{X \times \C/\C}(V^{-1}\scrM_f)/\DR_{X \times \C/\C}(V^0\scrM_f) \\
							&\simeq \DR_{X \times \C/\C}(\gr^{-1}_V \scrM_f) \\
							&\simeq \DR_X(\gr^{-1}_V \scrM_f) \\
							&\simeq \varphi_f(\scrL),
\end{align*}
where we have used the fact that $\gr^{-1}_V \scrM_f$ is
supported on $X \times 0$, and $\DR_X \simeq \DR_{X \times \C/\C}$
on $\scrD_{X \times \C}$-modules supported on $X \times 0$.

Let $U^{\bullet}\scrK_{f,0}$ denote the filtered
$V^{\bullet}\un{\scrD_{\C,0}}$-module of 
\cref{cor: scrK is filtered module}. 
Since $\gr^{-1}_{\sfV}(-) \simeq \varphi_{\Shv}(-)$, it suffices to 
show that $U^{\bullet}\scrK_{f,0}$
and $\sfV^{\bullet}\scrK_{f,0}$ are equivalent
objects of $\Mod_{V^{\bullet}\un{\scrD_{\C,0}}}\Fil(\Shv(X_0; \Sp))$.
By \cref{prop: uniqueness of sfV-filtration}, it suffices to
show that $U^{\bullet}\scrK_{f,0}$ and 
$\sfV^{\bullet}\scrK_{f,0}$ each satisfy the
properties enumerated therein. For $\sfV^{\bullet}\scrK_{f,0}$
this is the content of \cref{prop: derived KM filtration satisfies properties}.
For $U^{\bullet}\scrK_{f,0}$, this is the content of
\cref{cor: scrK is filtered module}.
\eproof

\bprop
\label{prop: major prop 2}
There is a natural equivalence,
\[{\DR(\wE_X^{-f/u} \otimes_{\O_X} \scrM)}|_{X_0} \simeq \muRH_{\Shv}(\varphi_f(\scrL), T_f)\]
as objects of $\Shv_c(X_0; \wE_{\C,0})$.
\eprop

\bproof
There is a natural morphism 
\beqn
\label{eqn: Sabbah 2}
\wih{\scrK}_f^{\bullet} := \wih{\scrE}_{\C,0} \otimes_{\scrD_{\C,0}} 
	\scrK_f^{\bullet} \to {\DR(\wE_X^{-f/u} \otimes_{\O_X} \scrM)} \simeq 
		\scrL \otimes_{\C} (\Omega_X^{\bullet}\llp u \rrp, ud + (-df \wedge -))),
\eeqn
obtained by identifying $u$ with
$u$, which is in fact a quasi-isomorphism (after
restriction to $X_0$) by
\cite[Lemma 3.15]{Sabbah}. As such,
we obtain the equivalence,
\begin{align}
\muRH(\varphi_f(\scrL), T_f) 				&\simeq \muRH(\gr^{-1}_U \scrK_{f,0}, \exp(-2\pi it\partial_t)) \label{monke1} \tag{a}\\
									&\simeq \muRH(\varphi_{\Shv} \scrK_{f,0}, \exp(-2\pi it\partial_t)) \label{monke2} \tag{b}\\
									&\simeq \wih{\scrK}_{f,0} \label{monke3} \tag{c}\\
									&\xrightarrow{\simeq} {\DR(\wE_X^{-f/u} \otimes_{\O_X} \scrM)}|_{X_0} \label{monke4} \tag{d}
\end{align}
where the equivalence \labelcref{monke2} is
the content of \cref{prop: two kinds of vanishing cycles}
above, and the equivalence \labelcref{monke3}
is shown in \cref{exmp: mu_shv = wih_shv}.
\eproof

\bcor
\label{cor: major corollary}
More generally, for any $c \in \C$, there is a natural equivalence,
\[{\DR(\wE_X^{-f/u} \otimes_{\O_X} \scrM)}|_{X_c} \simeq 
	\un{\wE}^{-c/u}_{X_c} \otimes_{\C\llp u \rrp} \muRH_{\Shv}(\varphi_{f-c}(\scrL), T_{f-c})\]
as objects of $\Shv_c(X_c; \wE_{\C,0}^{\poly})$.
\ecor

\bproof
Let $g := f-c$. Then by \cref{prop: major prop 2}, 
${\DR(\wE_X^{-g/u} \otimes_{\O_X} \scrM)}|_{g^{-1}(0)} 
\simeq \muRH_{\Shv}(\varphi_g(\scrL), T_g)$. Now observe
that the connection on $\DR(\wE_X^{-g/u} \otimes_{\O_X} \scrM)$
is given by $\nabla_{\partial_u} = \partial_u + g/u^2 =
\partial_u + f/u^2 - c/u^2$. Thus, we find that
$\DR(\wE_X^{-g/u} \otimes_{\O_X} \scrM) \simeq
\un{\wE}^{-c/u} \otimes_{\C\llp u \rrp} \DR(\wE_X^{-g/u} \otimes_{\O_X}
\scrM)$ as objects of $\Shv_c(X; \wE_{\C,0})$.
The constant sheaf $\un{\wE}^{\,c/u}$ has a multiplicative
inverse given by $\un{\wE}^{-c/u}$ under the tensor
product of sheaves valued in $\calD^b(\wE_{\C,0})$, since 
$\wE^{-c/u} \otimes_{\C\llp u \rrp} \wE^{\,c/u}
\simeq \C\llp u \rrp$, where the latter denotes the
$\C\llp u \rrp$-module with trivial $\partial_u$ action.
The corollary now follows.
\eproof

\appendix

\section{Constructible sheaves}
\label{sec: Appendix A}
In this section, we recall the basics of constructible sheaves in the
higher categorical setting. The main reference for this section
is \cite[Appendix A]{HA}.
Whenever possible, however, we
have used the notation found there.

\subsection{Locally constant sheaves}
Constructible sheaves are defined in terms
of locally constant sheaves, whose definition
we recall below.

\bdef[{\cite[Definition A.1.12]{HA}}]
\label{def: locally constant sheaf}
Let $\calX$ be an $\infty$-topos, and let $\F$ be an object of $\calX$. We will say that $\F$
is constant if it lies in the essential image of the geometric morphism $\pi^*: \calS \to \calX$. 
We will say that $\F$ is locally constant if there exists a 
small collection of objects $\{U_{\alpha} \in \calX\}_{\alpha \in S}$ 
such that the following conditions are satisfied:
\begin{enumerate}[(i)]
\item The objects $U_{\alpha}$ cover $\calX$: that is, there is an effective epimorphism 
$\coprod U_{\alpha} \to \bm{1}$, where $\bm{1}$ denotes the final object of $\calX$.
\item For each $\alpha \in S$, the product $\F \times U_{\alpha}$ 
is a constant object of the $\infty$-topos $\calX_{/U_{\alpha}}$.
\end{enumerate}
\edefn

Under certain conditions on $\calX$, the full subcategory
of locally constant objects of $\calX$ can be described as
functors out of a particular $\infty$-groupoid, called the
\emph{shape} of $\calX$.

\bdef[{\cite[Definition A.1.1]{HA}}]
Let $\calX$ be an $\infty$-topos, let $\pi_*: \calX \to \calS$ be the functor corepresented by the
final object of $\calX$, and let $\pi^*$ be the left adjoint to $\pi_*$. We will say that $\calX$ 
has constant shape if the composition $\pi_*\pi^*: \calS \to \calS$ is corepresentable.
\edefn

\bdef[{\cite[Definition A.1.5]{HA}}]
Let $\calX$ be an $\infty$-topos. We will say that an object 
$U \in \calX$ has constant shape if the $\infty$-topos $\calX_{/U}$ 
has constant shape. We will say that $\calX$ is locally of constant shape if every object
$U \in \calX$ has constant shape.
\edefn

If $\calX$ is locally of constant shape,
the shape of $\calX$ may
be identified with the object $\pi_!\bm{1} \in \calS$, where
$\pi_!$ is the left adjoint to $\pi^*$ (see \cite[Proposition A.1.8]{HA}).

\bthm[{\cite[Theorem A.1.15]{HA}}]
\label{thm: Lurie thm 1}
Let $\calX$ be an $\infty$-topos which is locally of constant shape, 
and let $\psi^*: \calS_{\pi_!\bm{1}} \to \calX$ be the functor of \cite[Proposition A.1.11]{HA}. 
Then $\psi^*$ is a fully faithful embedding, whose essential image
is the full subcategory of $\calX$ spanned by the locally constant objects of $\calX$.
\ethm

\subsubsection{Monodromy equivalence}
The classical monodromy equivalence states
that, for a topological space $X$, there is
an equivalence,
\[\Fun(\pi_{\leq 1}X, \Set) \simeq \Loc(X)\]
between representations of the fundamental
groupoid of $X$ and locally constant sheaves
of sets on $X$. This classical equivalence generalizes
to an equivalence between the \emph{full} $\infty$-groupoid
of $X$ and locally constant sheaves of \emph{spaces}
on $X$, as we now recall.

To every topological space $X$, the associated $\infty$-groupoid
is denoted by the $\Sing(X)$. The functor $\Sing: \Top \to \calS$
has a left adjoint called \emph{geometric realization}, denoted by
$|-|: \calS \to \Top$.

\bdef[{\cite[Definition A.4.9]{HA}}]
Let $X$ be a topological space. We will say that $X$ has singular shape if the
counit map $|\Sing(X)| \to X$ is a shape equivalence.
\edefn

\bdef[{\cite[Definition A.4.15]{HA}}]
We will say that topological space $X$ is locally 
of singular shape if every open set $U \subseteq X$ has singular shape.
\edefn

If $X$ is locally of singular shape, the $\infty$-topos
$\Shv(X;\calS)$ is locally of constant shape. It follows
from \cref{thm: Lurie thm 1} that the locally constant 
sheaves on $X$ are identified with representations of
the shape of $\calX$. Moreover the the shape of 
$\calX$ is identified with the homotopy type $\Sing(X)$ of $X$
itself (see \cite[\S A.4]{HA}).

\subsubsection{Presentable stable coefficients}
By imitating \cref{def: locally constant sheaf},
we obtain a notion of locally constant sheaves for
more general coefficients in presentable stable $\infty$-categories.

Indeed, by \cite[Remark I.1.3.2.8]{SAG}, any geometric
morphism of $\infty$-topoi $g^*: \calX \to \calY$ induces
a functor of spectral sheaves, $\Shv(\calX; \Sp) \to \Shv(\calY; \Sp)$,
which by abuse of notation we also denote by $g^*$. As also
noted in the remark, $g^*$ is left adjoint to the pushforward
functor $g_*: \Shv(\calY; \Sp) \to \Shv(\calX; \Sp)$, 
given by pointwise composition with $g^*: \calX \to \calY$.

More generally, given any presentable stable $\infty$-category
$\scrC$, $\Shv(\calX; \scrC) \simeq \Shv(\calX; \Sp) \otimes \scrC$.
By functoriality, we obtain morphisms $g^*: \Shv(\calX; \scrC) 
\to \Shv(\calY; \scrC)$ and $g_*: \Shv(\calY; \scrC) \to \Shv(calX; \scrC)$
associated to any geometric morphism $g^*$. In particular, given
$U \in \calX$, we may talk about the restriction of $\F \in \Shv(\calX; \scrC)$
to $\calX_{/U}$, denoted $\F|_U$.

\bdef
Let $\calX$ be an $\infty$-topos, and $\scrC$ a presentable
stable $\infty$-category. Let $\F$ be an object of $\Shv(\calX; \scrC)$. 
We will say that $\F$ is constant if it lies in the essential image 
of a functor $\pi^*: \scrC \to \Shv(\calX; \scrC)$ induced
by the geometric morphism $\pi^*: \calS \to \calX$. We will say that
$\F$ is locally constant if there exists a small collection of objects $\{U_{\alpha} \in \calX\}_{\alpha \in S}$ 
such that the following conditions are satisfied:
\begin{enumerate}[(i)]
\item The objects $U_{\alpha}$ cover $\calX$: that is, there is an effective epimorphism 
$\coprod U_{\alpha} \to \bm{1}$, where $\bm{1}$ denotes the final object of $\calX$.
\item For each $\alpha \in S$, the restriction $\F|_{U_{\alpha}}$ 
is a constant object of $\Shv(\calX_{/U_{\alpha}}; \scrC)$.
\end{enumerate}
\edefn

\subsection{Constructible sheaves}

\bdef[{\cite[Definition A.5.1]{HA}}]
Let $A$ be a partially ordered set. We will regard $A$ as a topological space, where
a subset $U \subseteq A$ is open if it is closed upwards: that is, if $x \leq y$ and $x \in U$ implies that $y \in U$.
Let $X$ be a topological space. An $A$-stratification of $X$ is a continuous map $f: X \to A$. Given
an $A$-stratification of a space $X$ and an element $a \in A$, we let $X_a$, $X_{\leq a}$, $X_{<a}$, $X_{\geq a}$, and $X>a$ denote
the subsets of $X$ consisting of those points $x \in X$ such that $f(x) = a$, $f(x) \leq a$, $f(x) < a$, $f(x) \geq a$,
and $f(x) > a$, respectively. 
\edefn

For a given $a \in A$, the subset $X_a$ is called a \emph{stratum}. 
A collection of such subsets for varying $a \in A$ are called \emph{strata}.

\bdef[{\cite[Definition A.5.2]{HA}}]
Let $A$ be a partially ordered set and let $X$ be a topological space equipped with
an $A$-stratification. We will say that an object $\F \in \Shv(X; \calS)$ is $A$-constructible 
if, for every element $a \in A$, the restriction $F|_{X_a}$ is a locally constant object of $\Shv(X_a; \calS)$. 
Here $\F|_{X_a}$ denotes the image
of $\F$ under the left adjoint to the pushforward functor $\Shv(X_a; \calS) \to \Shv(X; \calS)$.
\edefn

We let $\Shv^A(X; \calS)$ denote the full subcategory of $\Shv(X; \calS)$ 
spanned by $A$-constructible objects.

\subsubsection{Conically stratified spaces}
To ensure the theory of $A$-constructible sheaves
is well-behaved, Lurie introduces a regularity
condition on the stratification $f: X \to A$.
Given such an $A$-stratified space, we may
define an associated $A^{\lhd}$-stratified
space, denoted $C(X)$, called the open
cone of $X$. The definition of $C(X)$
is given in \cite[Definition A.5.3]{HA}.
The regularity condition on stratifications
is formulated using the notion of open cone.

\bdef[{\cite[Definition A.5.5]{HA}}]
Let $A$ be a partially ordered set, let $X$ be an $A$-stratified space, and let $x \in
X_a \subseteq X$ be a point of $X$. We will say that $X$ is conically 
stratified at the point $x$ if there exists an
$A_{>a}$-stratified topological space $Y$, 
a topological space $Z$, and an open embedding $Z \times C(Y) \hook X$
of $A$-stratified spaces whose image $U_x$ contains $x$. Here we 
regard $Z \times C(Y)$ as endowed with the
$A$-stratification determined by the $A^{\lhd}_{>a} \simeq A_{\geq a}$-stratification of $C(Y)$.
\edefn

We say $X$ is conically stratified if it is conically stratified at
every point $x \in X$.

\subsubsection{Whitney stratified spaces}
This are other regularity conditions on stratifications
available to us when $X$ is a complex analytic space
which comprise the notion of a Whitney stratification.

\bdef[{e.g.\ \cite[\S7.1.2]{MacPherson}}]
Let $X$ be a complex analytic space. Then an $A$-stratification
$f: X \to A$ is said to be a Whitney stratification if
\begin{enumerate}[(i)]
\item The strata $X_a$ are smooth manifolds.
\item The stratification is locally finite.
\item The closure of a stratum $\overline{X_a}$ is
a union of strata.
\item Any two strata satisfy Whitney's conditions\footnotemark A and B.
\end{enumerate}
\edefn
\footnotetext{These are easily searchable terms.}

The precise definition of a Whitney stratification
is not terribly relevant for our purposes. It is important,
however, to note that a Whitney $A$-stratified space is
conically stratified stratified. Indeed, around
each point in a Whitney stratified space, there
exist a neighborhood homeomorphic to the product
of a contractible space with a cone over a stratified
sphere (by, for example, \cite[Theorem 7.3]{MacPherson}).

\subsubsection{Exodromy equivalence}
In the previous section, we recalled that if a topological space
$X$ is locally of singular shape, the locally constant objects
of $\Shv(X; \calS)$ are equivalent to the functor category
$\Fun(\Sing(X), \calS)$. There is a similar description of
$A$-constructible sheaves as functors from a subcategory of
$\Sing(X)$.

\bdef
Let $X$ be a paracompact topological space,
and let $f: X \to A$ be a conical stratification of $X$.
Then $\Sing^A(X)$ is an $\infty$-category called
the \emph{exit-path category} of $X$ with respect to
$f: X \to A$. 
\edefn

\bthm[{\cite[Theorem A.9.3]{HA}}]
\label{thm: exit-path description}
Let $X$ be a paracompact topological space which is locally of singular shape
and is equipped with a conical $A$-stratification, where $A$ is a partially ordered set satisfying the
ascending chain condition.\footnotemark 
\footnotetext{We say that a partially ordered set $A$ satisfies the ascending chain condition
if every nonempty subset of A has a maximal element.}
Then there is a functor 
$\Psi_X: \Fun(\Sing^A(X), \calS) \to \Shv^A(X; \calS)$
which induces an equivalence of categories,
\[\Fun(\Sing^A(X), \calS) \xrightarrow{\simeq} \Shv^A(X;\calS).\]
\ethm

The equivalence proven in the above theorem
is usually called the (topological) \emph{exodromy equivalence}.

\brem
Any complex analytic space $X$ is paracompact
and locally of singular shape, and any Whitney
stratification $X \to A$ of it is automatically conical
and satisfies the ascending chain condition.
\erem

\subsubsection{Coefficients in presentable stable $\infty$-categories}
Given an arbitrary $\infty$-category $\scrC$,
the category of $A$-constructible $\scrC$-valued
sheaves on $X$ is defined in much the same way
as above.

\bdef
\label{def: A-constructible sheaves}
Let $X$ be an $A$-stratified topological space.
An object $\F \in \Shv(X;\scrC)$ is defined to be
$A$-constructible if, for every $a \in A$, the restriction
$\F|_{X_a}$ is a locally constant object of $\Shv(X_a;\scrC)$. 
\edefn

We will be interested in the case when $\scrC$ is
a presentable stable $\infty$-category.
In particular, we would like an analogue of
\cref{thm: exit-path description} to hold for $\scrC$-valued
sheaves. Bootstrapping from
\cref{thm: exit-path description} to such an
analogue is difficult.
In their recent preprint, \cite{PortaTeyssier},
Mauro Porta and Jean-Baptiste Teyssier
prove a very general exodromy equivalence for
stable presentable coefficients,
of which the following is a special case
of interest to us.

\bthm[{\cite[Theorem 5.17]{PortaTeyssier}}]
\label{thm: general exodromy}
Let $X$ be a paracompact topological space
locally of singular shape, and let $X \to A$ be
a conical stratification satisfying the ascending chain
condition. Then, 
\[\Shv^A(X;\scrC) \simeq \Fun(\Sing^A(X), \scrC).\]
\ethm

Throughout this work, we will consider
$A$-constructible sheaves valued in certain
presentable $\infty$-categories with a distinguished subcategory
of objects. In these cases, we would like to consider
only those $A$-constructible sheaves whose stalk
at each point of $X$ belongs to this distinguished
subcategory. In order to handle these cases uniformly,
we introduce the following formalism, which is a variant
of the formalism found in \cite{Nori}.

\bdef
\label{def: coefficient pair}
A \emph{coefficient pair} $(\scrN, \scrC)$ consists of 
a stable presentable $\infty$-category $\scrC$,
and a subcategory $\scrN \subset \scrC$ that is stable under
finite limits.
\edefn

To each coefficient pair,
there are two distinct categories 
one might consider, with somewhat
counter-intuitive names and notation.

\bdef
\label{def: weakly constructible sheaves}
Suppose that $X$ is a complex analytic space
and that $(\scrC, \scrN)$ is a coefficient pair.
We denote by $\Shv_{\wc}(X;\scrC)$ the full subcategory
of $\Shv(X;\scrC)$ spanned by sheaves that are $A$-constructible
with respect to \emph{some} Whitney stratification $X \to A$. 
\edefn

We call objects of $\Shv_{\wc}(X; \scrC)$ \emph{weakly constructible sheaves},
and the category itself the category of \emph{weakly constructible $\scrC$-valued
sheaves on $X$}.

\brem
Though an object of $\Shv_{\wc}(X; \scrC)$ is called weakly
constructible, it is, by definition, $A$-constructible for some
Whitney stratification $X \to A$.
\erem

\bdef
\label{def: constructible sheaves}
We denote by $\Shv_c(X; \scrC)$ the full subcategory of
$\Shv_{\wc}(X;\scrC)$ spanned by objects $\F$ such that
$\F_x \in \scrN$ for every $x \in X$.
\edefn

We call objects of $\Shv_c(X; \scrC)$ 
\emph{constructible sheaves} (without reference
to any stratification $X \to A$) with respect to the
coefficient pair $(\scrC, \scrN)$, and the category itself 
the category of \emph{constructible $\scrC$-valued
sheaves on $X$ with respect to $(\scrC, \scrN)$}.

We introduce one more category of sheaves,
which is not typically considered in the literature.

\bdef
We denote by $\Shv_c^{f.s.}(X;\scrC)$ the full
subcategory of $\Shv_c(X;\scrC)$ spanned by objects
$\F$ which are constructible with respect to some
Whitney stratification $X \to A$ where $A$ is a \emph{finite}
partially ordered set. We will sometimes call such objects
\emph{finitely constructible sheaves}.
\edefn

\subsection{Constructible sheaves on finite CW complexes}
Let $(\scrC, \scrN)$ be a coefficient pair, and let
$X$ be a topological space. While the stalks of a
constructible sheaf $\F \in \Shv_c(X; \scrC)$ lie in $\scrN$,
it is not the case that all sections of $\F$ take values in
$\scrN$. Nonetheless, sections of $\F$ over open subsets $U \subseteq X$
of finite homotopy type (i.e. homotopy equivalent to a finite
CW complex) \emph{do} lie in $\scrN$, as shown in the
following proposition.

\bprop
Suppose that $\F \in \Shv_c(X; \scrC)$ is $A$-constructible
with respect to some Whitney stratification $X \to A$, and that
$U \subseteq X$ is an open subset of finite homotopy
type such that the induced stratification on $U$
is finite. Then $\F(U) \in \scrN$.
\eprop

\bproof
By imitating the proof of \cite[Corollary 4.1.8]{Dimca}, we see that there exists a filtration 
$\mathcal{F}_0 \to \cdots \to \mathcal{F}_n = \mathcal{F}|_U$
of length $n = \{\# \, \text{of strata in} \, U\}$ of $\mathcal{F}|_U$ 
by objects of $\Shv_c(U; \scrC)$, such 
that the $i$th associated graded piece is the extension-by-zero of some local system 
$\mathcal{L}_i$ on the stratum $U_i$. The object $\mathcal{F}_0$ is locally constant
with stalks in $\scrN$. Assume that $\F_0(U) \in \scrN$. By induction, 
it then suffices to show that if $\mathcal{F}_i(U) 
\in \scrN$, then $\mathcal{F}_{i+1}(U) \in \scrN$. But $\mathcal{F}_{i+1}(U)$ 
is an extension of $\mathcal{F}_i(U)$ by $\mathcal{L}_i(U)$ and $\scrN$
is closed under extensions, so we win.

Thus, it remains to show that $\F_0(U) \in \scrN$. This is shown
by the classical proof that the singular cohomology of a finite
CW complex is bounded with finite dimensional cohomology
groups, which generalizes immediately to the setting of
local systems with stalks in $\scrN$.
\eproof

This observation suggests that, if $X$ has a basis
of open subsets of finite homotopy type, it might be
possible to produce functors on $\Shv_c(X; \scrC)$ 
from functors on $\scrN$. With this
in mind, we formulate the following convention.

\begin{conv}
\label{conv: finite homotopy type}
Assume that $X$ has finite homotopy type.
Let $\calB(X)$ denote the partially ordered
set associated to a fixed basis of open subsets
$B \subset X$ of finite homotopy type which
includes $X$ itself and is stable under finite intersections.
\end{conv}

\bdef
The partially ordered set $\calB(X)$ has a
natural Grothendieck topology on it.
Let $\Shv(\calB(X); \scrC)$ 
denote the $\infty$-category
of $\scrC$-valued sheaves on $\calB(X)$,
where $\scrC$ is an arbitrary $\infty$-category.
\edefn

\bdef
If $(\scrC, \scrN)$ is a coefficient pair, 
we define the $\infty$-subcategory,
$\Shv_c(\calB(X); \scrC)$ (resp.\ $\Shv_{\wc}(\calB(X); \scrC)$)
(resp.\ $\Shv_c^{f.s.}(\calB(X); \scrC)$) 
$\subset \Shv(\calB(X); \scrC)$
to be the essential image of $\Shv_c(X; \scrC)$
(resp.\ $\Shv_{\wc}(X; \scrC)$) (resp.\ $\Shv_c^{f.s.}(\calB(X); \scrC)$)
under the restriction functor
$\theta: \Shv(X; \scrC) \to \Shv(\calB(X); \scrC)$.
\edefn

Let $X$ and $\calB(X)$ be as in \cref{conv: finite homotopy type},
and let $(\scrC, \scrN)$ and $(\scrC',\scrN')$ be coefficient pairs.
Then any functor $\scrN \to \scrN'$ induces a functor
on basis sheaves,
\[F_{\calB(X)}: \Shv_c(\calB(X); \scrC) \to \Shv_c(\calB(X); \scrC')\] 
by pointwise composition.

\bprop
\label{prop: basis shv = real shv}
Let $X$ and $\calB(X)$ be as in \cref{conv: finite homotopy type},
and let $\scrC$ be an arbitrary category admitting
limits. Then the canonical restriction map,
\[\theta: \Shv(X; \scrC) \to 
	\Shv(\calB(X); \scrC)\]
is an equivalence of $\infty$-categories.
If, furthermore, $(\scrC, \scrN)$ is a coefficient pair,
then $\theta$ restricts to equivalences, 
\begin{align}
\Shv_c(X; \scrC) &\xrightarrow{\simeq} \Shv_c(\calB(X); \scrC) \\
\Shv_c^{f.s.}(X; \scrC) &\xrightarrow{\simeq} \Shv_c^{f.s.}(\calB(X); \scrC).
\end{align}
\eprop

\bproof
We begin by showing that the
restriction of $\infty$-topoi,
\[\theta_{\calS}: \Shv(X; \calS) \to \Shv(\calB(X); \calS),\]
is an equivalence.

The following argument
is based on the discussion
right before \cite[Warning 7.1.1.4]{HTT}.
Observe that the
category $\calB(X)$
has all finite products by 
\cite[Lemma 002O]{Stacks} since
it is closed under fiber products
(intersections) and because
it has a final object, ($X$).
Moreover, it is a small
$0$-category. The proof
of \cite[Proposition 6.4.5.7]{HA} 
shows that the $\infty$-topos, 
$\Shv(\calB(X); 
\mathcal{S})$,
is a $0$-localic. As such, it
is determined by the locale\footnotemark
\footnotetext{\cite[Definition 6.4.2.3]{HTT}}
of subobjects of the final
object $\underline{\ast} \in
\Shv(\calB(X); 
\mathcal{S}),$ which in
turn is equivalent to
the locale given by
the partially ordered set
$\calB(X)$ itself.
For the same reasons, the $\infty$-topos 
$\Shv(X; \calS)$
is determined by the partially
ordered set, $\calU(X)$, of opens in $X$.
Now observe that $\calB(X)$
being a basis means precisely
that the inclusion of partially
ordered sets, $\calB(X)
\subset \calU(X)$ induces
an isomorphism of locales.
It follows that the restriction map
$\theta_{\calS}$ is an equivalence
of $\infty$-topoi.

Now recall \cite[Proposition
I.1.3.1.7]{SAG}, which states that,
for a small $\infty$-category
$\mathcal{T}$ equipped with
a Grothendieck topology and $\scrE$
an arbitrary $\infty$-category
admitting limits,
\[\Shv_{\scrE}(\Shv(\mathcal{T}; 
\mathcal{S})) \xrightarrow{\simeq} \Shv(\mathcal{T}; \scrE),\]
Thus, in light of the equivalence
$\theta_{\mathcal{S}}$, the restriction map,
$\Shv(X; \scrC)
\xrightarrow{\theta'} \Shv(\calB(X); 
\scrC)$, 
is an equivalence.

Finally, by definition, the restriction of a weakly constructible
sheaf $X$ to a sheaf on $\calB(X)$ belongs
to $\Shv_{\wc}(\calB(X); \scrC)$
and the restriction of sheaves to $\calB(X)$ preserves
stalks because $\calB(X) \to \calU(X)$ is cofinal.
It follows that the equivalence $\theta'$ restricts to
a functor $\theta: \Shv_c(X; \scrC) 
\to \Shv_c(\calB(X); \scrC)$ which is
also an equivalence (e.g.\ by fully faithfulness and
essential surjectivity).
\eproof

\Cref{prop: basis shv = real shv} has several immediate
corollaries, which we list below.

\bcor
\label{cor: non-c induces shv functor}
Let $X$ and $\calB(X)$ be as in \cref{conv: finite homotopy type}. 
Let $(\scrC, \scrN)$ be a coefficient pair, and $\scrE$ an arbitrary
$\infty$-category that admits limits. Suppose that $F: \scrN \to \scrE$
is an arbitrary functor. Then $F$ induces a functor
\[F_{\Shv}: \Shv_c^{f.s.}(X; \scrC) \to \Shv(X; \scrE).\]
Moreover, if $\scrE'$ is another $\infty$-category that admits 
limits and $G: \scrE \to \scrE'$ is another functor, 
then $(G \circ F)_{\Shv}$ is naturally isomorphic
to $G_{\Shv} \circ F_{\Shv}$, where $G_{\Shv}: \Shv(X; \scrE) \to
\Shv(X; \scrE')$ is the functor given by the sheafification of
the $\scrE'$-valued presheaf given by pointwise composition with $G$.
\ecor

\bproof
By composition with $F$, we have an induced functor
$F_{\calB(X)}: \Shv_c^{f.s.}(\calB(X); \scrC) \to \PreShv(\calB(X); \scrE)$.
Composition with the sheafification functor $L: \PreShv(\calB(X); \scrE)
\to \Shv(\calB(X); \scrE)$ obtains the functor $L \circ F_{\calB(X)}$, 
which we also denote by $F_{\calB(X)}$.
By \cref{prop: basis shv = real shv}, the restriction functors
$\theta_{\scrC}$ and $\theta_{\scrE}$ are equivalences.
The composition $\theta_{\scrE}^{-1} 
\circ F_{\calB(X)} \circ \theta_{\scrC}$
furnishes the promised functor, $F_{\Shv}$.
Now let $G_{\calB(X)}: \Shv(\calB(X); \scrE) \to
\PreShv(\calB(X); \scrE') \xrightarrow{L} \Shv(\calB(X); \scrE')$ 
be the functor given by pointwise composition with $G$,
followed by sheafification. Clearly,
$(G \circ F)_{\calB(X)}$ is given by the composition
$G_{\calB(X)} \circ F_{\calB(X)}$, so we have 
\begin{align*}
(G \circ F)_{\Shv} 	&= \theta_{\scrE'}^{-1} \circ (G \circ F)_{\calB(X)} \circ \theta_{\scrC} \\ 
				&\simeq \theta_{\scrE'}^{-1} \circ G_{\calB(X)} \circ F_{\calB(X)} \circ \theta_{\scrC} \\
				&\simeq \theta_{\scrE'}^{-1} \circ G_{\calB(X)} \circ \theta_{\scrE} \circ \theta_{\scrE}^{-1} \circ F_{\calB(X)} \circ \theta_{\scrC} \\
				&= G_{\Shv} \circ F_{\Shv}.
\end{align*}
\eproof

\bcor
\label{cor: induces shv functor}
Let $X$ and $\calB(X)$ be as in \cref{conv: finite homotopy type}. 
Suppose that $(\scrC, \scrN)$ and $(\scrC', \scrN')$
are coefficient pairs, and let $F: \scrN \to \scrN'$
be a functor that preserves limits. Then $F$ induces a functor,
\[F_{\Shv}: \Shv_c^{f.s.}(X; \scrC) \to \Shv_c^{f.s.}(X; \scrC').\]
\ecor

\bproof
The proof is almost identical to that of 
\cref{cor: non-c induces shv functor} above.
By composition with $F$, we have an induced functor
$F_{\calB(X)}: \Shv_c^{f.s.}(\calB(X); \scrC) \to \Shv_c^{f.s.}(\calB(X); \scrC')$.
By \cref{prop: basis shv = real shv}, the restriction functors
$\theta_{\scrC}$ and $\theta_{\scrC'}$ are equivalences.
The composition $\theta_{\scrC'}^{-1} 
\circ F_{\calB(X)} \circ \theta_{\scrC}$
furnishes the promised functor, $F_{\Shv}$.
\eproof

\bcor
\label{cor: F_Shv compatible with pullback}
Let $(\scrC, \scrN)$ be a coefficient pair,
and let $F: \scrN \to \scrE$ be a 
functor to an $\infty$-category that admits limits.
Let $X$ and $Y$ be two complex analytic
spaces with bases $\calB(X)$ and $\calB(Y)$
as in \cref{conv: finite homotopy type},
and let $f: X \to Y$ be a map of analytic spaces.
Then $f^*$ preserves constructibility,
and $F_{\Shv}$ commutes with 
$f^*$ in the sense that the following
diagram commutes,
\[\begin{tikzcd}
	{\Shv_c^{f.s.}(Y; \scrC)} & {\Shv(Y; \scrE)} \\
	{\Shv_c^{f.s.}(X; \scrC)} & {\Shv(X; \scrE)} 
	\arrow["{f^*}", from=1-1, to=2-1]
	\arrow["{f^*}", from=1-2, to=2-2]
	\arrow["{F_{\Shv}}", from=1-1, to=1-2]
	\arrow["{F_{\Shv}}", from=2-1, to=2-2]
\end{tikzcd}\] 
\ecor

\bproof
Given any finite stratification of $Y \to B$, there exists a finite
stratification of $X \to A$ such that $f^{-1}(Y_b)$ is a union of
strata in $X$, for all $b \in B$. It is clear that if $\F \in \Shv(Y; \scrC)$
is $B$-constructible, then $f^*\F$ is $A$-constructible. Moreover,
$f^*$ preserves stalks. Thus, $f^*$ indeed restricts to a functor
$f^*: \Shv_c(Y; \scrC) \to \Shv_c(X;\scrC)$.
 
In order to prove that $F_{\Shv}$ commutes with $f^*$,
it suffices to show the commutativity of
\[\begin{tikzcd}
	{\Shv(\calB(Y); \scrN)} & {\PreShv(\calB(Y); \scrE)} \\
	{\Shv(\calB(X); \scrN)} & {\PreShv(\calB(X); \scrE)}
	\arrow["{{f^*}_{\scrN}}"', from=1-1, to=2-1]
	\arrow["{{f^*}_{\scrE}}"', from=1-2, to=2-2]
	\arrow["{F_{\calB(Y)}}", from=1-1, to=1-2]
	\arrow["{F_{\calB(X)}}", from=2-1, to=2-2].
\end{tikzcd}\]
Both ${f^*}_{\scrN}$ and ${f^*}_{\scrE}$ are induced
by the underlying geometric morphism of $\infty$-topoi,
$f^*: \calB_Y := \Shv(\calB(Y); \calS) \to \Shv(\calB(X); \calS) =: \calB_X$ via
composition with the limit-preserving functor ${f_*}^{\op}: 
\calB_X^{\op} \to \calB_Y^{\op}$.
As such, the commutativity of the above diagram comes
from its identification with the following one, which is canonically
commutative:
\[\begin{tikzcd}
	{\RFun(\calB_Y^{\op}, \scrN)} & {\Fun(\calB_Y^{\op}, \scrE)} \\
	{\RFun(\calB_X^{\op}, \scrN)} & {\Fun(\calB_X^{\op}, \scrE)}
	\arrow["{F \circ -}", from=1-1, to=1-2]
	\arrow["{F \circ -}", from=2-1, to=2-2]
	\arrow["{- \circ {f_*}^{\op}}", from=1-2, to=2-2]
	\arrow["{- \circ {f_*}^{\op}}"', from=1-1, to=2-1].
\end{tikzcd}\]
\eproof

\bcor
\label{cor: induces natural iso}
Let $X$ and $\calB(X)$ be as in \cref{conv: finite homotopy type}. 
Suppose that $(\scrC, \scrN)$ and $(\scrC', \scrN')$
are coefficient pairs, and that $F, \wit{F}: \scrN \to \scrN'$
are limit-preserving functors between 
distinguished subcategories. If $F$ and $\wit{F}$ are
naturally isomorphic, then $F_{\Shv} \simeq \wit{F}_{\Shv}$.
\ecor

\bproof
Suppose that $\calN_{F \to \wit{F}}: F \xrightarrow{\simeq} \wit{F}$
is a natural isomorphism. Then $\theta_{\scrC'}^{-1}
\circ \calN_{F \to \wit{F}} \circ \theta_{\scrC}$ is a natural
isomorphism $F_{\Shv} \xrightarrow{\simeq} \wit{F}_{\Shv}$.
\eproof

\subsection{Constructible sheaves of local systems}
As noted earlier in this work, the pair $(\Vect, \Perf)$
is a coefficient pair. Moreover, $\left(\Fun(S^1, \Vect), \Fun(S^1, \Perf)\right)$
is clearly a coefficient pair, as well.\footnotemark
\footnotetext{Note, for example, that $\Fun(S^1, \Vect) \simeq \calD(\C[t,t^{-1}])$.}
The goal of this section is to produce an equivalence,
\beqn
\label{eqn: S1 commutes with Shv}
\Fun(S^1, \Shv_c^{f.s.}(X;\C)) \simeq \Shv_c^{f.s.}(X; \Fun(S^1, \Vect)).
\eeqn

\blem
\label{lem: limit of coefficients}
Let $I$ be a small category, and
suppose that $F: I \to \Pr^R$ is an $I$-shaped diagram
of presentable categories $\scrC_{i \in I}$ such
that its structure morphisms reflect
limits. Let $\calX$ be an arbitrary $\infty$-topos. Then there
is an equivalence $\Shv(\calX; \lim_I \scrC_i) \xrightarrow{\simeq}
\lim_I \Shv(\calX; \scrC_i)$.
\elem

\bproof
Let $\scrC$ be an arbitrary $\infty$-category.
Recall that $\Shv(\calX; \scrC)$ is defined as
the fully subcategory of $\Fun(\calX^{\op}, \scrC)$
on functors which preserve limits. Since small
limits in $\Pr^R$ exist, all of the structure morphisms
$F_i: \lim_I \scrC_i \to \scrC_i$ preserve limits, as well.
Composition with $F_i$ therefore induces morphisms
$\Shv(\calX; \lim_I \scrC_i) \to \Shv(\calX; \scrC_i)$,
so we obtain, by the universal property, a functor
$\Shv(\calX; \lim_I \scrC_i) \to \Shv(\calX; \scrC_i)$.
On the other hand, $\lim_I \Shv(\calX; \scrC_i)$ is
a full subcategory of $\Fun(\calX; \lim_I \scrC_i)$.\footnotemark
\footnotetext{Recall that $\Fun(-,-)$ preserves limits in the second
variable.}
Because the structure morphisms of $F$ reflect limits,
functors in the subcategory $\lim_I \Shv(\calX; \scrC_i)$
necessarily preserve limits. Thus, $\lim_I \Shv(\calX; \scrC_i)
\subset \Shv(\calX; \lim_I \scrC_i)$. It is clear that this
inclusion is inverse to the universal morphism above.
\eproof

\blem
\label{lem: weakly constructible}
Let $I$ be a finite category, and
suppose that $F: I \to \Pr^R_{st}$ is an $I$-shaped diagram
of presentable stable $\infty$-categories $\scrC_{i \in I}$ such
that its structure morphisms reflect
limits. Let $X$ be a complex analytic space.
Then there is an equivalence $\Shv_{\wc}(X; \lim_I \scrC_i)
\xrightarrow{\simeq} \lim_I \Shv_{\wc}(X; \scrC_i)$.
\elem

\bproof
By \cref{lem: limit of coefficients}, we have an equivalence,
$\Shv(X; \lim_i \scrC_i) \xrightarrow{\simeq} \lim_I \Shv(X; \scrC_i)$.
It therefore suffices to show that the restriction of this equivalence
to $\Shv_{\wc}(X; \lim_i \scrC_i)$ is an equivalence. The essential
image of this restriction is obviously contained in $\lim_I \Shv(X; \scrC_i)$.\footnotemark
\footnotetext{The composition of a sheaf $\calU(X)^{\op} \to \scrC$, constructible
with respect to a particular stratification of $X$, with a functor $\scrC \to \scrD$
is constructible with respect that same stratification.}

It remains to show that an object of $\lim_I \Shv_{\wc}(X;\scrC_i)$
is weakly constructible as an object of $\Shv(X;\lim_I \scrC_i)$
under the obvious inclusion.
For each $i \in I$, the sheaf $\F_i$ is
constructible with respect to some Whitney stratification
$f_i: X \to A_i$, where $A_i$ is a partially
ordered set. Because $I$ is finite, this collection
of stratifications admits a common refinement
$f: X \to A$ with respect to which $\F_i$ is constructible.
Thus, it suffices to show that an object of
$\lim_I \Shv^A(X; \scrC_i)$ is $A$-constructible
under the inclusion $\lim_I \Shv^A(X; \scrC_i) \hook
\Shv(X; \lim_I \scrC_i)$.
By \cref{thm: general exodromy} 
(i.e.\ \cite[Theorem 5.17]{PortaTeyssier}), $\Shv^A(X; \lim_I \scrC_i) \simeq
\Fun(\Sing^A(X), \lim_I \scrC_i)$, where $\Sing^A(X)$ denotes
the exit-path category. It follows that
\begin{align*}
\lim_I \Shv^A(X; \scrC_i) &\simeq \lim_I \Fun(\Sing^A(X), \scrC_i) \\
					&\simeq \Fun(\Sing^A(X), \lim_I \scrC_i) \\
					&\simeq \Shv^A(X; \lim_I \scrC_i) \\
					&\subset \Shv(X; \lim_I \scrC_i).
\end{align*}
Moreover, after unraveling the
definitions of the functors involved, it 
is not hard to see that the composition of the chain of
functors above is the naturally the embedding
$\lim_I \Shv^A(X;\scrC_i) \hook \Shv(X; \lim_I \scrC_i)$,
so we are done.
\eproof

Before stating the next lemma, we
introduce the following definition.

\bdef
Let $(\scrC, \scrN)$ and $(\scrC', \scrN')$ be
coefficient pairs. Any functor $F: \scrC \to \scrC'$
such that the image of $\scrN$ is contained in $\scrN'$
is called a \emph{functor of coefficient pairs}.
\edefn

\blem
Suppose that $I$ is a finite category, and let
$F: I \to \Pr^R_{st}$ is an $I$-shaped diagram
such that the structure morphisms of $F$ are
functors of coefficient pairs $(\scrC_i, \scrN_i)_{i \in I}$.
Then the pair $(\lim_I \scrC_i, \lim_I \scrN_i)$
is a coefficient pair.
\elem 

\bproof
Since $\Pr^R_{st}$ admits limits, $\lim_I \scrC_i$
is a presentable stable $\infty$-category. It suffices therefore
to show that $\lim_I \scrN_i$ is a subcategory of $\scrC_i$
which is stable under finite limits and colimits. By universality,
we have a functor $\lim_I \scrN_i \to \lim_I \scrC_i$. Note that
the mapping space between two objects $a, b$ in an arbitrary
$\infty$-category $\scrC$ is given by 
$\Fun(\Delta^1, \scrC) \times_{\Fun(\partial \Delta^1, \scrC)} (a,b)$.
The functor $\lim_I \scrN_i \to \lim_I \scrC_i$ induces
a functor $\Fun(\Delta^1, \lim_I \scrN_i) \to \Fun(\Delta^1, \lim_I \scrC_i)$
by composition. Since $\Fun(-,-)$ preserves limits in the second
variable, we obtain a functor $\lim_I \Fun(\Delta^1, \scrN_i)
\to \lim_I \Fun(\Delta^1, \scrC_i)$. This functor is
the universal one induced by the morphism of $I$-shaped diagrams
whose terms are $\Fun(\Delta^1, \scrN_i)$ and $\Fun(\Delta^1, \scrC_i)$,
respectively. Since the functors $\Fun(\Delta^1, \scrN_i) \times_{\Fun(\partial\Delta^1, \scrN_i)} 
(a_i, b_i) \to \Fun(\Delta^1, \scrC_i) \times_{\Fun(\partial\Delta^1, \scrC_i)} (a_i,b_i)$ 
induced by the structure morphisms 
are equivalences by the assumption
that $\scrN_i \subset \scrC_i$ is fully faithful, we obtain that
\[\Fun(\Delta^1, \lim_I \scrN_i) 
	\times_{\Fun(\partial \Delta^1, \lim_I \scrN_i)} (a,b) 
		\to \Fun(\Delta^1, \lim_I \scrC_i) 
			\times_{\Fun(\partial \Delta^1, \lim_I \scrC_i)} (a,b),\]
for $a,b \in \lim_I \scrN_i$,
is an equivalence as well, showing that
the universal functor $\lim_I \scrN_i
\to \lim_I \scrC_i$ is fully faithful, as well.
\eproof

\blem
\label{lem: constructible}
Let $I$ be a finite category, and
suppose that $F: I \to \Pr^R_{st}$ is an $I$-shaped diagram
of presentable stable $\infty$-categories $\scrC_{i \in I}$ with
subcategories $\scrN_i \subset \scrC_i$ such that
$(\scrC_i, \scrN_i)$ is a coefficient pair. Assume that
the structure morphisms of $F$ are functors
of coefficient pairs and reflect
limits. Denote by $(\scrC, \scrN)$ the coefficient
pair $(\lim_I \scrC_i, \lim_I \scrN_i)$.
Then if $X$ is a complex analytic space,
there is an equivalence $\Shv_c(X; \scrC)
\xrightarrow{\simeq} \lim_I \Shv_c(X; \scrC_i)$.
\elem


\bcor
\label{cor: S1 commutes with Shv}
Let $(\scrC, \scrN)$ be a coefficient pair.
There is an equivalence of categories,
\[\Fun(S^1, \Shv_c^{f.s.}(X; \scrC)) \simeq \Shv_c^{f.s.}(X; \Fun(S^1, \scrC)),\]
where the latter category denotes constructible
sheaves on $X$ with respect to the coefficient pair
$(\Fun(S^1, \scrC), \Fun(S^1, \scrN))$.
\ecor

\bproof
Recall that $S^1 \simeq \colim_{S^1} \ast$. Because
$\Fun(-,-)$ takes colimits in the first variable to limits,
we have $\Fun(S^1, \Shv_c(X;\C)) \simeq \lim_{S^1} \Shv_c(X;\C)$.
On the other hand, we have 
\begin{align*}
\Shv_c(X; \Fun(S^1, \Vect)) &\simeq \Shv_c(X; \lim_{S^1} \Vect) \\ 
							&\simeq \lim_{S^1} \Shv_c(X; \Vect),
\end{align*}
where we have used \cref{lem: constructible}.
Restricting the equivalence 
\[\Shv_c(X; \Fun(S^1, \scrC)) \xrightarrow{\simeq} \Fun(S^1, \Shv_c(X; \scrC))\]
obtained in this manner
to the subcategory $\Shv_c^{f.s.}(X; \Fun(S^1, \Vect))$
furnishes the desired equivalence.
\eproof

The desired equivalence \labelcref{eqn: S1 commutes
with Shv} now follows from \cref{cor: S1 commutes
with Shv} by taking the coefficient pair 
$(\scrC, \scrN)$ to be $(\Vect, \Vect^b)$.

\section{Derived local systems on $\bfAone \setminus 0$}
\label{sec: Appendix B}
We gather some results on derived local systems on $\bfAone \setminus 0$,
i.e. objects of $\calD^b(\Loc(\bfAone \setminus 0))$.

\bdef
Let $X$ be a complex algebraic variety or
complex analytic space. We denote by $\lLoc(X)$ the
full subcategory of $\calA(X;\C)$ spanned
by locally constant sheaves of finite rank.
\edefn

\brem
\label{rem: comparison of Locs}
If $\dim X = n$, shift the cohomological
degree by $n$ gives an equivalence
of categories $(-)[n]: \lLoc(X) \xrightarrow{\simeq} \Loc(X)$.
\erem

Let $\calD^b_{\lLoc}(X; \C)$ denote the
full subcategory of $\Shv_c(X;\C)$
spanned by objects whose cohomology sheaves
are objects of $\lLoc(X)$. The fully faithful
embedding $\lLoc(X) \hook \calA(X;\C)$
induces a functor $\calD^b(\lLoc(X)) \to \calD^b_{\lLoc}(X;\C)$.

\blem
\label{lem: beeg lemma}
The functor $\calD^b(\lLoc(\Aone \setminus 0)) 
\to \calD^b_{\lLoc}(\Aone \setminus 0; \C)$
is a t-exact equivalence of stable $\infty$-categories.
\elem

\bproof
The image of the functor is clearly generated by
local systems, so it suffices to show that the functor
is fully faithful. In order to show that it is fully faithful,
it suffices to show that for any $\calL \in \lLoc(\Aone \setminus 0)$,
the functor $\Ext^q_{\Shv(\Aone \setminus 0; \C)}(\calL, -)|_{\lLoc(\Aone \setminus 0)}$
is effaceable for every $q >0$. In order to do so, we
follow the strategy employed by Nori in his proof of
\cite[Theorem 3]{Nori}.

We introduce the following notation. Let $\Delta: \Aone \to \bbA^2$ be the diagonal
map, $\pr_i: \bbA^2 \to \Aone$ be the canonical projections,
and $j: \Aone \setminus 0 \hook \Aone$ and $i: 0 \hook \Aone$
be the canonical inclusions. For the remainder of the proof,
all functors are underived unless otherwise indicated. 
Now let $\scrL \in \lLoc(\Aone \setminus 0)$, and 
define $\calG$ by the short exact sequence,
\beqn
\label{eqn: ses}
0 \to \calG \to \pr_1^*j_!\scrL \to \Delta_*j_!\scrL \to 0.
\eeqn
The derived functor of ${\pr_2}_*$ induces a morphism,
$j_!\scrL = {\pr_2}_*\Delta_*j_!\scrL \to \F := R^1{\pr_2}_*\calG$,
which an monomorphism of sheaves by \cite[Proposition 2.2]{Nori}.
Applying the exact functor $j^*$ to this morphism, we obtain
a morphism,
\[u: \scrL \simeq j^*j_!\scrL \to j^*\F.\]
We claim that $j^*\F \in \lLoc(\Aone \setminus 0)$,
$u$ is a monomorphism,
and that $H^q(\Aone \setminus 0; j^*\F) = 0$ for all $q \geq 0$.

In order to show that $j^*\F$ is a local system, we use
\cite[Remark 1.5]{Nori}. Let $V := \pr_1^{-1}(0) \cup 
\Delta(\Aone) \subset \bbA^2$. The
restriction $\pr_2|_{V}$ is clearly a finite surjective morphism.
From the definition of $\calG$, we see that 
the restriction of $\calG$ to the complement of $V$ is
locally constant, and $\calG|_V = 0$. We may also express
$\Aone$ as the disjoint union of $X_0 = 0$ and $X_1 = \Aone \setminus 0$,
so that $\pr_2: (\pr_1^{-1}X_i \cap V)_{\red} \to X_i$, for $i =0,1$, is a finite
\'etale morphism. Then by \emph{loc. cit.}, $j^*\F
\:= R^1{\pr_2}_*\calG|_{X_1}$ is locally constant. 
By \cite[Proposition 2.2 (3)]{Nori}, $\F$ is constructible,
so its stalks are finite dimensional. Since $j^*$ preserves
stalks, it follows that $j^*\F$ has finite dimensional stalks
and therefore belongs to $\lLoc(\Aone \setminus 0)$, as desired.

By \cite[Proposition 2.2 (1)]{Nori}, the morphism $j_!\scrL
\to \F$ is a monomorphism of sheaves. It follows from 
the exactness of the functor $j^*$ that $u$ is a monomorphism.

In order to prove that $H^q(\Aone \setminus 0; j^*\F) = 0$ for all $q \geq 0$,
we imitate the proof of \cite[Proposition 2.2 (2)]{Nori}. Let $\jmath:
\bbA^2 \setminus \pr_2^{-1}(0) \hook \bbA^2$ be the canonical inclusion. Applying
the exact functor $\jmath^*$ to \labelcref{eqn: ses} 
and using smooth base change, we obtain
another short exact sequence,
\beqn
\label{eqn: ses 2}
0 \to \jmath^*\calG \to {\pr'_1}^*\scrL \to \Delta'_*\scrL \to 0,
\eeqn
where $\pr'_1$ denotes ${\pr_1}|_{\bbA^2 \setminus \pr_2^{-1}(0)}$
and $\Delta'$ denotes $\Delta|_{\Aone \setminus 0}$.
From the long exact sequence obtained by 
Applying the derived functor $R{\pr'_1}_*$ to  
\labelcref{eqn: ses 2}, we obtain the long exact sequence,
\[0 \to {\pr'_1}_*\jmath^*\calG \to \cdots \to R^q{\pr'_1}_*{\pr'_1}^*\scrL \to R^q{\pr'_1}_*\Delta'_*\scrL \to R^{q+1}{\pr'_1}_*\jmath^*\scrL \to \cdots\]
We note, however, that $R{\pr'_1}_*\Delta'_* = R(\pr'_1 \circ \Delta')_* = {\id_{\Aone \setminus 0}}_*$,
and that ${\pr'_1}_*{\pr'_1}^*$ is exact because there
exists a homotopical isomorphism over $\Aone$ from
$\bbA^2 \setminus \pr_2^{-1}(0)$ to a trivial circle bundle
over $\Aone$, whose projection map is clearly proper.
Consequently, we see that $R^q{\pr_1}_*(\jmath^*\calG) = 0$ for all 
$q>0$. From the Leray spectral sequence applied
to the terminal map, $\pt = \pt \circ \pr'_1$, we
see that $R\pt_*(\jmath^*\G) = 0$ for all $q>0$.

On the other hand, by \cite[Proposition 1.3A]{Nori}, 
we have that $R^q{\pr'_2}_*(\jmath^*\calG) = 0$
for $q>1$. Another application of the Leray spectral
sequence yields an equality, $R^q\pt_*R^1{\pr'_2}_*(\jmath^*\calG) =
R^{q+1}\pt_*(\jmath^*\calG) =0$ for all $q \geq 0$.
Now, using \cite[Corollary 1.3B]{Nori}, we
note that $R^1{\pr'_2}_*(\jmath^*\calG) \simeq
j^*R^1{\pr_2}_*\calG =: j^*\F$, so we have shown
that $H^q(\Aone \setminus 0; j^*\F) = 0$ for $q \geq 0$
as desired.

Altogether, this shows that the functor 
${H^q(\Aone \setminus 0; -)}|_{\lLoc(\Aone \setminus 0)}$
is effaceable for all $q \geq 0$. By \cite[Remark 3.8]{Nori},
we see that the functor ${H^q(\Aone \setminus 0; \sExt^p(\calL, -))}|_{\lLoc(\Aone \setminus 0)}$
is also effaceable for all $p,q \geq 0$. Because $\calL$ is locally
constant sheaf with projective (i.e. finite dimensional)
stalks, $\sExt^q(\calL, \calH) = 0$ for all $q>0$
and for any sheaf $\calH$. Consequently,
the local-to-global spectral sequence relating
$\Ext$ and $\sExt$ gives a natural isomorphism, 
\[{\Ext^q_{\Shv(\Aone \setminus 0; \C)}(\calL, -)}|_{\lLoc(\Aone \setminus 0)}
	\simeq H^q(\Aone \setminus 0; {\sHom(\calL, -)}|_{\lLoc(\Aone \setminus 0)}).\]
The proof is concluded by observing 
that $\sHom(\calL, -) = \sExt^0(\calL,-)$.
\eproof

Note that the embedding $\Loc(\Aone \setminus 0) \hook \Perv(\Aone \setminus 0)$
similarly induces a functor $\calD^b(\Loc(\Aone \setminus 0)) \to \Shv_{\Loc}(\Aone \setminus 0; \C)$,.
where now $\Shv_{\Loc}(\Aone \setminus 0; \C)$ denotes the full
subcategory of $\Shv_c(\Aone \setminus 0; \C)$ spanned by objects
whose \emph{perverse} cohomology sheaves belong to $\Loc(\Aone \setminus 0)$.
The following corollary follows easily from \cref{lem: beeg lemma} and 
\cref{rem: comparison of Locs} above.

\bcor
\label{cor: beeg corollary}
The functor $\calD^b(\Loc(\Aone \setminus 0)) 
\to \calD^b_{\Loc}(\Aone \setminus 0; \C)$
is a t-exact equivalence of stable $\infty$-categories,
where the t-structure on the right-hand category is
induced from the perverse t-structure on 
$\Shv_c(\Aone \setminus 0; \C)$.
\ecor

\brem
Of course, if $X$ is a complex algebraic variety,
then $\Loc(X) \simeq \Loc(X^{\an})$ canonically,
so the results above hold \emph{mutatis mutandis}
for $\bfAone \setminus 0$ in place of $\Aone \setminus 0$.
\erem

\blem
\label{lem: locally const}
Let $X$ be a complex analytic space.
The stable $\infty$-category $\calD^b_{\lLoc}(X;\C)$ is
equivalent to the $\infty$-category of locally
constant $\Perf$-valued sheaves,
as defined by Lurie in \cite[Appendix A.1]{HA}.
\elem

\bproof
The category of locally constant $\Perf$-valued
sheaves is a full subcategory of $\calD^b(X;\C)$.
Moreover, it is clear that the cohomology
objects of a locally constant $\Perf$-valued
sheaf are local systems,
since the pullback of sheaves is left exact.
Thus, the category of locally constant $\Perf$-valued
sheaves embeds fully faithfully into $\calD^b_{\lLoc}(X;\C)$.
It remains to show that this embedding is essentially
surjective. Let $\F \in \calD^b_{\lLoc}(X;\C)$, for any
point $x \in X$, there exists a neighborhood of $x$, $U_i(x) \subset X$
such that ${(H^i\F)}|_{U_i(x)}$ is constant. If $H^i\F \simeq 0$,
we let $U_i(x)$ be $X$ itself.
Let $U(x) := \cap_{i \in \bbZ} U_i(x)$. Since
$H^i\F$ is non-zero for only finitely many
$i \in \bbZ$, $U(x)$ is an open subset of $X$.
Since any object in $V \in \Vect$ is equivalent
to the direct such of its cohomologies $\oplus_i H^iV[i]$,
it follows that $\F|_{U(x)}$ is constant, as well. Since
$x$ was arbitrary, we conclude that $\F$ is locally
constant, which concludes the proof. 
\eproof

\blem
\label{lem: pairs}
Suppose that $X$ is a topological space
locally of singular shape.
Then $\calD^b_{\lLoc}(X; \C)$
is equivalent as a stable $\infty$-category to the category of
functors $\Fun(\Sing(X), \Perf)$.
\elem

\bproof
Since $X$ is locally of singular shape, the
$\infty$-topos associated to $X$, $\Shv_{\calS}(X)$,
is locally of constant shape, and the shape of
$X$ can be identified with the space $\Sing(X)$.
By \cite[Theorem A.1.15]{HA}, the fully 
$\infty$-subcategory of locally constant
objects of $\Shv_{\calS}(X)$ is equivalent to
the overcategory $\calS_{/\Sing(X)}$. 
Locally constant sheaves valued in $\Perf$
are given by $\Perf$-valued
sheaves on the $\infty$-topos $\calS_{/\Sing(X)}$,
i.e. the full subcategory of 
$\Fun({\calS_{/\Sing(X)}}^{\op}, \Perf)$ 
spanned by functors which preserve limits.
On the other hand, $\calS_{/\Sing(X)}$ is generated
under small colimits by $\Sing(X)$, so
$\Shv_{\Perf}(\calS_{/\Sing(X)}) \simeq
\Fun(\Sing(X), \Perf)$. Now use
\cref{lem: locally const} to conclude.
\eproof

In particular, if $X$ is an Eilenberg-MacLane space, $K(\bbZ, 1)$,
$\Fun(\Sing(X), \Perf)$ identifies with $\Fun(S^1, \Perf)$, 
where $S^1$ here is the space $S^1$ treated as an 
$\infty$-category with a single object. Informally,
objects of $\Fun(S^1, \Perf)$ are pairs
of an object $E \in \Perf$ and an automorphism 
$T \in \End_{\Vect_{\C}}(E)$.

\brem
Both \cref{lem: locally const} and \cref{lem: pairs} hold
with $\calD^b_{\Loc}(X;\C)$ in place of $\calD^b_{\lLoc}(X;\C)$.
\erem

\section{Filtered objects}
\label{sec: Appendix C}

\subsection{Filtered objects in abelian categories}
Our presentation in this section and the next cleaves very
closely to that found in \cite{GwilliamPavlov}, except 
that we modify some of their definitions to account
for the fact that we work exclusively with decreasing filtrations
in this work. To this point, the presentation for filtered
objects in abelian categories is separate from presentation for 
filtered objects in stable $\infty$-categories. Briefly: the
former are defined as sequences of monomorphisms in
an abelian category, while the latter are defined as the localization
of the category of sequences in a stable $\infty$-category at
a certain collection of morphisms. In particular, the maps of
a filtration in the stable $\infty$-categorical setting are not
required to be monomorphisms.

\subsubsection{}

Let $\calJ := \bbZ^{\op}$, where $\bbZ$ is the linearly 
ordered set of integers regarded as an ordinary category.
Let $\calA$ denote an abelian category.

\bdef[{\cite[Definition 1.2]{GwilliamPavlov}}]
\label{def: Seq(A)}
The category of \emph{sequences in $\calA$} is the functor
category $\Fun(\calJ, \calA)$, which we denote by
$\Seq(\calA)$.
\edefn

\bdef[{\cite[Definition 1.3]{GwilliamPavlov}}]
\label{def: Fil(A)}
The \emph{filtered category of $\calA$}, denoted
$\Fil(\calA)$ is the full subcategory of $\Seq(\calA)$ spanned by
functors $F: \calJ \to \calA$ satisfying that condition that $F(m \to n)$
is a monomorphism for every $n \leq m$. Given an filtered object
$F: \calJ \to \calA$, we will say that $F$ is \emph{a filtration} of
$F(\infty) := \colim_{\calJ} F(n)$.
\edefn

Thus, the only filtrations we consider here are ``exhaustive" filtrations
in usual terminology. We will often denote $F(n)$ by $F_n$ and call it 
the $n$th filtered piece of $F_{\infty}$. For each $n \in \calJ$ and
$F \in \Fil(\calA)$, we let $\gr_n F := \coker F(n+1) \to F(n)$.
We let $\gr: \Fil(\calA) \to \prod_{n \in \calJ} \gr_n$ denote the
\emph{associated graded} functor which sends a filtered object
$F \in \Fil(\calA)$ to the graded object $\prod_{n \in \calJ} \gr_n F$
in $\calA$.

The filtered category $\Fil(\calA)$, is additive and admits kernels and
cokernels, but it is famously not an abelian category in general.
Though derived categories of such categories can still be considered
(see e.g. \cite{SchapiraSchneiders}), it is more natural in this
context to consider the so-called ``filtered" derived category.
Let $\Ch(\calA)$ denote the abelian category of unbounded chain
complexes in $\calA$.

\bdef[{\cite[Definition 1.7]{GwilliamPavlov}}]
\label{def: filtered derived category}
The \emph{filtered derived category} of $\calA$, denoted $D^{\fil}(\calA)$, 
is the localization of $\Ch(\Fil(\calA))$ ($\simeq \Fil(\Ch(\calA))$) with respect 
to the collection filtered weak equivalences, which are maps of sequences $f : F \to G$
such that $\gr(f): \gr F \to \gr G$ is an indexwise quasi-isomorphism.
\edefn

\subsection{Filtered objects in stable $\infty$-categories}
By abuse of notation, we let $\calJ := \bbZ^{\op}$, where $\bbZ$ is the linearly 
ordered set of integers regarded as an $\infty$-category.

\bdef[{c.f. \cite[\S2.2]{GwilliamPavlov}}]
\label{def: Seq(C)}
Let $\scrC$ be a stable $\infty$-category. 
We define $\Seq(\scrC)$ to be the functor
category $\Fun(\calJ, \scrC)$. An object of $\Seq(\scrC)$
is called a \emph{sequence} in $\scrC$.
\edefn

To each $n \in \calJ$, there is a functor $\ev_n: \Seq(\scrC)
\to \scrC$ given by evaluation at $n$. For $F \in \Seq(\scrC)$,
we denote the composition $\ev_n \circ F$ either by $F_n$ or $F(n)$.
For each $n \in \calJ$ and $F \in \Seq(\scrC)$, we let $\gr^n F$ denote
the cofiber of the map $F_{n+1} \to F_n$.

\bdef[{\cite[Definition 2.1]{GwilliamPavlov}}]
The \emph{associated graded functor} $\gr: \Seq(\scrC) \to \prod_{\bbZ} \scrC$
is the functor sending a sequence $F$ to the graded object $\{\gr^n F\}_{n \in \bbZ}$.
We denote by $\pr_n: \prod_{\bbZ} \scrC \to \scrC$ the canonical projection
onto the $n$th factor. Clearly, $\gr^n \simeq \pr_n \circ \gr^{\bullet}$.
\edefn

We say that a morphism of sequences, $f: F \to G \in \Mor(\Seq(\scrC))$,
is a graded equivalence if $\gr(f)$ is an equivalence in $\scrC$. Let
$\calW_{\gr}$ denote the collection of graded equivalences in $\Seq(\scrC)$.

\bdef[{\cite[Definition 2.3]{GwilliamPavlov}}]
\label{def: Fil(C)}
Let $\scrC$ be a stable $\infty$-category.
The \emph{filtered $\infty$-category} of $\scrC$ is defined to be
the localization $\Seq(\scrC)[\calW_{\gr}^{-1}]$, which we denote
by $\Fil(\scrC)$. We will say that an object $F \in \Fil(\scrC)$ is a 
filtration of the object $\lim_{\calJ} F \in \scrC$.
\edefn

As shown in \cite[Lemma 2.15]{GwilliamPavlov}, if it exists, $\Fil(\scrC)$
may be identified with the full subcategory of $\calW_{\gr}$-local 
objects in $\Seq(\scrC)$, meaning that any functor on $\Seq(\scrC)$
restricts to a functor on $\Fil(\scrC)$. In particular, we if a sequence
$F$ lies in $\Fil(\scrC)$, we call $F_n$ or $F(n)$, the \emph{$n$th
filtered piece} of $F$.

We have reason in this work to consider sequences
whose morphisms are monomorphisms,\footnotemark like in the case
of an abelian category.
\footnotetext{In the sense of \cite[\S5.5.6]{HTT}. 
Taken from \textit{loc. cit.}: a monomorphism is a morphism
$f: C \to D$ which is $(-1)$-truncated as an object of the mapping space; this
is equivalent to the assertion that the functor $\scrC_{/f} \to \scrC_{/D}$ is fully faithful.}

\bdef
We denote by $\Fil^{\mon}(\scrC)$ the full subcategory of $\Fil(\scrC)$ spanned by
sequences whose structure morphisms are monomorphisms.
\edefn

The following two theorems of Gwilliam--Pavlov guarantee the existence
of the localization at $\calW_{\gr}$ for a large class of stable $\infty$-categories,
and illuminate the connection of their construction with the classical
filtered derived category.

\bthm[{\cite[Theorem 2.5]{GwilliamPavlov}}]
If $\scrC$ is a stable $\infty$-category that 
admits sequential limits, then $\Fil(\scrC)$ exists and
is a stable $\infty$-category.
\ethm

\bthm[{\cite[Theorem 2.6]{GwilliamPavlov}}]
Let $\calA$ be a Grothendieck abelian category.
Then the homotopy category of $\Fil(\calD(\calA))$
is equivalent to the classical filtered derived 
category $D^{\fil}(\calA)$.
\ethm

With this in mind, we called $\Fil(\calD(\calA))$ the
filtered derived \emph{$\infty$-category} of $\calA$.

\subsubsection{Filtered categories of presentable categories}
Suppose that $\scrC$ is additionally presentable. Then $\Seq(\scrC)$
is a presentable as well, by \cite[Proposition 5.5.3.6]{HTT}. Moreover,
$\Fil(\scrC)$ is an accessible left exact localization of $\Seq(\scrC)$ by
\cite[Proposition 2.14]{GwilliamPavlov} by the parenthetical remark 
contained in the proof of \cite[Theorem 3.9]{GwilliamPavlov}. Altogether,
this shows the following lemma.

\blem
Suppose that $\scrC$ is a presentable stable $\infty$-category.
Then $\Fil(\scrC)$ is also a presentable stable $\infty$-category.
\elem

\subsubsection{Symmetric monoidal structure on filtered categories}
It turns out that if $\scrC$ is a presentable, closed symmetric
monoidal stable $\infty$-category, whose monoidal product we
denote by $\otimes$, then $\Seq(\scrC)$ is as well under the
Day convolution product, which we denote by $\star$. Concretely,
given $F, G \in \Seq(\scrC)$, the $n$th pieces of the Day convolution
is given by
\[F \star G(n) \simeq \colim_{p+q \geq n} F(p) \otimes G(q).\]
We direct the reader to \S2.23 of \cite{GwilliamPavlov} for more
details.

The symmetric monoidal structure on $\Seq(\scrC)$ induces
a closed symmetric monoidal structure on $\Fil(\scrC)$ by completion.
Let $\oblv: \Fil(\scrC) \to \Seq(\scrC)$ denote the forgetful functor.
This has a left adjoint called the \emph{completion} functor which
we denote by $\comp$, which is the homotopical counterpart
of the classical completion of filtrations. We recall the following
definition/theorem.

\bdef[{\cite[Theorem 2.25]{GwilliamPavlov} and }]
The filtered $\infty$-category $\Fil(\scrC)$ is closed symmetric
monoidal $\infty$-category under the \emph{completed convolution
product} $\wih{\star}$, defined by
\[F \wstar G := \comp(\oblv(F) \star \oblv(G)),\]
for any $F, G \in \Fil(\scrC)$.
\edefn

Moreover, the associated graded functor is strong monoidal
by \cite[Proposition 2.26]{GwilliamPavlov}, intertwining the
completed convolution product on $\Fil(\scrC)$ with the Day
convolution product on $\prod_{n \in \calJ} \scrC$, which
we denote by $\star_{\gr}$.

\subsection{Modules over a filtered ring}
As an application of the formalism of Gwilliam--Pavlov,
we consider the case of filtered modules over a filtered ring.
We first state the definition of a filtered algebra in the higher
categorical context.

\bdef
Let $\calO^{\otimes}$ be an $\infty$-operad
in the sense of \cite{HA},
and let $\scrC^{\otimes}$ be a symmetric monoidal
presentable stable $\infty$-category. A 
\emph{filtered $\calO$-algebra of $\scrC$} 
is a map of $\infty$-operads, $\calO^{\otimes} \to \Fil(\scrC)^{\otimes}$,
where the $\infty$-operad $\Fil(\scrC)^{\otimes}$ is given by Day
convolution.
\edefn

\brem
If $\calO^{\otimes}$ is the associative operad $\Assoc$
(i.e. $\bbE_1^{\otimes}$), we will say that the map 
$\calO^{\otimes} \to \Fil(\scrC)^{\otimes}$
is a filtered associative algebra in $\scrC$.
\erem


\begin{interlude}
At this point, we recall the ``left module" $\infty$-operad
$\LM^{\otimes}$ introduced by Lurie in \S4.2.1 of \cite{HA}.
In what follows, we try to use the notation of that section.
One of the basic features of $\LM^{\otimes}$ is that is
contains $\Assoc^{\otimes}$ as a subcategory and also 
admits a map of $\infty$-operads $\LM^{\otimes} \to \Assoc^{\otimes}$.
Let $\scrC$ denote a monoidal $\infty$-category, given by a fibration
of $\infty$-operads $\scrC^{\otimes} \to \Assoc^{\otimes}$.
The category $\LMod(\scrC)$ is defined to be the category
$\Alg_{\LM/\Assoc}(\scrC)$.
An object $M \in \LMod(\scrC)$ is given by
a functor $M: \LM^{\otimes} \to \scrC^{\otimes}$,
and the restriction $M|_{\Assoc^{\otimes}}$ is an associative
algebra of $\scrC$. There is a categorical fibration $\LMod(\scrC) \to
\Alg(\scrC)$, and the category of left modules over a given algebra $A \in
\Alg(\scrC)$ is defined to be the fiber product $\LMod(\scrC) \times_{\Alg(\scrC)} \{A\}$,
which we denote by $\LMod_A(\scrC)$. We call $\LMod_A(\scrC)$ 
\emph{the category of left $A$-modules in $\scrC$}.
\end{interlude}

Suppose that $R^{\bullet}$ is a classical filtered associative $\C$-algebra
in the sense that $R^i \cdot R^j \subset R^{i+j}$. Then
$R$ can be considered as a filtered associative algebra in the category
of spectra, $\Fil(\Sp)$. As such, we may consider
its category of modules, $\Mod_R(\Fil(\Sp))$. On the other
hand, since $R^{\bullet}$ is classical, we may consider module
in $\Fil(\Ab)$ over it.

\bprop
\label{prop: filtered R-modules into filtered R-module spectra}
There is an embedding,
\[\Mod_R(\Fil(\Ab)) \hook \Mod_R(\Fil(\Sp)),\]
which preserves filtered colimits and sends exact
sequences in $\Mod_R(\Fil(\Ab))$ to fiber sequences
in $\Mod_R(\Fil(\Sp))$.
\eprop

\bproof
We identify $\Ab$ with the heart of the canonical 
t-structure on $\Sp$ and let $i: \Ab \to \Sp$ denote its inclusion.
We note that this inclusion is lax symmetric monoidal.
Composition with $i$ induces a functor on the category
of sequences, $\Seq(\Ab) \to \Seq(\Sp)$.
Restriction of this functor to the subcategory $\Fil(\Ab) \subset \Seq(\Ab)$
followed by composition with the localization functor $\Seq(\Sp) \to
\Fil(\Sp)$ obtains a functor,
\[\Fil(\Ab) \to \Fil(\Sp).\]
Since $i$ is lax symmetric monoidal, so is the
functor $\Fil(\Ab) \to \Fil(\Sp)$ with respect to
the symmetric monoidal structures on each induced by
Day convolution.
As noted above in the special case of $R$, if
$\Assoc^{\otimes} \to \Fil(\Ab)^{\otimes}$ is a filtered associative algebra in $\Ab$,
we obtain a filtered associative algebra in $\Sp$ by composition
with the above functor. In this case, we obtain a commutative
diagram of $\infty$-operads,
\[\begin{tikzcd}
	& {\Fil(\Sp)^{\otimes}} \\
	{\Assoc^{\otimes}} & {\Fil(\Ab)^{\otimes}}
	\arrow[from=2-1, to=2-2]
	\arrow[from=2-2, to=1-2]
	\arrow[from=2-1, to=1-2].
\end{tikzcd}\]
Now, by the general formalism of $\infty$-operads,
we obtain a functor of left module categories,
\[\Mod_R(\Fil(\Ab)) \to \Mod_R(\Fil(\Sp)),\]
given heuristically by taking a left $R$-module,
$\LM^{\otimes} \to \Fil(\Ab))^{\otimes}$ and composing it with
the map of $\infty$-operads $\Fil(\Ab)^{\otimes} \to \Fil(\Sp)^{\otimes}$. 
It remains to show that the functor $\Mod_R(\Fil(\Ab)) \to \Mod_R(\Fil(\Sp))$
preserves limits and take exact sequences to fiber sequences. This follows
immediately from the corresponding statement for $i$. 
\eproof

\brem
\label{rem: ev_k on modules}
Given an integer $k \in \bbZ$, the evaluation functor $\ev_k$
lifts to a functor $\ev_k: \Mod_R(\Fil(\Sp)) \to \Mod_{R^0}(\Sp)$,
which we denote by the same. 
\erem

\brem
\label{rem: gr on R-modules}
Similarly, the associated graded functor $\gr^{\bullet}: \Fil(\Sp)
\to \prod_{\bbZ} \Sp$ lifts to a functor $\gr^{\bullet}: \Mod_R(\Fil(\Sp)) 
\to \prod_{\bbZ} \Mod_{R^0}$, which we denote by the same. 
\erem

\subsection{Sheaves valued in the filtered category}

We end this appendix with a useful proposition, which may be known,
but is not present in \cite{GwilliamPavlov}.

\bprop
\label{prop: Seq commutes with Shv}
Let $\calX$ be an $\infty$-topos. Given a stable $\infty$-category 
$\scrC$ that admits sequential limits,
there is an equivalence of categories,
\[\Shv(\calX; \Seq(\scrC)) \xrightarrow{\simeq} \Seq(\Shv(\calX; \scrC)).\]
\eprop

\bproof
Recall that $\Shv(\calX; \Seq(\scrC))$ is the full subcategory of
$\Fun(\calX^{\op}, \Seq(\scrC))$ that preserves limits.
Since $\Seq(\scrC) := \Fun(\calJ, \scrC)$, we obtain by adjunction
an equivalence 
\[\Fun(\calX^{\op}, \Seq(\scrC)) \xrightarrow[\simeq]{\Ad} \Fun(\calJ, \Fun(\calX^{\op}, \scrC))\]
Note that, under $\Ad$, an object $F \in \Fun(\calX^{\op}, \Seq(\scrC))$
corresponds to a sequence in $\Fun(\calX^{\op}, \scrC)$ whose
$n$th piece is the functor, $\ev_n \circ F$. 
We claim that $\Ad$ sends
limit-preserving functors in $\Fun(\calX^{\op}, \Seq(\scrC))$ to sequences 
in $\Fun(\calX^{\op}, \scrC)$ whose pieces are limit-preserving functors.

Suppose we are given a small diagram $D: I \to \calX^{\op}$, and suppose that
$F \in \Fun(\calX^{\op}, \Seq(\scrC))$ preserves the limit over $D$.
Composition with $D$ gives a diagram of sequences 
$F \circ D \in \Fun(I, \Seq(\scrC))$. 
Since limits in functor categories are computed
pointwise, $\lim_I (F \circ D)$ is 
equivalent to the functor $\calJ \to \scrC$ given
by $n \mapsto \lim_I ({\ev_n \circ F} \circ D)$. On the other hand, because $F$ is
assumed to preserve the limit over $D$, $\lim_I (F \circ D) \simeq F(\lim_I D)$.
Altogether, we have
\begin{align*}
\Ad(F)(\lim D)_n 		&\simeq (\ev_n \circ F)(\lim D) \\
					&\simeq \ev_n(F(\lim D)) \\
					&\simeq \ev_n(\lim (F \circ D)) \\
					&\simeq \lim (\ev_n \circ F \circ D) \\
					&\simeq \lim (\Ad(F)_n \circ D).
\end{align*}
Since $D$ was an arbitrary small diagram, the lemma follows.
\eproof

The following corollary follows easily
from the proof of the above proposition
by restricting to the subcategory of limit-preserving
functors $\calX^{\op} \to \Seq(\scrC)$ which factor
through the fully faithful embedding $i: \Fil(\scrC) \hook \Seq(\scrC)$.

\bcor
\label{cor: Fil commutes with Shv}
Let $\calX$ be an $\infty$-topos.
Given a stable $\infty$-category $\scrC$ that admits sequential limits,
the equivalence established in \cref{prop: Seq commutes with Shv}
restricts to a fully faithful functor,
\[\imath: \Shv(\calX; \Fil(\scrC)) \xrightarrow{} \Fil(\Shv(\calX; \scrC)).\]
\ecor

\bproof
Note that, since $i: \Fil(\scrC) \hook \Seq(\scrC)$ preserves limits,
composition with $i$ induces a functor $\Shv(\calX; \Fil(\scrC)) \to
\Shv(\calX; \Seq(\scrC))$. It is clear that this functor is fully faithful,
and by \cref{prop: Seq commutes with Shv}, we obtain
a fully faithful functor embedding $\Shv(\calX; \Fil(\scrC))
\to \Seq(\Shv(\calX; \scrC))$. It remains to show that this
functor sends equivalences in $\Shv(\calX; \Fil(\scrC))$ to
graded equivalences in $\Seq(\Shv(\calX; \scrC))$.
But this follows from the corresponding statement
for the embedding $i: \Shv(X; \Fil(\scrC)) \to
\Shv(X; \Seq(\scrC))$ and the following commutative
diagram of functor categories,
\[\begin{tikzcd}
	{\Fun(\calX^{\op}, \Seq(\scrC))} & {\Fun(\calJ, \Fun(\calX^{\op}, \scrC))} \\
	{\Fun(\calX^{\op}; \prod_{\bbZ} \scrC)} & {\prod_{\bbZ} \Fun(\calX^{\op};\scrC)}
	\arrow["{\gr^{\bullet} \circ -}", from=1-1, to=2-1]
	\arrow["{\gr^{\bullet}}", from=1-2, to=2-2]
	\arrow["\Ad", from=1-1, to=1-2]
	\arrow["{\simeq }", from=2-1, to=2-2].
\end{tikzcd}\]
\eproof

\printbibliography

\end{document}